\documentclass{ejm}

\usepackage{amsmath}
\usepackage{afterpage}
\usepackage{amssymb}
\usepackage{graphicx}
\usepackage{amsfonts}
\usepackage{url}
\usepackage{epstopdf}
\usepackage{color}
\usepackage{tikz-cd}
\usepackage{bm}
\usepackage[colorlinks=true,linkcolor=blue,citecolor=red]{hyperref}%
\usepackage{multirow}
\usepackage{rotating}
\usepackage{cite}

\def\e{{\bf e}}
\def\G{{\bf G}}
\def\E{{\bf E}}
\def\A{{\bf A}}
\def\I{{\bf I}}
\def\Q{{\bf Q}}
\def\q{{\bf q}}
\def\a{{\bf a}}
\def\b{{\bf b}}
\def\B{{\bf B}}
\def\K{{\bf K}}
\def\D{{\bf D}}
\def\M{{\bf M}}
\def\C{{\bf C}}

\def\Cmu{{\mathcal C}}

\def\nnu{{\boldsymbol\nu}}
\def\xxi{{\boldsymbol\xi}}
\def\eeta{{\boldsymbol\eta}}

\def\det{\mathrm{det}}
\def\trace{\mathrm{trace}}
\def\rank{\mathrm{rank}}

\def\ve{\varepsilon}

\def\pa{\partial\Omega}
\def\H{{\mathbb H}}

\def\R{{\mathbb R}}

\newtheorem{remark}{Remark}

\def\AA{{\bf A}}
\def\MM{{\bf M}}

\def\x{\bm{x}}
\def\y{\bm{y}}

\def\atanh{\mathrm{atanh}}

\title{Competition of small targets in planar domains: 
from Dirichlet to Robin and Steklov boundary condition}

\author[Denis~S.~Grebenkov, Michael J. Ward]{%
  D. \ns S. \ns G\ls R\ls E\ls B\ls E\ls N\ls K\ls O\ls V, \ns    M.\ns J. \ns W\ls A\ls R\ls D}
\affiliation{
Denis S. Grebenkov; 
CNRS -- Universit\'e de Montr\'eal CRM -- CNRS,
6128 succ Centre-Ville, Montr\'eal QC H3C 3J7, Canada; \\
Laboratoire de Physique de la Mati\`{e}re Condens\'{e}e (UMR 7643), \\ 
CNRS -- Ecole Polytechnique, Institut Polytechnique de Paris, 91120 Palaiseau, France \\
(corresponding author, email: denis.grebenkov@polytechnique.edu)

Michael J. Ward; Department of Mathematics, University of British Columbia, Vancouver,
BC, V6T 1Z2, Canada \\ (email: ward@math.ubc.ca)}

\date{\today}

\begin{document}

\label{firstpage}
\maketitle

\baselineskip=12pt

\begin{abstract}
We consider steady-state diffusion in a bounded planar domain with
multiple small targets on a smooth boundary.  Using the method of
matched asymptotic expansions, we investigate the competition of these
targets for a diffusing particle and the crucial role of surface
reactions on the targets.  We start from the classical problem of
splitting probabilities for perfectly reactive targets with Dirichlet
boundary condition and improve some earlier results.  We discuss how
this approach can be generalized to partially reactive targets
characterized by a Robin boundary condition.  In particular, we show
how partial reactivity reduces the effective size of the target.  In
addition, we consider more intricate surface reactions modeled by
mixed Steklov-Neumann or Steklov-Neumann-Dirichlet problems.  We
provide the first derivation of the asymptotic behavior of the
eigenvalues and eigenfunctions for these spectral problems in
the small-target limit. Finally, we show how our asymptotic approach
can be extended to interior targets in the bulk and to exterior
problems where diffusion occurs in an unbounded planar domain outside
a compact set.  Direct applications of these results to
diffusion-controlled reactions are discussed.
\end{abstract}


\baselineskip=16pt 

\setcounter{equation}{0}
\setcounter{section}{0}

\section{Introduction}
\label{sec:intro}

Diffusive search for hidden targets is critically important for
various physical, chemical and biological systems
\cite{Redner,Schuss,Metzler,Masoliver,Lindenberg,Dagdug}.  In the most
basic setting, a point-like particle (e.g., a molecule, an ion, a
protein, a virus, a bacterium, etc.) undergoes diffusive motion inside
a confining environment and searches for an immobile target (e.g., a
catalytic site on a solid surface, a channel on a plasma membrane, a
specific site on the DNA, a cell, etc.).  If the target is hidden in
the bulk, it is often called an interior trap or a sink, whereas a
target on the boundary is referred to as a reactive patch or an escape
window.  In both cases, if the target is small, one usually speaks
about the narrow escape problem \cite{Holcman13,Holcman14}, bearing in
mind the picture of an open window, through which the particle can
leave the domain and never return.  Most former works were dedicated
to finding and even optimizing the {\it mean} first-passage time (FPT)
to a single target or to a given arrangement of multiple targets
\cite{Singer06a,Schuss07,Pillay2010,Cheviakov10,Chen11,Grebenkov16,Lindsay17,Grebenkov17a,Iyaniwura21,Guerin23}.
Other relevant characteristics of the diffusive search such as the
whole distribution of the FPT
\cite{Benichou08,Godec16,Grebenkov18a,Grebenkov19a,Cherry22} and
Laplacian eigenvalues \cite{Kolokolnikov05,Coombs09,Cheviakov11}, were
also studied.

A common limitation of most former works is their emphasis either on a
single target, or on multiple targets of the same type.  In turn, many
biochemical applications involve targets of different types.  For
instance, signal transduction between neurons relies on diffusive
search by calcium ions of a sensor protein on the vesicle with
neurotransmitters inside the presynaptic bouton
\cite{Sala90,Neher08,Holcman13,Guerrier18,Reva21}.  While the sensor
protein is the primary target, calcium ions can reversibly bind to
buffer molecules inside the confining domain or leave it through
calcium channels on its boundary.  Both buffer molecules and channels
play the role of auxiliary targets that compete for calcium ions and
thus allow to control the signal transduction.  More generally, the
successful reaction of a diffusing particle on a ``primary'' target
may fail due to its eventual capture by other targets, or its escape.

When all targets are perfect (i.e., the reaction occurs instantly upon
the first arrival), the competition between targets for a diffusing
particle is characterized via diffusive fluxes, splitting
probabilities and conditional first-passage times
\cite{Traytak96,Traytak97,Felici03,Chevalier11,Berezhkovskii12,Delgado15,Kurella15,Grebenkov19f,Bressloff20,Grebenkov20f}.
In particular, the asymptotic behavior of these quantities for small
interior traps or absorbing patches on the boundary and the dependence
on their spatial arrangement have been studied in depth.
However, as the targets are not perfectly reactive in most
applications
\cite{Collins49,Sano79,Sapoval94,Erban07,Lawley15,Galanti16b,Grebenkov19b,Grebenkov20f,Piazza22,Bressloff22,Grebenkov23b},
their competition also depends on their reactivities.  The role of
partially reactive traps, as modeled by a Robin condition condition,
is not nearly as well understood, especially in the two-dimensional
case.

The problem becomes even more challenging for more intricate
surface reactions, which cannot be described by the conventional Robin
boundary condition on targets. We will refer to such targets as
{\em imperfect}.  For instance, the target reactivity can be
progressively increased or decreased by encounters with a diffusing
particle.  Such activation or passivation processes are described
within the encounter-based approach
\cite{Grebenkov20,Grebenkov20c,Grebenkov23a,Bressloff23d,Bressloff23e}.
In probabilistic terms, the reaction event occurs when the
number of reaction attempts upon each arrival onto the target exceeds
some random threshold.  The probability distribution of the threshold
characterizes the reaction mechanism (see details in
\cite{Grebenkov20}).  For instance, the particular case of the
exponential distribution corresponds to a partially reactive target
with a constant reactivity, and its probabilistic description is
equivalent to solving the diffusion equation with the Robin boundary
condition.  In turn, other distributions of the threshold describe
more intricate surface reactions and involve integral-type boundary
conditions.  As shown in \cite{Grebenkov20}, such PDE problems can be
solved by employing spectral expansions based on the {\em Steklov
problem} (see Sec. \ref{sec:SteklovN} and \ref{sec:SteklovND}
for its formulation and basic properties).  In particular, the
Steklov eigenfunctions turn out to be particularly suitable for
dealing with diffusive motion in the confining domain between
successive arrivals onto an imperfect target.
The peculiar feature of the Steklov problem that distinguishes
it from common spectral problems for the Laplacian, is that the
spectral parameter appears in the boundary condition.  Various
properties of the Steklov problem have been thoroughly investigated
(see \cite{Levitin,Behrndt15,Hassell17,Girouard17,Colbois24} and
references therein).  When imperfect targets are located on the inert
impenetrable boundary, one needs to combine Steklov and Neumann
boundary conditions.  Such a mixed Steklov-Neumann problem was already
known in hydrodynamics, where it is referred to as the sloshing
problem \cite{Henrici70,Fox83,Kozlov04,Levitin22}.  In the case of a
single target, the asymptotic behavior of its eigenvalues and
eigenfunctions in the small-target limit was recently studied
\cite{Grebenkov25}.  However, the scaling arguments and related
analysis from \cite{Grebenkov25} are not directly applicable to the
case of multiple targets.  The asymptotic behavior of the spectrum of
the mixed Steklov-Neumann problem is thus unknown, despite the
importance of its potential applications.  Yet another unstudied
setting concerns a single imperfect target with Steklov condition in
the presence of multiple escape windows with Dirichlet condition.  A
mathematical framework for studying such an escape problem relies on
the mixed Steklov-Neumann-Dirichlet problem \cite{Grebenkov23}.  To
our knowledge, the asymptotic behavior of its eigenvalues and
eigenfunctions in the small-target limit has not been studied
previously.

In this paper, we progressively fill the gap between perfect and
imperfect targets.  In Sec. \ref{sec:Dirichlet}, we start with the
conventional setting of $N$ absorbing sinks and study their splitting
probabilities, i.e., the probability of hitting one sink before any
other.  This relatively simple setting allows us to introduce in a
didactic way many notions and tools that will be employed throughout
the manuscript.  Even though this problem was studied in the past (see
\cite{Chevalier11,Bressloff20} and references therein), we succeed in
improving and generalizing some earlier results.  Section
\ref{sec:Robin} presents an extension to partially reactive targets,
in which the Dirichlet boundary condition is replaced by a Robin
condition.  We show how partial reactivity effectively reduces the
target size.  The major contributions of the paper are presented in
Secs. \ref{sec:SteklovN} and \ref{sec:SteklovND}.  In
Sec. \ref{sec:SteklovN}, we consider the mixed Steklov-Neumann problem
for $N$ imperfect targets.  For this novel problem, we obtain the
asymptotic behavior of its eigenvalues and eigenfunctions in the
small-target limit.  In turn, Sec. \ref{sec:SteklovND} focuses on the
mixed Steklov-Neumann-Dirichlet problem, in which one target is
imperfect (with Steklov condition), whereas the other targets are
perfect (with Dirichlet condition).  We apply matched asymptotic
expansion techniques to investigate the asymptotic behavior in the
small-target limit.  For all considered cases, the accuracy of the
derived asymptotic formulas is illustrated on two examples: the case
of two patches in an arbitrary domain and the case of $N$
equally-spaced patches on the boundary of a disk.  Our analytical
results are compared with numerical solutions obtained by a
finite-element method in Matlab (its home-made implementation for
Steklov problems is described in \cite{Chaigneau24}).  In
Sec. \ref{sec:extensions}, we discuss two further extensions of the
present analysis: the case of interior targets (or traps), and
exterior problems for which diffusion occurs outside a compact set.
In this way, we cover a broad variety of settings, in which multiple
small targets of different types compete for diffusing particles in
planar domains.  We summarize our main results in Section
\ref{sec:conclusion}.

\section{Splitting probabilities on Dirichlet patches}
\label{sec:Dirichlet}

To introduce the theoretical framework and tools, we begin by
revisiting the classical problem of splitting probabilities, which are
commonly used to characterize competition between multiple perfectly
reactive targets for a diffusing particle.  Although this
problem has been studied previously (see
\cite{Chevalier11,Bressloff20} and references therein), we will
improve and generalize some earlier results.

\begin{figure}
\begin{center}
\includegraphics[width=60mm]{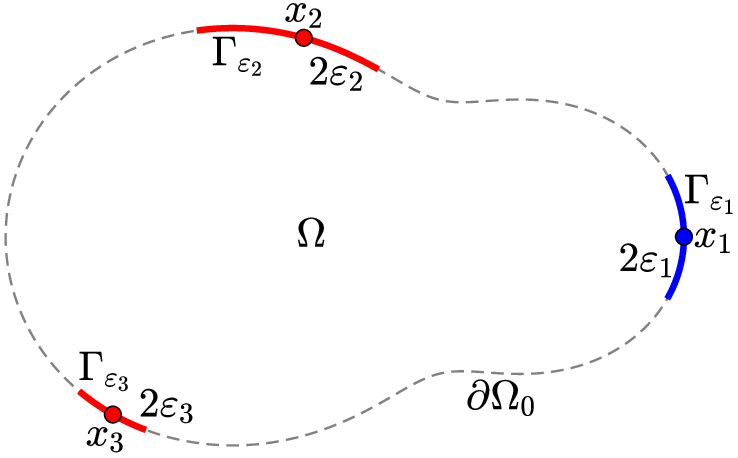} 
\end{center}
\caption{
Illustration of a bounded domain $\Omega \subset\R^2$ with a smooth
boundary $\pa$ split into three absorbing patches $\Gamma_{\ve_i}$ of
length $2\ve_i$ (in red and blue), and the remaining reflecting part
$\pa_0$ (gray dashed line).  For a particle starting from a point
$\x\in\Omega$, the splitting probability $S_1(\x)$ is the probability
of hitting the blue patch $\Gamma_{\ve_1}$ first.}
\label{fig:scheme}
\end{figure}

Let $\Omega \subset \R^2$ be a bounded planar domain with a smooth
boundary $\pa$.  Let $\{\Gamma_{\ve_1}, \ldots, \Gamma_{\ve_N}\}$ be
$N$ disjoint subsets of the boundary $\pa$ that represent multiple
patches of lengths $2\ve_1, \ldots, 2\ve_N$ that are centered at
boundary points $\x_1, \ldots, \x_N$ (each patch $\Gamma_{\ve_j}$ is
connected).  The remaining part of the boundary, denoted as
$\pa_0 = \pa \backslash (\Gamma_{\ve_1} \cup \cdots \cup
\Gamma_{\ve_N})$, is reflecting (Fig. \ref{fig:scheme}).  We are
interested in the small-target limit when all patches are small and
comparable (i.e., $\ve_1 \sim o(1)$ and $\ve_j/\ve_1 \sim {\mathcal
O}(1)$).  We assume that the patches are well-separated in the sense
that $|\x_i - \x_j| = {\mathcal O}(1)$ for all $i\ne j$.

In this section, we consider that all patches $\Gamma_{\ve_j}$
are absorbing sinks (i.e., perfectly reactive targets).  For a
particle started from a point $\x \in \Omega$, we aim at determining
the splitting probability $S_k(\x)$ ($k = 1,\ldots,N$), i.e., the
probability of the arrival onto the Dirichlet patch $\Gamma_{\ve_k}$
{\it before} hitting any other patch.  This probability satisfies the
boundary value problem (BVP)
\begin{subequations}\label{sk:full_no}
\begin{align}
\Delta S_k & = 0 \quad \textrm{in}~\Omega\,, \\  \label{eq:Sj_1}
  S_k &= \delta_{j,k} \quad \textrm{on}~ \Gamma_{\ve_j}\,, \quad
        j\in\lbrace{1,\ldots,N\rbrace}\,, \\
\partial_n S_k & = 0 \quad \textrm{on}~ \pa_0\,,
\end{align}
\end{subequations}
where $\Delta$ is the Laplacian, $\partial_n$ is the normal derivative
oriented outward to the domain $\Omega$, and $\delta_{j,k}$ is the
Kronecker symbol.  In the analysis below,
$k\in\lbrace{1,\ldots,N\rbrace}$ is fixed. In the small-target
  limit $\ve_j\to 0$ for each $j\in \lbrace{1,\ldots,N\rbrace}$, we
  will use the method of matched asymptotic expansions for problems
  with logarithmic interactions \cite{Ward93b} to approximate
  solutions to (\ref{sk:full_no}) that are accurate to all powers of
  ${1/\ln(\ve_j)}$.

\subsection{Inner solutions}

The inner solution near each Dirichlet patch $\Gamma_{\ve_j}$ can be
found by introducing the local coordinates $\y = \ve_j^{-1} \Q_j
(\x-\x_j)$, where $\Q_j$ is an appropriate rotation matrix to restrict
$\y = (y_1,y_2)$ to the upper half-plane $\H_2 = \R \times \R_+$ (the
matrix $\Q_j$ plays no role since $\Q_j^\dagger \Q_j = \I$ and $|\y| =
\ve_j^{-1} |\x-\x_j|$).  We look for an inner solution near the
patch $\Gamma_{\ve_j}$ in the form
\begin{equation}  \label{eq:V2}
  S_k(\x_j + \ve_j \Q_j^\dagger \y) = \delta_{j,k} + A_j g_\infty(\y)\,,  \quad
  j\in\lbrace{1,\ldots,N\rbrace}\,,
\end{equation}
where $A_j$ is an unknown constant, and $g_\infty(\y)$ is the Green's
function satisfying the canonical BVP given by
\begin{subequations}  \label{eq:w_problem}
\begin{align} \label{eq:wc_Laplace}
\Delta g_\infty & = 0  \quad \textrm{in}~ \H_2\,,   \\
\partial_n g_\infty & = 0  \quad \textrm{on}~ y_2 = 0\,, ~ |y_1| \geq 1\,; \qquad
g_\infty  = 0 \quad \textrm{on}~ y_2 =0\,, ~ |y_1| < 1 \,,\\
g_\infty & \sim \ln |\y| + {\mathcal O}(1) \quad \textrm{as}~ |\y| \to \infty\,,
\end{align}
\end{subequations}
(the subscript $\infty$ highlights infinite reactivity of the
perfect patch, see below).  The exact solution of this classical
problem is given in Appendix \ref{sec:wc} for completeness.  The
analysis below will require only the knowledge of the asymptotic
behavior of $g_\infty(\y)$ at infinity.  We recall that the constant
term in this behavior,
\begin{equation}  \label{eq:wc_infinity}
g_\infty(\y) \sim \ln |\y| - \ln(d) + o(1) \quad \textrm{as}~ |\y| \to \infty\,,
\end{equation}
is determined by the logarithmic capacity $d$ of the interval
$(-1,1)$, which is simply $d={1/2}$.
Now putting $\y = \ve_j^{-1} \Q_j (\x-\x_j)$, we get that the
far-field behavior of the inner solution is
\begin{equation}
  S_k \sim \delta_{j,k} + A_j \biggl[\ln |\x-\x_j| - \ln(\ve_j d) +
  o(1)\biggr] \quad \textrm{as}~ \x\to \x_j\,.
\end{equation}
Setting
\begin{equation}
  \nu_j = - \frac{1}{\ln(\ve_j d)} \,, \quad \mbox{with} \quad
    d=\frac{1}{2} \,,
\end{equation}
we rewrite this far-field behavior, for each
$j \in \lbrace{1,\ldots,N\rbrace}$, as
\begin{equation}  \label{eq:Sk_farfield}
  S_k \sim \delta_{j,k} + A_j \biggl[\ln |\x-\x_j| + 1/\nu_j + o(1)\biggr] \quad
  \textrm{as}\,\,\x\to \x_j\,.
\end{equation}

\subsection{Outer solution and matching conditions}
\label{sec:outer_solution}

Now we consider the outer problem for $S_k(\x)$ given by
\begin{subequations}  \label{eq:V_outer}
\begin{align}
  \Delta S_k & = 0  \quad \textrm{in}~ \Omega\,,   \\
  S_k & \sim \delta_{j,k} + A_j \bigl[\ln|\x-\x_j| + 1/\nu_j \bigr]  \quad
  \textrm{as}~ \x \to \x_j\in \pa\,, \quad j\in\lbrace{1,\ldots,N\rbrace}\,, \\
 \partial_n S_k & = 0 \quad \textrm{on}~ \pa \backslash \{\x_1,\ldots,\x_N\}\,. 
\end{align}
\end{subequations}
To find this solution, we introduce the surface Neumann Green's
function $G(\x,\xxi)$, which satisfies
\begin{subequations}  \label{eq:surfNeumann}
\begin{align}
  \Delta G & = \frac{1}{|\Omega|}  \quad \textrm{in}~ \Omega\,; \quad
 \partial_n G  = 0 \quad \textrm{on}~ \pa \backslash\{\xxi\}\,;
   \quad \int\limits_{\Omega} G(\x,\xxi)\,d\x  = 0\,,            
  \label{eq:surfNeumann_asympt}\\
  G(\x,\xxi) & \sim - \frac{1}{\pi} \ln|\x-\xxi| + R(\xxi) + o(1) \quad
  \textrm{as}~ \x \to \xxi\in \pa\,, 
\end{align}
\end{subequations}
where $|\Omega|$ is the area of $\Omega$, and $R(\xxi)$ is the regular
part of $G(\x,\xxi)$, defined by
\begin{equation}
  R(\xxi) = \lim\limits_{\x\to \xxi} \biggl(G(\x,\xxi) +
  \frac{1}{\pi} \ln |\x-\xxi|\biggr)\,.
\end{equation}

\begin{remark} For a disk, $G$ and $R$ are known analytically from
  \cite{Kolokolnikov05} and \cite{Pillay2010} (see
  (\ref{eq:Green_disk}) below). For a square domain, they can be
  represented in terms of rapidly converging infinite series
  representations (see Sec. 3.3 of \cite{Pillay2010}). Similar
  representations can be obtained for rectangles and ellipses from the
  results for the interior Neumann Green's function in Sec. 4.2 of
  \cite{kolok_split} and in Sec. 5 of \cite{Iyaniwura21},
  respectively, by allowing the interior source point to tend to the
  domain boundary (see Sec.~\ref{sec:NGreen_ellipse_int} for this
  derivation for an ellipse and Table \ref{tab:NGreen}). A numerical
  method to compute $G$ and $R$ in arbitrary planar domains is
  described in Sec.~3.3 of \cite{Pillay2010}.
\end{remark}
  
The divergence theorem applied to Eqs. (\ref{eq:V_outer}) yields
\begin{equation}  \label{eq:Asum}
\sum\limits_{j=1}^N A_j = 0\,.
\end{equation}
Under this condition, the solution to Eqs. (\ref{eq:V_outer}) can be
written as the linear combination
\begin{equation}  \label{eq:solution}
S_k(\x) = \chi_k - \sum\limits_{i=1}^N \pi A_i G(\x,\x_i) \,,
\end{equation}
where $\chi_k$ is a constant to be found.  Matching the inner
and outer asymptotic expansions as $\x\to\x_j$ gives
\begin{equation}  \label{eq:Aj}  
A_j + \nu_j\biggl[A_j R_j + \sum\limits_{i=1 \atop i\ne j}^N G_{j,i} A_i\biggr] 
= \chi_k \nu_j - \nu_k \delta_{j,k}\,,
\end{equation}
for $j\in \lbrace{1,\ldots,N\rbrace}$, where we have defined
\begin{equation} \label{eq:Gij}  
G_{j,i} = \pi G(\x_j,\x_i) \quad (i\ne j)\,, \qquad R_j = \pi R(\x_j)\,.  
\end{equation}
Together with the compatibility condition (\ref{eq:Asum}), we obtain a
system of $N+1$ linear equations that determine the unknown
coefficients $A_1, \ldots,A_N$ and $\chi_k$.

\subsection{Matrix reformulation and general solution}
\label{sec:matrix}

To proceed, we introduce the following vectors and matrices of sizes
$N\times 1$ and $N\times N$:
\begin{equation}  \label{eq:matrices_def}
\e  = \left(\begin{array}{c} 1 \\ 1 \\ \cdots \\ 1  \\ \end{array}\right), \quad
\nnu = \left(\begin{array}{c c c c} \nu_1 & 0 & \cdots & 0 \\ 
0 & \nu_2 & \cdots & 0 \\ 
\cdots&\cdots&\cdots&\cdots  \\ 
0 & 0 &\cdots& \nu_N  \\ \end{array}\right)\,,  \quad
\G  = \left(\begin{array}{c c c c} R_1 & G_{1,2} & \cdots & G_{1,N} \\ 
G_{2,1} & R_2 & \cdots & G_{2,N} \\ 
\cdots&\cdots&\cdots&\cdots  \\ 
G_{N,1} & G_{N,2} &\cdots& R_N  \\ \end{array}\right)\, \quad
\e_k = \left(\begin{array}{c} 0 \\ 1 \\ \cdots \\ 0  \\ \end{array}\right)\,,
\end{equation}
where $1$ stands on the $k$-th row of the vector $\e_k$.  

In terms of this notation, the system (\ref{eq:Aj}) is written in
matrix form as
\begin{equation}\label{eq:Avec_eq}
\A + \nnu \G \A = \chi_k \nnu \e - \nnu \e_k \,.
\end{equation}
Applying $\e^\dagger$ on the left, we isolate $\chi_k$ as
\begin{equation}  \label{eq:chiS}
\chi_k = \frac{\e^\dagger \nnu \G \A}{\bar{\nu}} + \frac{\nu_k}{\bar{\nu}}  \,,
\end{equation}
where we used $\e^\dagger \A = 0$ due to Eq. (\ref{eq:Asum}), and
defined
\begin{equation}  \label{eq:barnu}
\bar{\nu} = \e^\dagger \nnu \e = \sum\limits_{j=1}^N \nu_j \,.
\end{equation}
Eliminating $\chi_k$ from Eq. (\ref{eq:Avec_eq}), we get a
matrix equation for $\A$ given by
\begin{equation}\label{eq:Amat_full}
  (\I + \nnu \G) \A - \nnu \e \, \frac{\e^\dagger \nnu \G \A}{\bar{\nu}}
  = \frac{\nu_k}{\bar{\nu}} \nnu \e - \nu_k \e_k \,.
\end{equation}
By introducing the matrices $\M_0$ and $\E$ by
\begin{equation}  \label{eq:M0}
\M_0  = \I + \biggl(\I - \frac{\nnu \E}{\bar{\nu}}\biggr) \nnu \G\,, \qquad
\E  = \e \e^\dagger \,,
\end{equation}
the solution to Eq. (\ref{eq:Amat_full}) is
\begin{equation} \label{eq:A_splitting}
\A = \frac{\nu_k}{\bar{\nu}}  \M_0^{-1} \bigl(\nnu \e - \bar{\nu} \e_k \bigr)\,,
\end{equation}
where $\M_0^{-1}$ is the inverse of $\M_0$.  In fact, in the
small-target limit, all $\nu_j \ll 1$ so that $\M_0$ is a small
perturbation of the identity matrix $\I$ and is thus invertible.  Once
the coefficients $A_i$ are found, the constant $\chi_k$ follows from
Eq. (\ref{eq:chiS}).  As a consequence, the splitting probability
$S_k(\x)$ is fully determined via the representation
(\ref{eq:solution}).  This is the main result of this section.  The
shape of the confining domain $\Omega$ and the arrangement of
Dirichlet patches are captured by the matrix $\G$, whereas the sizes
of patches are accounted for via the matrix $\nnu$.  When the surface
Neumann Green's function is known analytically, a numerical
computation of the coefficients $A_i$ and $\chi_k$ is fast, at least
if the number of targets is not too large. We emphasize that the
  solution $\A$ in (\ref{eq:A_splitting}) to the linear system
  (\ref{eq:Amat_full}) has accounted for all logarithmic correction
  terms in the asymptotic expansion of the splitting probability. This
  technique for effectively summing what otherwise would be an
  infinite logarithmic expansion in powers of $\nu_k$ was developed in
  \cite{Ward93b}, and has been used in other contexts (see
  \cite{Coombs09}, \cite{Kurella15}, \cite{Pillay2010}).

We further emphasize that Eq. (\ref{eq:solution}) is only applicable in the
outer region, i.e., when $|\x-\x_j| \gg {\mathcal O}(\ve_j)$ for all
$j \in \lbrace{1,\ldots,N\rbrace}$.  In turn, if the starting point
$\x$ is too close to $\x_j$, this asymptotic formula may give wrong
values (e.g., negative or exceeding $1$).  In practice, the outer
solution can be capped by $0$ and $1$ to avoid such invalid values,
i.e., one can use $\max\{0, \min\{1, S_k(\x)\}\}$ instead of
$S_k(\x)$.  We remark that if an accurate approximation of the
splitting probability is needed near the patch, one has to use the
corresponding inner solution.

We also note that the constant $\chi_k$ can be interpreted as
the volume-averaged splitting probability.  In fact, if the starting
point $\x$ is not fixed but uniformly distributed in $\Omega$, the
average over the starting point yields
\begin{equation}\label{eq:Sk_bar}
  \overline{S}_k = \frac{1}{|\Omega|} \int\limits_{\Omega} S_k(\x) \,
  d\x = \chi_k \,,
\end{equation}
where we used $\int_{\Omega} G(\x,\x_i)\, d\x=0$.

\subsection{Example of two patches}

In the case of two targets ($N = 2$), the matrix $\M_0$ from
Eq. (\ref{eq:M0}) reads
\begin{equation}  
\M_0 = \I + \gamma   
\left(\begin{array}{c c} R_1 - G_{1,2} & G_{1,2} - R_2 \\ G_{1,2} - R_1 & R_2 - G_{1,2} \\  \end{array}\right), 
\end{equation}
where $\gamma = \nu_1 \nu_2/(\nu_1 + \nu_2)$.  The inverse of
this matrix is
\begin{equation}
\M_0^{-1}  = \frac{1}{1 + \gamma(R_1 + R_2 - 2G_{1,2})}  \left( \I + 
\gamma \left(\begin{array}{c c} R_2 - G_{1,2} & R_2 -G_{1,2} \\ R_1 -G_{1,2}
     & R_1 - G_{1,2} \\  \end{array}\right)  \right) \,. \label{eq:M0_inv_two}
\end{equation}
Substituting this expression into Eq. (\ref{eq:A_splitting}), we find
for $k = 1$ that
\begin{equation}  \label{eq:disk_split_A2} 
  -A_1 = A_2 = \biggl(\frac{1}{\nu_1} + \frac{1}{\nu_2} +
  \bigl[R_1 + R_2 - 2G_{1,2}\bigr]\biggr)^{-1} \,.
\end{equation}
Then, by using Eq. (\ref{eq:chiS}) we determine $\chi_1$ as
\begin{equation}   
  \chi_1 = \frac{\nu_1 - A_2 \nu_1 (R_1 - G_{1,2}) +
    A_2\nu_2 (R_2 - G_{1,2})}{\nu_1+\nu_2} \,,
\end{equation}
which can be further simplified as
\begin{equation}  \label{eq:disk_split_chi} 
  \chi_1 = \frac{1/\nu_2 + (R_2 - G_{1,2})}{1/\nu_1 +
    1/\nu_2 + (R_1 + R_2 - 2G_{1,2})} \,.
\end{equation}

For instance, if $\Omega$ is the unit disk, the surface Neumann
Green's function is well known \cite{Pillay2010}:
\begin{equation}  \label{eq:Green_disk}
G(\x,\xxi) = - \frac{1}{\pi} \ln |\x-\xxi| + \frac{|\x|^2}{4\pi} - \frac{1}{8\pi} \,,
\quad R(\xxi) = \frac{1}{8\pi} \,.
\end{equation}
Substitution of these expressions into Eqs. (\ref{eq:disk_split_A2},
\ref{eq:disk_split_chi}) yields
\begin{subequations}  \label{eq:A2chi_two}
\begin{align}
  -A_1 = A_2 & = \biggl(-\ln(\ve_1\ve_2) + 2\ln(2) +
               2\ln|\x_1-\x_2|\biggr)^{-1}\,,   \\
  \chi_1 & = \frac{-\ln(\ve_2/2) + \ln|\x_1-\x_2|}
          {-\ln(\ve_1\ve_2/4) + 2\ln|\x_1-\x_2|}  \,,
\end{align}
\end{subequations}
and we conclude that
\begin{equation}  \label{eq:Sk_two}
S_1(\x) = \chi_1 + A_2 \ln \biggl(\frac{|\x-\x_2|}{|\x-\x_1|} \biggr)\,.
\end{equation}
We can easily check that $S_1(\x)$ approaches $0$ [resp., $1$] as
$\ve_1\to 0$ [resp., $\ve_2 \to 0$], as expected.  
%

In the special case of two identical targets, $\nu_1 = \nu_2 = \nu$,
one has $\chi_1 = 1/2$ and $A_2 = 1/(2/\nu + 2\ln|\x_1-\x_2|)$ so that
\begin{equation} \label{eq:Sk_twoequal}
  S_1(\x) = \frac12 \biggl[1 + \frac{\nu}{1 + \nu \ln |\x_1-\x_2|} \ln
  \biggl(\frac{|\x-\x_2|}{|\x-\x_1|}\biggr)\biggr]\,.
\end{equation}
We remark that if we were to expand the denominator of the second term
into a Taylor series in powers of $\nu \ll 1$ up to
${\mathcal O}(\nu^2)$, we would recover the truncated approximation
given in Eq. (98) from \cite{Chevalier11}.  However, as
$\nu = -1/\ln(\ve/2)$ is not necessarily small enough, our new result
Eq. (\ref{eq:Sk_twoequal}) that incorporates all logarithmic terms is
preferable than using the previous truncated approximation from
\cite{Chevalier11}.

\begin{figure}
\begin{center}
\includegraphics[width=88mm]{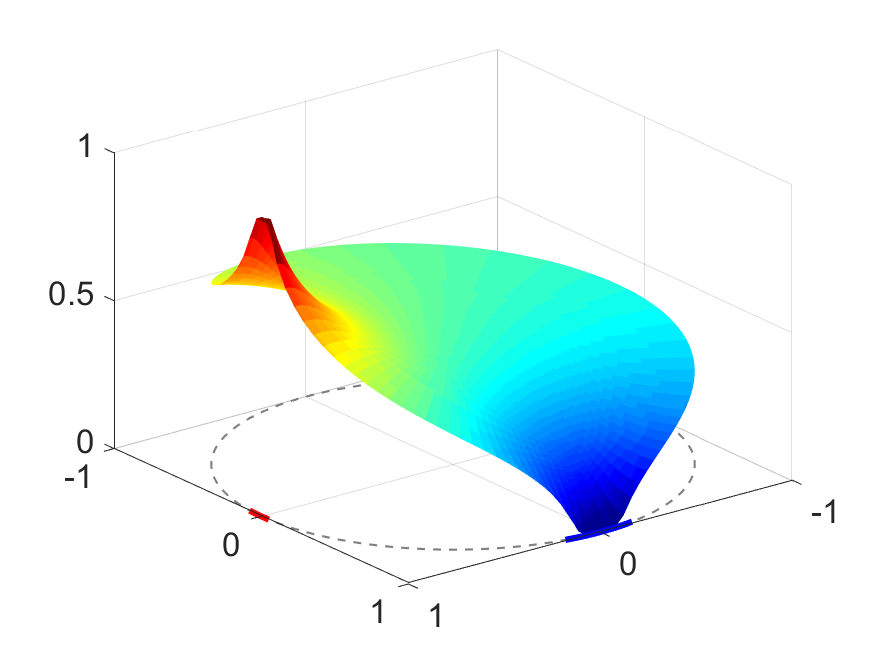} 
\end{center}
\caption{
Splitting probability $S_1(\x)$, given by Eq. (\ref{eq:Sk_two}), for
the unit disk with two Dirichlet patches of length $2\ve_1 = 0.2$
(red) and $2\ve_2 = 0.4$ (blue).  Note that $S_1(\x)$ was capped by
$0$ and $1$, i.e., we plotted $\max\{0, \min\{1, S_1(\x)\}\}$.}
\label{fig:disk_splitting}
\end{figure}

Figure \ref{fig:disk_splitting} illustrates the splitting probability
$S_1(\x)$ from Eq. (\ref{eq:Sk_two}).  As explained earlier, we plot
the capped version of this quantity, $\max\{0, \min\{1, S_1(\x)\}\}$,
to avoid invalid values near two patches.  Expectedly, $S_1(\x)$
increases as $\x$ gets closer to the first patch (red arc) and
decreases as $\x$ gets closer to the second patch (blue arc).
However, as the outer solution (\ref{eq:Sk_two}) is not applicable in
the vicinity of these patches, we can observe some discrepancy, e.g.,
$S_1(\x)$ does not vanish on the second patch, as it should.  To amend
this discrepancy, we can use the inner solution when $|\x-\x_j|
\lesssim \ve_j$.

\subsection{Example of equally-spaced identical patches on the boundary of the unit disk}
\label{sec:Dirichlet_N}

If there are $N$ identical targets, one has $\nu_j = \nu$ so that
$\nnu = \nu \I$, $\bar{\nu} = N \nu$, and thus $\M_0 = \I + \nu (\I -
\E/N) \G$ from Eq. (\ref{eq:M0}).  Moreover, if the patches are
equally-spaced on the boundary of the unit circle, the matrix $\G$
defined in Eq. (\ref{eq:matrices_def}) is circulant and symmetric.  As
a consequence, its eigenvectors can be written as
\begin{equation}  \label{eq:qj_def}
\q_j = \frac{1}{\sqrt{N}} \bigl(\omega^{-j}, \omega^{-2j}, \ldots , \omega^{-Nj}\bigr)^\dagger\,,  \qquad j\in \lbrace{1,\ldots,N\rbrace}\,,
\end{equation}
where $\omega = e^{2\pi i/N}$, and the transposition $\dagger$ now
denotes the Hermitian conjugate.  Upon taking the real and imaginary
parts of $\q_j$, the resulting real-valued eigenvectors form an
orthonormal basis in $\R^{N}$ since $\G$ is symmetric.  Let us denote
by $\kappa_j$ the associated eigenvalues of $\G$:
\begin{equation} \label{eq:kappa_def} 
  \kappa_j = R_1 + \sum\limits_{m=1}^{N-1} \omega^{mj} \, G_{1,1+m}\,,
  \qquad j\in \lbrace{1,\ldots,N\rbrace}\,.
\end{equation}
For any $j \in \lbrace{1,\ldots,N-1\rbrace}$, we get
\begin{equation}
  \M_0 \q_j = \bigl(\I + \nu [\I - \E/N]\bigr) \kappa_j \q_j =
  (1 + \nu \kappa_j) \q_j \,,
\end{equation}
because $\E \q_j = \e \e^\dagger \q_j = 0$ for any $1 \leq j \leq N-1$
due to orthogonality of $\q_j$ to $\q_N = \e/\sqrt{N}$.  As a
consequence, each $\q_j$ with $j\in \lbrace{1,\ldots,N-1\rbrace}$ is
also the eigenvector of $\M_0$, associated to the eigenvalue $1 + \nu
\kappa_j$.  In addition, we have
\begin{equation}
\M_0 \q_N = \bigl(\I + \nu [\I - \E/N]\bigr) \kappa_N \q_N = \q_N \,,
\end{equation}
since $\E \q_N /N = \q_N$. Therefore, $\q_N$ is the eigenvector of
$\M_0$ associated with the eigenvalue $1$.  We use this spectral
information to invert the matrix $\M_0$ as
\begin{equation}  \label{eq:M0_inv_N}
  \M_0^{-1} = \q_N \q_N^\dagger + \sum\limits_{j=1}^{N-1} \q_j
  (1 + \nu \kappa_j)^{-1} \q_j^\dagger \,.
\end{equation}
Substituting this spectral representation into
Eq. (\ref{eq:A_splitting}), we get
\begin{equation} \label{eq:A_splitting_N}
\A = - \nu \sum\limits_{j=1}^{N-1} \q_j (1 + \nu \kappa_j)^{-1} \q_j^\dagger \e_k\,,
\end{equation}
where we used the orthogonality of the eigenvectors $\q_j$.
Substitution of this expression into Eq. (\ref{eq:chiS}) yields 
$\chi_k = 1/N$.  This is consistent with the interpretation of
$\chi_k$ as the volume-averaged splitting probability: when all
patches are identical and equally-spaced on the boundary of the unit
disk, they are equivalent from the uniformly distributed starting
point, so that $\overline{S}_k = 1/N$ from Eq. (\ref{eq:Sk_bar}).

To complete this example, we will simplify Eq. (\ref{eq:kappa_def}) by
using the explicit form (\ref{eq:Green_disk}) of the surface Neumann
Green's function for the unit disk.  Since the centers of the patches
$\x_j$ are equally-spaced on the domain boundary, we have $\x_j =
e^{2\pi i(j-1)/N} = \omega^{j-1}$ for
$j\in\lbrace{1,\ldots,N\rbrace}$.  In this way, substituting
Eq. (\ref{eq:Green_disk}) into Eq. (\ref{eq:kappa_def}), we get
\begin{equation}\label{eq:kappa_orig}
  \kappa_j = \frac{1}{8}\sum_{m=0}^{N-1}\omega^{jm}
  - \sum\limits_{m=1}^{N-1} \omega^{mj} \ln\vert 1 - \omega^m\vert \,,
  \qquad j\in \lbrace{1,\ldots,N\rbrace} \,,
\end{equation}
where we interpret points as complex numbers and $|z|$ as the
modulus of $z$.  For $j=N$, for which $\omega^{Nm}=1$, we get
\begin{equation}\label{eq:kappaN}
  \kappa_N = \frac{N}{8} - \ln\left\vert \prod_{m=1}^{N-1}(1-\omega^m)\right\vert
  = \frac{N}{8}-\ln{N} \,,
\end{equation}
where the product in Eq. (\ref{eq:kappaN}) was evaluated by using the
roots of unity together with L'Hopital's rule to get $\lim_{z\to 1}
{(z^N-1)/(z-1)}=N=\prod_{m=1}^{N-1}(1-\omega^m)$.  For $j<N$, the
first term in Eq. (\ref{eq:kappa_orig}) vanishes and we obtain
\begin{equation}  \label{eq:kappaj}
  \kappa_j = -  \sum\limits_{m=1}^{N-1} \omega^{mj} \ln|1 - \omega^m|\,, \qquad
  j\in \lbrace{1,\ldots,N-1\rbrace} \,,
\end{equation}
which reduces after some simplifications to
\begin{equation}\label{eq:kappaj_new}
\kappa_j = \ln{2} - \sum_{m=1}^{N-1}\cos\left(\frac{2\pi j m}{N}\right)
  \ln\left[ \sin\left(\frac{\pi m}{N}\right) \right]\,,
\end{equation}
for $j\in \lbrace{1,\ldots,N-1\rbrace}$.  The asymptotic behavior of
$\kappa_j$ for large $N$ is derived in Appendix \ref{sec:kappa}.

\section{Splitting probability on Robin patches}
\label{sec:Robin}

In most applications, targets are not perfectly reactive
\cite{Collins49,Sano79,Sapoval94,Erban07,Lawley15,Galanti16b,Grebenkov19b,Grebenkov20f,Piazza22,Bressloff22,Grebenkov23b}.
Starting from Collins and Kimball \cite{Collins49}, partial reactivity
is usually implemented by replacing a Dirichlet boundary condition by
a Robin condition.  In the case of splitting probabilities, a
straightforward generalization of the previous setting consists in
replacing Dirichlet boundary condition (\ref{eq:Sj_1}) by the Robin
boundary condition:
\begin{equation}  \label{eq:Sk_RobinBC}
  \partial_n S_k + q_j S_k = q_j \delta_{j,k} \quad \textrm{on}~ \Gamma_{\ve_j}\,,
  \qquad j\in \lbrace{1,\ldots,N\rbrace}\,,
\end{equation}
where the constant $0 < q_j < \infty$ characterizes the reactivity of
the $j$-th patch $\Gamma_{\ve_j}$.  Here we excluded the limit $q_j =
0$ that would correspond to an inert patch that could be treated
as a part of the reflecting boundary $\pa_0$.  The Dirichlet
condition is recovered in the limit $q_j \to +\infty$.

The change of the boundary condition on the patch is a local effect
that does not impact the outer solution.  In turn, the inner solution
near each patch $\Gamma_{\ve_j}$ in Eq. (\ref{eq:V2}) should now be
replaced by
\begin{equation}  \label{eq:V2b}
  S_k(\x_j + \ve_j \Q_j^\dagger \y) = \delta_{j,k} + A_j g_{\ve_j q_j}(\y)\,,  \qquad
  j\in \lbrace{1,\ldots,N\rbrace},
\end{equation}
where $g_{\mu}(\y)$ is the Robin Green's function, which satisfies
\begin{subequations}  \label{eq:Steklov_inner}
\begin{align}
\Delta g_{\mu} & = 0  \quad \textrm{in}~ \H_2\,,   \\
  \partial_n g_{\mu} & = 0  \quad \textrm{on}~ y_2 = 0\,, ~ |y_1| \geq 1\,;\qquad
\partial_n g_{\mu} + \mu g_{\mu}  = 0 \quad \textrm{on}~ y_2 =0\,, ~ |y_1| < 1\,,
                        \label{eq:Steklov_inner_BC}\\
g_{\mu} & \sim \ln |\y| + {\mathcal O}(1) \quad \textrm{as}~ |\y| \to \infty\,.
\end{align}
\end{subequations}
In Appendix \ref{sec:Cmu_asympt}, we derive a spectral expansion for
this Green's function: 
\begin{equation}  \label{eq:gR}
  g_{\mu}(\y) = g_\infty(\y) + \pi \sum\limits_{k=0}^\infty
  \frac{\Psi_{2k}(\infty)}{\mu_{2k} + \mu} \Psi_{2k}(\y)\,,
\end{equation}
for any $\mu \notin \bigcup_{k=0}^\infty \{- \mu_{2k}\}$.  Here
$\mu_k$ and $\Psi_k(\y)$ are the eigenvalues and eigenfunctions of the
auxiliary Steklov-Neumann problem in the upper half-plane:
\begin{subequations} \label{eq:Vkint}
\begin{align}   \label{eq:Vkint_Laplace}
\Delta \Psi_k & = 0 \quad \textrm{in}~ \H_2 \,,\\
  \partial_n \Psi_k & = \mu_k \Psi_k \quad \textrm{on}~ y_2 = 0\,, ~ |y_1| < 1\,;
\qquad \partial_n \Psi_k  = 0 \quad \textrm{on}~ y_2 = 0\,, ~ |y_1| \geq 1\,, \\
\Psi_k & = {\mathcal O}(1) \quad \textrm{as}~ |\y|\to \infty\,.
\end{align}
\end{subequations}
An efficient numerical procedure for computing these eigenmodes is
summarized in Appendix \ref{sec:Steklov_interval}.  Once tabulated,
these eigenfunctions play a role of ``special functions'', like
orthogonal polynomials.

\begin{figure}
\begin{center}
\includegraphics[width=88mm]{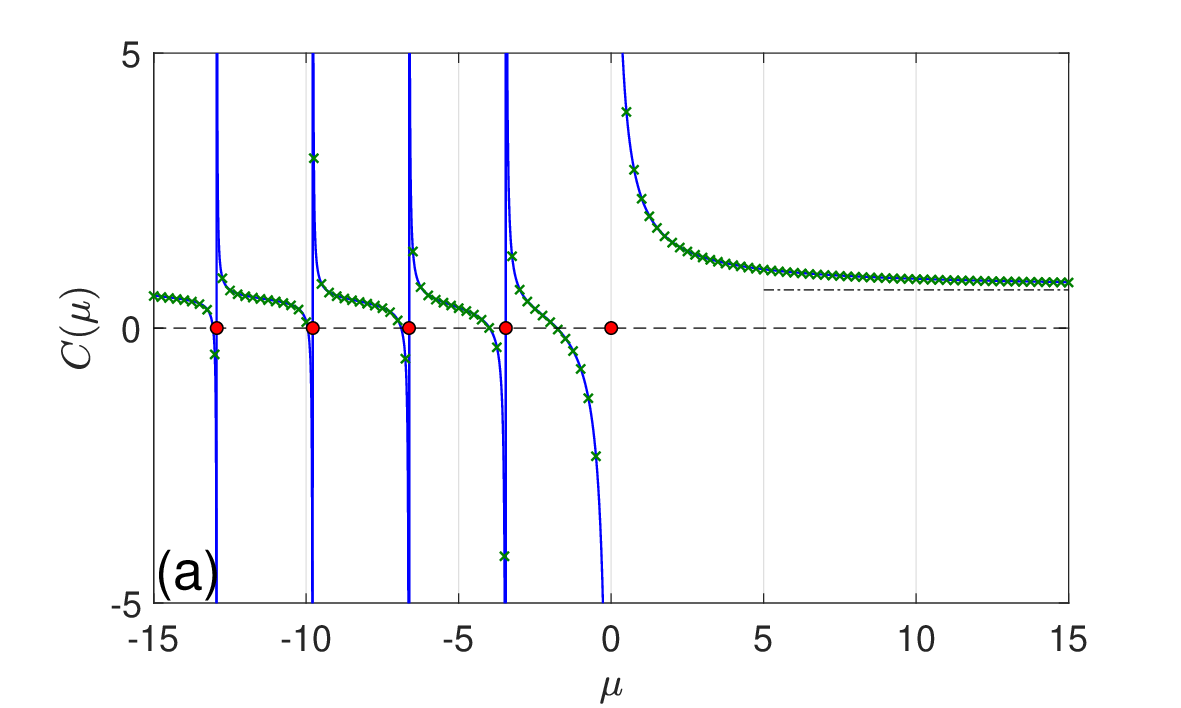} 
\includegraphics[width=88mm]{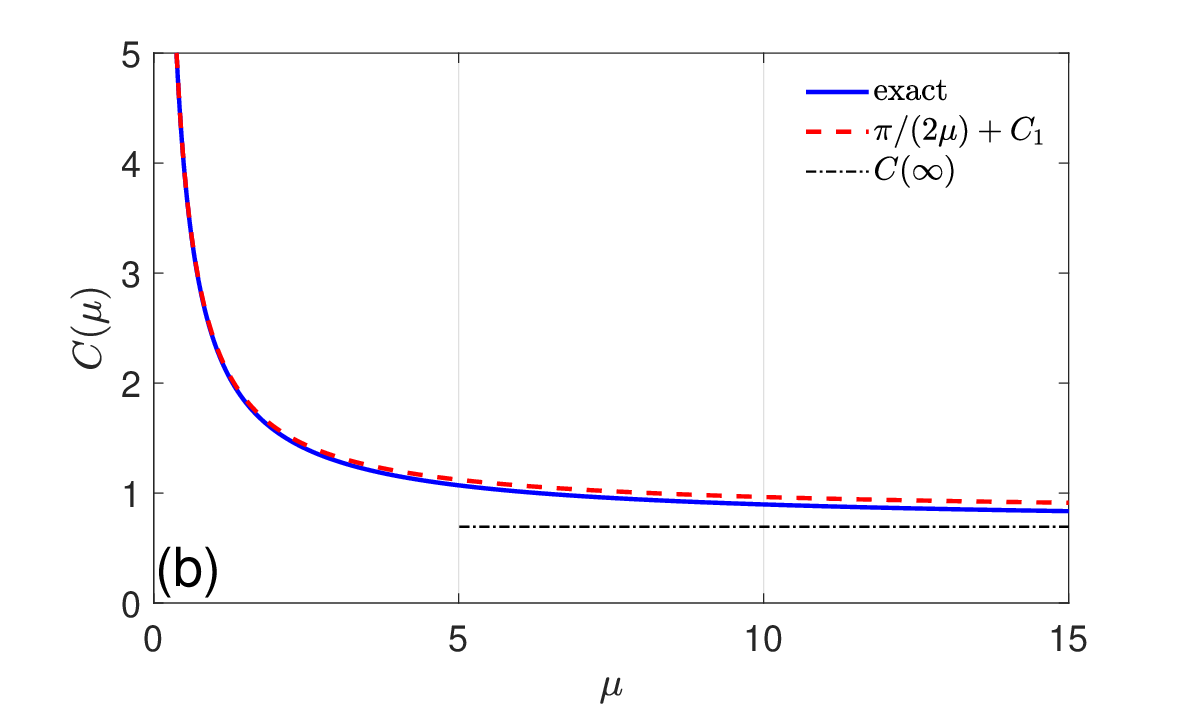} 
\end{center}
\caption{ 
{\bf (a)} Function $\Cmu(\mu)$ from Eq. (\ref{eq:Cmu_spectral}), 
in which the infinite series is truncated either to 50 terms (solid
line) or to 10 terms (crosses), to highlight the accuracy of both
truncations. 
Filled circles indicate the values $-\mu_{2k}$, at which $\Cmu(\mu)$
diverges.  Dash-dotted line outlines the asymptotic limit $\ln(2)$ of
$\Cmu(\mu)$ as $\mu\to \infty$.  {\bf (b)} Comparison of
$\Cmu(\mu)$ and its approximation (\ref{eq:Cmu_approx}), which is
accurate over a broad range of $\mu$. }
\label{fig:Cmu}
\end{figure}

As a consequence, Eq. (\ref{eq:gR}) determines the constant term
$\Cmu(\mu)$ in the asymptotic behavior of $g_\mu(\y)$ at infinity,
defined by
\begin{equation}  \label{eq:ws_asympt}
g_\mu(\y) \sim \ln |\y| + \Cmu(\mu) + o(1) \quad \textrm{as}~ |\y| \to \infty\,.
\end{equation}
In the limit $|\y| \to \infty$, we get 
\begin{equation}  \label{eq:Cmu_spectral}
  \Cmu(\mu) = \ln (2) + \frac{\pi}{2\mu} + \pi \sum\limits_{k=1}^\infty
  \frac{[\Psi_{2k}(\infty)]^2}{\mu_{2k} + \mu} \,,
\end{equation}
for any $\mu \notin \bigcup_{k=0}^\infty \{- \mu_{2k}\}$, where
we used $\mu_0 = 0$ and $\Psi_0(\infty) = 1/\sqrt{2}$ (see Appendix
\ref{sec:Cmu_asympt} for more details).  The first ten coefficients
contributing to $\Cmu(\mu)$ are listed in Table
\ref{tab:mu_int}, while Fig. \ref{fig:Cmu} illustrates the behavior of
the function $\Cmu(\mu)$.

\begin{figure}
\begin{center}
\includegraphics[width=88mm]{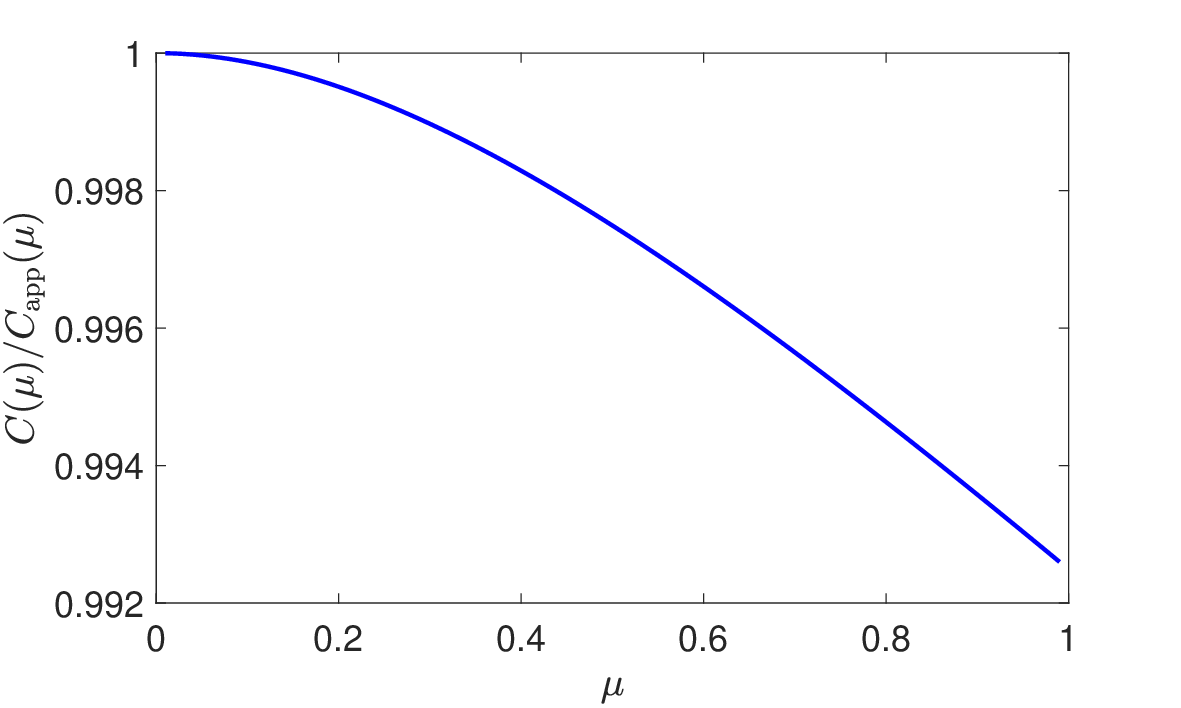} 
\end{center}
\caption{The ratio of $\Cmu(\mu)$ with its approximation (\ref{eq:Cmu_approx})
  is very close to unity on the range $0<\mu<1$.}\label{fig:Cmu_c}
\end{figure}

A Taylor expansion of Eq. (\ref{eq:Cmu_spectral}) near $\mu=0$ yields
\begin{equation}  \label{eq:Cmu_Taylor}
\Cmu(\mu) = \frac{\pi}{2\mu} + \sum\limits_{n=0}^\infty (-\mu)^n \, C_{n+1} \,, 
\end{equation}
where the coefficients $C_n$ can be expressed in terms of $\mu_{2k}$
and $\Psi_{2k}(\infty)$ (see Appendix \ref{sec:Cmu_asympt}).
Moreover, we  calculated in Appendix \ref{sec:Cmu} the exact
values of the first two coefficients as
\begin{equation}\label{eq:C0_C1}
C_1  = \frac32 - \ln (2) \approx 0.8069 \,, \qquad
C_2  = \frac{21 - 2\pi^2}{18 \pi} \approx 0.0223\,.
\end{equation}
Since the second and higher-order coefficients turn out to be small,
the following small-$\mu$ approximation,
\begin{equation} \label{eq:Cmu_approx} \Cmu(\mu) \approx {\mathcal
    C}_{\rm app}(\mu) \equiv \frac{\pi}{2\mu} + C_1 \quad \mbox{for} \quad
  \mu\ll 1 \,,
\end{equation}
is remarkably accurate as seen in both Fig. \ref{fig:Cmu}(b) and
Fig.~\ref{fig:Cmu_c}. This approximation is one of the key results
needed for sections \ref{sec:SteklovN} and \ref{sec:SteklovND} below.

The relation (\ref{eq:ws_asympt}) implies 
\begin{equation}  \label{eq:V_A1}
S_k \sim \delta_{j,k} + A_j \biggl[\ln|\x-\x_j| + 1/\nu_j + o(1) \biggr] \quad \textrm{as}~ \x\to \x_j\,,
\end{equation}
where we now redefine $\nu_j$ as
\begin{equation}  \label{eq:nuj_Robin}
\nu_j = \frac{1}{-\ln(\ve_j) + \Cmu(\ve_j q_j)} \,.
\end{equation}
It follows that the far-field behavior of the inner solution is
identical to that in Eq. (\ref{eq:Sk_farfield}) for Dirichlet patches,
whereas the partial reactivity is fully taken into account through the
new definition (\ref{eq:nuj_Robin}) of $\nu_j$.  As a consequence, we
retrieve the same representation (\ref{eq:solution}) for the splitting
probability $S_k(\x)$, with the coefficients $A_i$ and $\chi_k$
given by Eqs. (\ref{eq:chiS}, \ref{eq:A_splitting}).  This equivalence
shows that partially reactive targets with $q_j > 0$ can still be
treated as the perfect ones but with the reduced effective length,
defined by
\begin{equation}  \label{eq:eps_eff}
\ve_j^{\rm eff} = \ve_j \exp\bigl(\ln (2) - \Cmu(\ve_j q_j)\bigr) \,.
\end{equation}
Together with the spectral expansion (\ref{eq:Cmu_spectral}), this is
the main result of this section.

From Eq. (\ref{eq:Cmu_spectral}), the function $\Cmu(\mu)$ decreases
monotonically from $+\infty$ to $\ln(2)$ on the range $\mu>0$ (see
Fig. \ref{fig:Cmu}).  As a consequence, $\nu_j$ in
Eq. (\ref{eq:nuj_Robin}) decreases monotonically from
${-1/\ln(\ve_j/2)}$ (this is the former definition of $\nu_j$ for the
Dirichlet patch) to $0$, whereas $\ve_j^{\rm eff}$ decreases from
$\ve_j$ to $0$ as the reactivity $q_j$ drops from infinity to $0$.
This shows that a target with a smaller reactivity has less chance to
capture the diffusing particle.
When $q_j \sim {\mathcal O}(1)$, one has $\ve_j q_j \ll 1$, so that
the approximation (\ref{eq:Cmu_approx}) is applicable.  This yields,
that $\nu_j \approx 2\ve_j q_j/\pi \ll 1$ and so to leading order
\begin{equation}
\ve_j^{\rm eff} \approx \ve_j e^{-\pi/(2q_j \ve_j)}  \qquad (\ve_j q_j \ll 1) \,.
\end{equation}

For weakly reactive targets (i.e., if $q_j \ve_j \ll 1$ for all
  $j = 1,\ldots,N$), we get $\nu_j \approx 2\ve_j q_j/\pi \ll 1$.
According to Eq. (\ref{eq:A_splitting}), all the coefficients $A_i$
are small (of the order of $\ve$) so that the first term in
Eq. (\ref{eq:chiS}) can be neglected, yielding
\begin{equation}  \label{eq:chi_approx0}
\chi_k \approx \frac{\nu_k}{\nu_1 + \ldots + \nu_N}  \,,
\end{equation}
and thus
\begin{equation}
\overline{S}_k \approx \frac{\ve_k q_k}{\ve_1 q_1 + \ldots + \ve_N q_N} \,,
\end{equation}
independently of the location of the patches.  We emphasize that the
approximation (\ref{eq:chi_approx0}) is generally not accurate for
perfect targets: even if $\ve_j$ are very small, the gauge function
$\nu_j = -1/\ln(\ve_j/2)$ may not be small enough to neglect
higher-order terms in powers of $\nu_j$.

\begin{figure}
\begin{center}
\includegraphics[width=88mm]{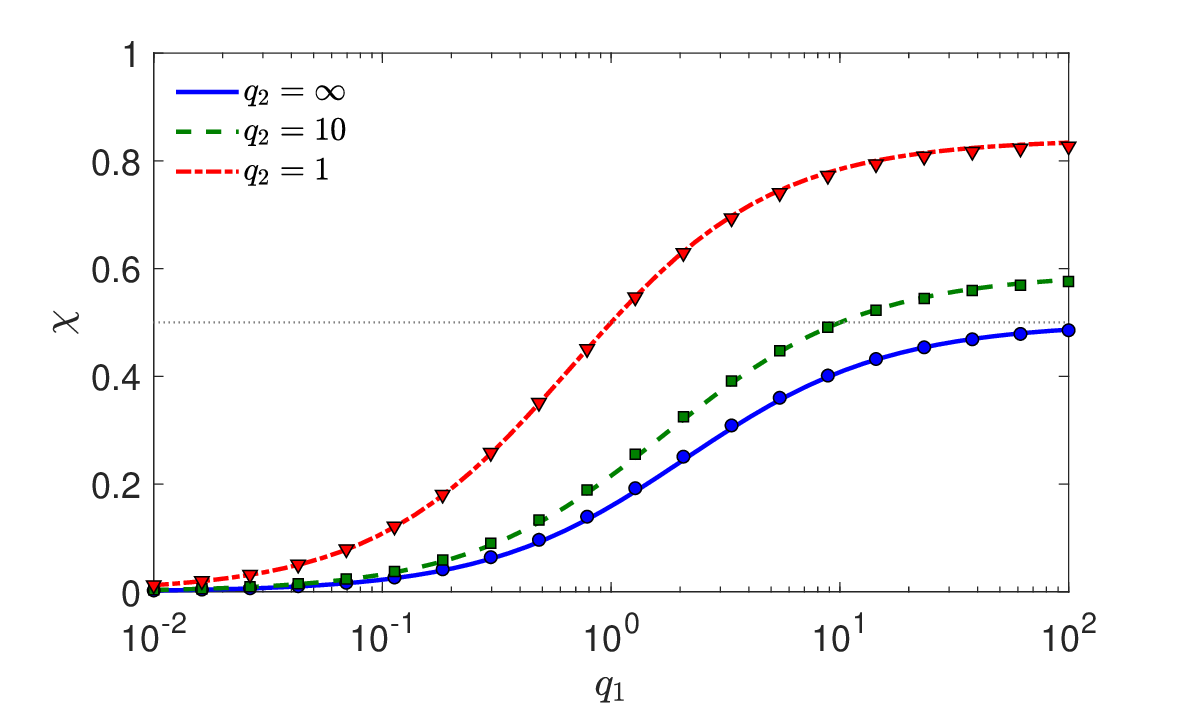} 
\end{center}
\caption{ 
Volume-averaged splitting probability $\overline{S}_1 = \chi_1$ for
the unit disk, calculated from (\ref{eq:A2chi_twoRobin_b}),
with two patches of equal length $2\ve = 0.2$ located at boundary
points $(\pm 1,0)$.  Three curves correspond to three values of the
reactivity parameter $q_2$ of the second patch.  Symbols present
the numerical solution of the BVP (\ref{eq:Sj_1}) with Robin boundary
condition (\ref{eq:Sk_RobinBC}) by a finite-element method in Matlab
PDEtool, with the maximal meshsize $0.02$.}
\label{fig:disk_Robin}
\end{figure}

For the case of two partially reactive patches on the boundary of the
unit disk, substitution of the effective lengths $\ve_j^{\rm eff}$
from Eq. (\ref{eq:eps_eff}) into Eqs. (\ref{eq:A2chi_two}) yields
\begin{subequations}  \label{eq:A2chi_twoRobin}
\begin{align}
A_2 & = \biggl(-\ln(\ve_1\ve_2) + \Cmu(q_1\ve_1) + \Cmu(q_2\ve_2) + 2\ln|\x_1-\x_2|\biggr)^{-1}\,,   \\
\chi_1 & = \frac{-\ln(\ve_2) + \Cmu(q_2\ve_2) + \ln|\x_1-\x_2|}
   {-\ln(\ve_1\ve_2) + \Cmu(q_1\ve_1) + \Cmu(q_2\ve_2)+ 2\ln|\x_1-\x_2|}  \,.
         \label{eq:A2chi_twoRobin_b}
\end{align}
\end{subequations}
%
Figure \ref{fig:disk_Robin} shows the behavior of $\chi_1$ for two
patches of equal length. In this figure, the high accuracy of our asymptotic
solution (\ref{eq:A2chi_twoRobin_b}) is confirmed by comparison with a
numerical solution of the BVP (\ref{eq:Sj_1}) with Robin boundary
condition (\ref{eq:Sk_RobinBC}) by a finite-element method.

\subsection{The mean first-reaction time}\label{sec:mfrt}

Although our asymptotic analysis has focused on calculating
splitting probabilities it can be easily modified to calculate the
mean first-reaction time (MFRT).

The dimensionless MFRT $u(\x)$ satisfies a Poisson equation with mixed
Neumann-Robin boundary conditions:
\begin{subequations} \label{mfpt:ssp}
\begin{align}
  \Delta u & = - 1 \,, \quad \x \in \Omega \,,
                  \label{mfpt:ssp_1}\\
  \partial_{n} u + q_i u & = 0\,, \quad \x \in \Gamma_{\ve_i}
                          \,, \quad i\in \lbrace{1,\ldots,N\rbrace} \,, \\
  \partial_{n} u & = 0 \,, \quad \x \in \partial \Omega_0
                   =\pa \backslash (\Gamma_{\ve_1} \cup \cdots \cup
\Gamma_{\ve_N})\,. \label{1:ssp_2}
\end{align}
\end{subequations}
As previously, each reactive boundary patch $\Gamma_{\ve_i}$ has
length $2\ve_i$, reactivity parameter $q_i$, and is centered at
$\x_i\in \partial\Omega$.

The matched asymptotic analysis of Eq. (\ref{mfpt:ssp}) in the
small-target limit $\ve_i\ll 1$ is very similar to that for analyzing
the splitting probability. The inner solution near the $j$-th patch in
terms of an unknown coefficient $A_j$ is
\begin{equation}\label{mfrt:inn}
  V_{j}(\y)= u(\x_j+\ve_j \Q_j^\dagger \y) = A_j g_{\mu_j}(\y)\,,  \qquad
  j\in \lbrace{1,\ldots,N\rbrace},
\end{equation}
where $\mu_j\equiv \ve_j q_j$ and $g_{\mu}(\y)$ is the Robin Green's
satisfying Eq. (\ref{eq:Steklov_inner}) with far-field behavior
Eq. (\ref{eq:ws_asympt}). Upon matching the far-field behavior of
$V_{j}(\y)$ to the outer solution, we find that to within all
logarithmic terms the outer solution satisfies
\begin{subequations}  \label{mfrt:U_outer}
\begin{align}
  \Delta u & = -1  \quad \textrm{in}~ \Omega\,,   \\
  u & \sim  A_j \bigl[\ln|\x-\x_j| + 1/\nu_j \bigr]  \quad
      \textrm{as}~ \x \to \x_j\in \pa\,, \quad j\in\lbrace{1,\ldots,N\rbrace}\,,
  \label{mfrt:U_outer_2}\\
 \partial_n u & = 0 \quad \textrm{on}~ \pa \backslash \{\x_1,\ldots,\x_N\}\,, 
\end{align}
\end{subequations}
where $\nu_j$ is defined by (\ref{eq:nuj_Robin}).
%
The solvability condition for (\ref{mfrt:U_outer}) is that
$\sum_{j=1}^{N} A_j={|\Omega|/\pi}$. We then represent $u$ in terms of
the surface Neumann Green's function and the volume average
$\overline{u}_0=|\Omega|^{-1}\int_{\Omega} u(\x)\, d\x$ as
\begin{equation}  \label{mfrt:solution}
 u(\x) = \overline{u}_0 - \pi \sum\limits_{i=1}^N A_i G(\x,\x_i) \,.
\end{equation}
Imposing the singularity behavior in Eq. (\ref{mfrt:U_outer_2}), we
obtain an $(N+1)$-dimensional linear algebraic system for
$\overline{u}_0$ and $A_1,\ldots,A_N$ given by
\begin{equation}  \label{mfrt:Aj}  
A_j + \nu_j\biggl[A_j R_j + \sum\limits_{i=1 \atop i\ne j}^N G_{j,i} A_i\biggr] 
= \overline{u}_0 \nu_j\,, \quad j\in\lbrace{1,\ldots,N\rbrace}\,; \qquad
\sum_{j=1}^{N} A_j = \frac{|\Omega|}{\pi} \,,
\end{equation}
where $G_{j,i}$ and $R_i$ were defined in (\ref{eq:Gij}).
In matrix form Eq. (\ref{mfrt:Aj}) is written for
$\A\equiv (A_1,\ldots,A_N)^\dagger$ as
\begin{equation}  \label{mfrt:mat}
  \A + \nnu \G \A = \overline{u}_0 \nnu \e \,, \quad \e^{T} \A =
  \frac{|\Omega|}{\pi} \,,
\end{equation}
where $\e$, $\nnu$ and the Green's matrix $\G$ were defined in
Eq. (\ref{eq:matrices_def}). By eliminating $\overline{u}_0$ in
Eq. (\ref{mfrt:mat}), we conclude that
\begin{equation}  \label{mfrt:u0}
  \overline{u}_0 = \frac{|\Omega|}{\pi \bar{\nu}} +
  \frac{\e^\dagger \nnu \G \A}{\bar{\nu}} \,,
\end{equation}
where $\bar{\nu}\equiv \sum\limits_{j=1}^N \nu_j$, while $\A$ is the solution
to
\begin{equation}\label{mfrt:asolve}
  \M_0 \A = \frac{|\Omega|}{\pi \bar{\nu}} \nnu \,,
\end{equation}
with $\M_0$ and $\E$ being defined in Eq. (\ref{eq:M0}).

By using Eq. (\ref{eq:nuj_Robin}) for $\nu_j$, which involves the
local reactivity parameter $q_j$ on the patch, one can invert $\M_0$
in Eq. (\ref{mfrt:asolve}) to get the coefficients $A_j$.  As a
consequence, Eq. (\ref{mfrt:u0}) gives access to the volume-averaged
MFRT $\overline{u}_0$, whereas Eq. (\ref{mfrt:solution}) determines
the MFRT $u(\x)$ for any well-separated spatial configuration of
partially reactive patches. This result generalizes that in
\cite{Pillay2010}, where perfect reactivities ($q_j=\infty$) were
assumed.

\section{Steklov-Neumann problem}
\label{sec:SteklovN}

As discussed in Sec. \ref{sec:intro}, the Robin boundary condition
describes targets with a constant reactivity.  In turn, more
sophisticated surface reactions can be incorporated by using the
encounter-based approach \cite{Grebenkov20,Grebenkov20c,Grebenkov23a},
which relies on the mixed Steklov-Neumann problem.  In this section,
we apply the tools described above to derive the 
asymptotic properties of this spectral problem in the
small-target limit.

As before, we consider a bounded planar domain $\Omega$ with a smooth
boundary $\pa$, which has $N$ small well-separated patches
$\Gamma_{\ve_j}$, and $\pa_0 = \pa \backslash (\Gamma_{\ve_1} \cup
\ldots \cup \Gamma_{\ve_N})$.  We study the mixed Steklov-Neumann
spectral problem:
\begin{subequations} \label{eq:SteklovN_original}
\begin{align}
\Delta V & = 0 \quad \textrm{in}~\Omega \,,\\  \label{eq:SteklovN}
  \partial_n V & = \sigma V \quad \textrm{on}~ \Gamma_{\ve_1} \cup \cdots
                 \cup \Gamma_{\ve_N}  \,, \\
\partial_n V & = 0 \quad \textrm{on} ~ \pa_0 \,.
\end{align}
\end{subequations}
This spectral problem is known to have a discrete positive spectrum
\cite{Levitin}, i.e., infinitely many eigenpairs $\{\sigma_k, V_k\}$
that are enumerated by $k = 0,1,\ldots$ to form an increasing sequence
of eigenvalues: $0 = \sigma_0 < \sigma_1 \leq \sigma_2 \leq \ldots
\nearrow \infty$.  Note that the Steklov boundary condition
(\ref{eq:SteklovN}) with a nonnegative $\sigma$ differs from the
previous Robin condition by the opposite sign.  We aim at determining
the asymptotic behavior of the eigenvalues $\sigma_k$ and the
associated eigenfunctions $V_k$ in the small-target limit.

In the case of a single Steklov patch ($N = 1$), the small-$\ve_1$
asymptotic behavior of the eigenvalues and eigenfunctions was analyzed
in \cite{Grebenkov25}.  In fact, a simple scaling argument suggests
that $\sigma_j \approx \mu_j/\ve_1$ ($j=1,2,\ldots$) to leading
order, where $\mu_j$ are the eigenvalues of the mixed Steklov-Neumann
problem (\ref{eq:Vkint}) for the interval in the upper half-plane
(see also Appendix \ref{sec:Cmu_asympt}).

If there are two well-separated Steklov patches, it is tempting to
apply the same scaling argument in the vicinity of each patch.  In
this way, we can expect that the spectrum of the problem
(\ref{eq:SteklovN_original}) with $N = 2$ is composed of two sequences
of eigenvalues: $\{\mu_j/\ve_1\}$ from the first patch of length
$2\ve_1$, and $\{\mu_j/\ve_2\}$ from the second patch of length
$2\ve_2$.  In other words, the two patches might be expected to not
interact with each other in the small-target limit as $\ve_j \to 0$.
This intuitive argument turns out to be correct for all the
eigenvalues, except for the first nontrivial eigenvalue $\sigma_1$.
Indeed, if two patches could be treated as independent, the eigenvalue
$\sigma_1$ would have to be zero, as $\sigma_0$.  However, the zero
eigenvalue can only correspond to a constant eigenfunction, so that if
$\sigma_1$ was zero, one would have $V_1 = const = V_0$, which is
impossible.  We conclude that even if the patches are extremely small,
the eigenvalue $\sigma_1$ must be strictly positive, and its
asymptotic behavior must result from long-range interactions between
two patches.

In this section, we adapt the analysis from Sec. \ref{sec:Robin} to
the case of $N$ Steklov patches and determine the asymptotic behavior
of the first $N-1$ eigenvalues $\sigma_j$, for
$j \in\lbrace{1,\ldots,N-1\rbrace}$, in this setting.

\subsection{Matched asymptotic analysis}

As before, we look at the inner solution near each Steklov patch
$\Gamma_{\ve_j}$.  Upon comparing the Robin and Steklov conditions
(\ref{eq:Sk_RobinBC}, \ref{eq:SteklovN}), we notice two differences:
(i) $q_j$ is replaced by $-\sigma$, and (ii) there is no inhomogeneous
term $q_j \delta_{j,k}$ in the right-hand side.  Apart from these two
points, the boundary value problems for $S_k(\x)$ and $V(\x)$ are
identical.  As a consequence, we can immediately rewrite the
asymptotic behavior (\ref{eq:V_A1}) for each
$j\in\lbrace{1,\ldots,N\rbrace}$ as
\begin{equation}\label{eq:SN_sing}  
  V \sim A_j \biggl[\ln|\x-\x_j| + 1/\nu_j + \Cmu(-\sigma \ve_j) +
  o(1) \biggr] \quad \textrm{as}~ \x\to \x_j\,,
\end{equation}
where we now redefine $\nu_j$ as
\begin{equation}
\nu_j = - {1/\ln(\ve_j)} \,, 
\end{equation}
and where $A_j$ is an unknown coefficient.  Note that the constant
term $\Cmu(-\sigma \ve_j)$ is not incorporated into the new definition
of $\nu_j$, as we did earlier in the Robin case.

As before, the outer solution is represented as
\begin{equation}  \label{eq:Vsolution}
V = \chi - \pi \sum\limits_{i=1}^N A_i G(\x,\x_i)\,,
\end{equation}
with an unknown constant $\chi$.  The divergence theorem still ensures
the compatibility condition (\ref{eq:Asum}).
Upon enforcing the singularity behavior (\ref{eq:SN_sing}) for the
solution in Eq. (\ref{eq:Vsolution}), we get
\begin{equation}\label{eq:SN_mat}
  \chi - (\G \A)_j = A_j (1/\nu_j + \Cmu(-\sigma \ve_j))\,,  \quad
  j\in \lbrace{1,\ldots,N\rbrace} \,,
\end{equation}
where we used the matrix notations introduced in
Eq. (\ref{eq:matrices_def}) of Sec. \ref{sec:matrix}.  Multiplying
this equation by $\nu_j$ and introducing the diagonal matrix $\C$
formed by $\{\Cmu(-\sigma \ve_1), \ldots, \Cmu(-\sigma \ve_N)\}$, we
rewrite Eq. (\ref{eq:SN_mat}) as
\begin{equation} \label{eq:auxil12}
\bigl(\I + \nnu \G + \nnu \C\bigr) \A = \chi \nnu \e \,.
\end{equation}
Left-multiplying this equation by $\e^\dagger$, and using $\e^\dagger
\A = 0$, we isolate $\chi$ as
\begin{equation}  \label{eq:chi_Steklov_SN}
\chi = \frac{1}{\bar{\nu}} \bigl(\e^\dagger \nnu\G + \e^\dagger \nnu\C\bigr)\A \,,
\end{equation}
where $\bar{\nu}$ was defined by Eq. (\ref{eq:barnu}).
Substituting this expression back into Eq. (\ref{eq:auxil12}), we
obtain that
\begin{equation}  \label{eq:SN_A_matrix}
  \biggl(\I + \biggl(\I - \frac{\nnu \E}{\bar{\nu}}\biggr) \nnu \left(
    \G + \C\right) \biggr) \A = {\bf 0} \,.
\end{equation}
The necessary and sufficient condition for the existence of a
nontrivial solution to this matrix equation is
\begin{equation} \label{eq:detSN_equation}
\det\biggl(\I + \biggl(\I - \frac{\nnu \E}{\bar{\nu}}\biggr) \nnu (\G + \C)\biggr) = {\bf 0} \,.
\end{equation}
The matrices $\nnu$ and $\G$ are determined by the sizes and
arrangement of the Steklov patches, while the matrix $\C$ is formed by
$\{\Cmu(-\sigma \ve_1),\ldots, \Cmu(-\sigma \ve_N)\}$, with the
function $\Cmu(\mu)$ given by Eq. (\ref{eq:Cmu_spectral}).  As a
consequence, Eq. (\ref{eq:detSN_equation}) determines the unknown
parameter $\sigma$.  Moreover, the functional form
(\ref{eq:Cmu_spectral}) implies that there are infinitely many {\it
negative} solutions, denoted as $-\sigma_j$, which are actually
small-target approximations of the Steklov eigenvalues.  

Let us first provide qualitative insights on these solutions.  If the
matrix $\C$ was fixed in the small-target limit, the condition
$\nu_j\ll 1$ would imply the smallness of the second matrix term in
Eq. (\ref{eq:detSN_equation}) as compared to the identity matrix $\I$,
thus ensuring the positivity of the determinant.  To compensate the
smallness of the matrix $\nnu$, the matrix $\C$ must therefore be
large in the small-target limit.  This is possible when at least one
$\sigma \ve_i$ is close to $- \mu_{2k}$ for some $k$.  In other words,
one can expect that a solution $-\sigma_j$ of
Eq. (\ref{eq:detSN_equation}) is close to $\mu_{2k}/\ve_i$ for some
$i$ and $k$.  This intuitive picture suggests that an eigenvalue
$\sigma_j$ of the mixed Steklov-Neumann problem with $N$ patches can
be approximated by that on a single patch (say, $\Gamma_{\ve_i}$), as
if there were no other patches and the associated eigenfunction was
localized on $\Gamma_{\ve_i}$.  The situation is, however, more subtle
in the vicinity of $\mu_0 = 0$.
In the analysis below, we focus on the asymptotic behavior of the
first $N$ eigenvalues $\sigma_j$.

\subsection{First $N$ eigenvalues}

Let us assume that  $- \mu = \sigma \ve_j \ll 1$ for all $j =
1,\ldots,N$, so that we can apply the approximate relation
(\ref{eq:Cmu_approx}) for the function $\Cmu(\mu)$.  Introducing the
diagonal matrix $\eeta$ formed by $\{ {\pi/(2\ve_1)}, \ldots,
{\pi/(2\ve_N)}\}$, we get $\C \approx - \eeta/\sigma + C_1 \I$.  Its
substitution into Eq. (\ref{eq:detSN_equation}) yields
\begin{equation}\label{eq:SN_deteig}
  \det\biggl(\I + \biggl(\I - \frac{\nnu \E}{\bar{\nu}}\biggr) \nnu
  \bigl[\G - \eeta/\sigma + C_1 \I\bigr]\biggr) = {\bf 0} \,.
\end{equation}
Upon defining $\M_1$ and $\B_1$ by
\begin{equation}\label{eq:M1_B1}
  \M_1  = \I + \biggl(\I - \frac{\nnu \E}{\bar{\nu}}\biggr) \nnu
  \bigl(\G + C_1 \I\bigr) \,, \qquad
\B_1  = \biggl(\I - \frac{\nnu \E}{\bar{\nu}}\biggr) \nnu \eeta\,,
\end{equation}
we can rewrite Eq. (\ref{eq:SN_deteig}) in a more compact form as
\begin{equation}
\det(\sigma \M_1 - \B_1) = 0\,.
\end{equation}
In the small-target limit, one has $\nu_j \ll 1$ so that $\M_1$ is
invertible since it is a small perturbation of the identity matrix.
As a consequence, we get that
\begin{equation}\label{eq:Mat_main}
\det(\sigma \I - \M_1^{-1} \B_1) = 0\,.
\end{equation}

We now prove that the eigenvalues $\hat{\sigma}_j$ for
$j=0,\ldots,N-1$ of the matrix $\M_1^{-1} \B_1$ are real.  To do so,
we first write $\M_1$ in Eq. (\ref{eq:M1_B1}) as
\begin{equation}\label{eq:Mat_m1_1}
  \M_1 = \I - \hat{\B}_1 \hat{\G} \,, \,\,
\end{equation}
where we define the symmetric matrices $\hat{\B}_1$ and $\hat{\G}$ by
\begin{equation}\label{eq:Mat_m1_2}
  \hat{\B}_1=\B_1\eeta^{-1} = \left(\I - \frac{\nnu \E}{\bar{\nu}}\right)\nnu
  \,, \qquad \hat{\G}=-\G - C_1 \I \,.
\end{equation}
By using the Neumann series to calculate $\M_1^{-1}$, which converges
since $\nu_j\ll 1$, we obtain that
\begin{equation}
 \M_1^{-1} \hat{\B}_1 = \hat{\B}_1 + \hat{\B}_1\left(\sum_{n=1}^{\infty} \K_n\right)
  \hat{\B}_1 \,, \qquad \K_n = \left( \hat{\G} \hat{\B}_1\right)^{n-1}
    \hat{\G} \,.
\end{equation}
Since $\K_n$ is symmetric for each $n=1,2,\ldots$, it follows that
$\M_1^{-1}\hat{\B}_1$ is symmetric.  Finally, we denote
$\D=\M_1^{-1}\B_1 =\M_{1}^{-1}\hat{\B}_1 \eeta$ and introduce
$\hat{\D} = \eeta^{\frac12} \D \eeta^{-\frac12} =
\eeta^{\frac12} \M_1^{-1} \hat{\B}_1 \eeta^{\frac12}$, which is
symmetric, so that its eigenvalues are real.  Since $\eeta$ is
positive definite, $\hat{\D}$ is related to $\D$ by a similarity
transformation so that the latter must also have real eigenvalues
$\hat{\sigma}_j$ for $j=0,\ldots,N-1$.

Moreover, we note that $\hat{\sigma}_0 = 0$.  This follows since
$\B_1^\dagger \e = 0$, so that $\e$ is an eigenvector of
$\B_1^\dagger$ and thus of $[\M_1^\dagger]^{-1} \B_1^\dagger$.  The
associated eigenvalue $0$ is thus an eigenvalue of $\B_1 \M_1^{-1}$ as
well as of $\M_1^{-1} \B_1$.

In summary, the eigenvalues of the matrix $\M_1^{-1} \B_1$ provide the
leading terms in the asymptotic behavior of the first $N$ eigenvalues
$\sigma_j$ of the Steklov-Neumann problem:
\begin{equation}
\sigma_j \approx \hat{\sigma}_j\,, \qquad  j\in\lbrace{1,\ldots,N-1\rbrace}\,.
\end{equation}
This is the main result of this section.  We observe that $\sigma_0 =
\hat{\sigma}_0 = 0$, as expected.

\subsection{Associated eigenfunctions}

In addition, we construct the associated eigenfunctions of the first
$N$ eigenvalues.  For this purpose, let us rewrite
Eq. (\ref{eq:SN_A_matrix}) as
\begin{equation}
(\sigma \I - \M_1^{-1} \B_1) \A = {\bf 0}\,.
\end{equation} 
While each eigenvalue $\hat{\sigma}_j$ of the matrix $\M_1^{-1} \B_1$
approximates the $j$-th eigenvalue of the mixed Steklov-Neumann
problem, the corresponding eigenvector of this matrix is the vector of
coefficients $A_i$ determining the associated eigenfunction $V_j$ via
Eq. (\ref{eq:Vsolution}), up to a multiplicative factor; we recall
that $\chi$ is given by Eq. (\ref{eq:chi_Steklov_SN}).

The missing multiplicative factor can be fixed by imposing an
appropriate normalization of eigenfunctions.  For the Steklov problem,
the natural normalization is
\begin{equation}
\int\limits_{\Gamma} V^2\, ds = 1 \,,
\end{equation}
where $\Gamma = \Gamma_{\ve_1}\cup \cdots \cup \Gamma_{\ve_N}$.  For
the trivial eigenvalue $\sigma_0 = 0$, one has a constant
eigenfunction $V_0$, whose normalization yields: $V_0^2 = 1/|\Gamma|$.
In the following, we assume that $\sigma > 0$.
 
To proceed, we recall that the inner solution near the $i$-th patch
reads in local coordinates is
\begin{equation}
V(\x_i + \ve_i \Q_i^\dagger \y) \approx A_i g_{-\sigma \ve_i}(\y)\,.
\end{equation}
Using the representation (\ref{eq:gR}) of the Green's function
$g_{\mu}(\y)$, the restriction of $V$ onto $\Gamma_i$ becomes
\begin{equation} \label{eq:SteklovSN_Vi}
V|_{\Gamma_i}(y_1)  = V(\x_i + \ve_i \Q_i^\dagger (y_1,0)^\dagger) 
\approx \pi A_i \sum\limits_{k=0}^\infty \frac{\Psi_{2k}(\infty)
  \Psi_{2k}(y_1,0)}{\mu_{2k} - \sigma \ve_i} \,.
\end{equation}
This equation helps to deduce the required condition on the
coefficients $A_i$:
\begin{align*}
1 & = \sum\limits_{i=1}^N \int\limits_{\Gamma_i} V^2 \, ds
    = \sum\limits_{i=1}^N \ve_i \pi^2 A_i^2 \sum\limits_{k=0}^\infty
    \frac{[\Psi_{2k}(\infty)]^2}{(\mu_{2k} - \sigma \ve_i)^2} \,,
\end{align*}
where we used the orthogonality of the eigenfunctions $\Psi_k$ (see
Appendix \ref{sec:Cmu_asympt}).  The last sum can be re-written as the
derivative of $\Cmu(\mu)$, denoted as $\Cmu^{\prime}(\mu)$:
\begin{equation}
1 = - \pi \sum\limits_{i=1}^N \ve_i A_i^2 \Cmu^{\prime}(-\sigma \ve_i) \,.
\end{equation} 
When $0 < \sigma \ve_i \ll 1$, the Taylor expansion
(\ref{eq:Cmu_Taylor}) implies $\Cmu^{\prime}(\mu) \approx -\pi/(2\mu^2)$ and
thus
\begin{equation}  \label{eq:SteklovSN_Ai_norm}
\frac{2\sigma^2}{\pi^2} \approx  \sum\limits_{i=1}^N \frac{A_i^2}{\ve_i } \,.
\end{equation} 

To complete this section, let us briefly discuss the positivity of
Steklov eigenfunctions on patches $\Gamma_{\ve_j}$.  For a single
patch, all Steklov eigenfunctions $V_j$ must change sign on the patch
due to their orthogonality to $V_0 = 1/\sqrt{|\Gamma|}$.  When there
are $N$ Steklov patches, the orthogonality still holds so that any
eigenfunction $V_j$ with $j > 0$ must change sign on the union of
patches $\Gamma = \Gamma_{\ve_1} \cup \cdots \cup \Gamma_{\ve_N}$.
However, it is generally unknown whether $V_j$ changes the sign or not
on each patch $\Gamma_i$.
Looking at Eq. (\ref{eq:SteklovSN_Vi}), one can expect that if
$\sigma\ve_i$ is small enough, the eigenfunction $V$ does not change
sign on the patch $\Gamma_i$ (i.e., it is either positive, or negative
on it).  Indeed, the term $1/(-2\sigma \ve_i)$ of the sum in
Eq. (\ref{eq:SteklovSN_Vi}) that corresponds to $k = 0$, is expected
to provide the dominant contribution as compared to the remaining
terms.  This property follows from the conjectured inequality
(\ref{eq:gR_negative}).  In other words, if the patches are small
enough, the first $N$ eigenfunctions do not change their signs
on each patch.  This conjecture is confirmed by several numerical
examples (not shown).

\subsection{Example of two patches}

When $N = 2$, we calculate that
\begin{equation}
(\I - \nnu \E/\bar{\nu})\nnu = \gamma \left(\begin{array}{cc} 1 & -1 \\ -1 & 1 \\ \end{array} \right),
\end{equation}
where we label $\gamma = \nu_1\nu_2/(\nu_1+ \nu_2)$.  Then, from
Eq. (\ref{eq:M1_B1}), we get
\begin{equation*} 
\M_1  = \I + \gamma 
\left(\begin{array}{cc} (R_1 - G_{1,2}) + C_1 & (G_{1,2} - R_2)-C_1 \\ 
(G_{1,2}-R_1) - C_1 & (R_2 - G_{1,2}) + C_1 \\ \end{array} \right) \,, \qquad
\B_1  = \frac{\pi}{2} \gamma \left(\begin{array}{cc} 1/\ve_1
         & -1/\ve_2 \\ -1/\ve_1 & 1/\ve_2 \\ \end{array} \right)\,.
\end{equation*}
Since $\e$ is a left-eigenvector of $\M_1$ with eigenvalue one, the
second eigenvalue of $\M_1$ is simply $\trace(\M_1)-1$.  As a
consequence, we find
\begin{equation} 
\det(\M_1) = \trace(\M_1)-1=1 + \gamma[(R_1 + R_2 - 2G_{1,2}) + 2C_1] \,,
\end{equation}
and
\begin{equation} 
\M_1^{-1} = \frac{1}{\det(\M_1)}\biggl[\I + \gamma C_1 \E + \gamma 
\left(\begin{array}{cc} R_2 - G_{1,2} & R_2-G_{1,2} \\ 
R_1 - G_{1,2} & R_1 - G_{1,2} \\ \end{array} \right)\biggr] \,,
\end{equation}
from which we calculate
\begin{equation}
\M_1^{-1} \B_1 = \frac{\pi \gamma/2}{\det(\M_1)} 
\left(\begin{array}{cc} 1/\ve_1 & - 1/\ve_2 \\ -1/\ve_1 & 1/\ve_2 \\
      \end{array} \right)\,.
\end{equation}
The two eigenvalues of this matrix are $\hat{\sigma}_0 = 0$ and
\begin{equation}
\hat{\sigma}_1 = \frac{\pi \gamma(1/\ve_1 + 1/\ve_2)}{2\,\det(\M_1)} \,,
\end{equation}
so that upon solving for $\hat{\sigma}_1^{-1}$, we get
\begin{equation} 
\frac{1}{\hat{\sigma}_1} = \frac{2\ve_1\ve_2}{\pi(\ve_1+\ve_2)} 
\biggl[\frac{1}{\gamma} + (R_1 + R_2 - 2G_{1,2}) + 2C_1\biggr] \,.
\end{equation}
To simplify this expression, we use
\begin{equation}
\gamma = \frac{1}{1/\nu_1 + 1/\nu_2} = \frac{1}{-\ln(\ve_1\ve_2)} \,.
\end{equation}
This yields the following asymptotic behavior for the Steklov
eigenvalue:
\begin{align}  
  \frac{1}{\sigma_1} \approx \frac{2\ve_1\ve_2}{\pi(\ve_1+\ve_2)}
  \biggl[-\ln(\ve_1\ve_2) + 2C_1 + (R_1 + R_2 - 2G_{1,2})\biggr]\,.
\end{align}

For instance, if $\Omega$ is the unit disk, Eq. (\ref{eq:Green_disk})
yields
\begin{align}
  \frac{1}{\sigma_1} \approx \frac{2\ve_1\ve_2}{\pi(\ve_1+\ve_2)}
  \biggl(-\ln(\ve_1\ve_2) + 2C_1 + 2\ln |\x_1-\x_2|\biggr)\,.
\end{align}

Moreover, the corresponding eigenfunction $V_1$ can be easily found by
noting that $A_1 = -A_2$ from Eq. (\ref{eq:Asum}), whereas
Eq. (\ref{eq:SteklovSN_Ai_norm}) implies
\begin{equation}
A_1 = -A_2 \approx \frac{\sqrt{2} \sigma_1}{\pi \sqrt{1/\ve_1 + 1/\ve_2}} \,.
\end{equation}
These coefficients determine $(V_1)|_{\Gamma_i}$ via
Eq. (\ref{eq:SteklovSN_Vi}).

\subsection{Example of identical equally-spaced patches on the
  boundary of the unit disk}

When all patches are of the same size, $\ve_j = \ve$, we have $\nnu =
\nu \I$ and $\eeta = \I \pi/(2\ve)$, so that
\begin{equation}
\M_1  = \I + \nu \biggl(\I - \frac{\E}{N}\biggr) (\G + C_1 \I)\,, \qquad
\B_1  = \frac{\pi \nu}{2\ve} \biggl(\I - \frac{\E}{N}\biggr) \,.
\end{equation}
If the patches are equally-spaced on the boundary of the unit disk,
$\G$ is circulant and symmetric, and its eigenvectors and eigenvalues
were given in Eqs. (\ref{eq:qj_def}, \ref{eq:kappa_def}).  As shown
earlier in Sec. \ref{sec:Dirichlet_N}, one has $(\I - \E/N) \q_j =
\q_j$ for any $j = 1,2,\ldots,N-1$, and $(\I - \E/N) \q_N = 0$.  As a
consequence, $\M_1 \q_N = \q_N$ and $\M_1 \q_j = (1 + \nu
(\kappa_j+C_1)) \q_j$, where $\kappa_j$ are given explicitly by
Eq. (\ref{eq:kappaj}).  Upon calculating $\M_1^{-1}\q_j$, we readily
find that
\begin{equation*}
  \B_1 \M_1^{-1} \q_N = 0\,, \qquad \B_1 \M_1^{-1} \q_j =
  \frac{\pi \nu}{2\ve} \bigl(1 + \nu (\kappa_j+C_1)\bigr)^{-1} \q_j \,.
\end{equation*}
Since the eigenvalues of the matrices $\B_1 \M_1^{-1}$ and
$\M_1^{-1}\B_1$ are identical, we have from Eq. (\ref{eq:Mat_main})
that $\hat{\sigma}_0 = 0$ and
\begin{equation}\label{eq:hat_sigma}
  \hat{\sigma}_j = \frac{\pi \nu}{2\ve[1 + \nu (\kappa_j+C_1)]}\,,  \qquad
  j\in \lbrace{1,\ldots,N-1\rbrace}\,.
\end{equation}
Substituting $C_1={3/2}-\ln{2}$ from Eq. (\ref{eq:C0_C1}) and
$\nu={-1/\ln\ve}$ into Eq. (\ref{eq:hat_sigma}), we obtain the
following small-$\ve$ asymptotic result for the first $N-1$
eigenvalues of the mixed Steklov-Neumann problem:
\begin{equation}
  \frac{1}{\ve \sigma_j} \approx \frac{2}{\pi} \biggl(-\ln(\ve) +
  \frac32 - \ln(2) + \kappa_j \biggr)\,,
\end{equation}
for $j = 1,2,\ldots,N-1$.  When the number of patches is large,
i.e. $N \gg 1$, while still enforcing the well-separated patch
assumption $\ve N\ll \pi$, we can use the asymptotic relation
(\ref{eqb:kappa_approx}) for $\kappa_j$, valid for $j\ll N$, to
conclude that 
\begin{equation}
  \frac{1}{\ve \sigma_j} \approx \frac{2}{\pi}
\biggl( \frac{N}{2j} -\ln\left(\frac{N\ve}{\pi}\right) - \frac{1}{3}\biggr)\,.
\end{equation}
If $N\ve$ is not too small, the logarithmic and constant terms can be
neglected to yield to a first approximation
\begin{equation}
\sigma_j \approx \frac{\pi j}{N\ve} \quad \mbox{for} \quad j\ll N\,,
\end{equation} 
when $N\gg 1$.  It is instructive to compare this approximation to the
case of a single Steklov patch of half-length $\ve_1 = N\ve$, for
which $\sigma_j \approx \mu_j/\ve_1 \approx (\pi j/2)/(N\ve)$, where
we used the asymptotic relation (\ref{eq:muk_asympt}).  As a
consequence, the configuration with a single patch of half-length
$N\ve$ yields approximately twice smaller eigenvalues.  This suggests
that the fragmentation of a patch will increase the eigenvalues.

\section{Steklov-Neumann-Dirichlet problem}
\label{sec:SteklovND}

In this section, we consider the last setting of a single Steklov
patch $\Gamma_{\ve_1}$ and $N-1$ Dirichlet patches $\Gamma_{\ve_j}$
($j = 2,\ldots,N$).  This is a typical situation when the diffusing
particle needs to react on $\Gamma_{\ve_1}$ {\it before} escaping the
domain $\Omega$ through multiple opening windows $\Gamma_{\ve_2},
\ldots, \Gamma_{\ve_N}$.  A formal solution of such an escape problem
was provided in \cite{Grebenkov23} on the basis of the mixed
Steklov-Neumann-Dirichlet spectral problem, formulated as
\begin{subequations} \label{eq:Steklov_original}
\begin{align}
\Delta V & = 0 \quad \textrm{in}~\Omega \,,\\  \label{eq:Steklov}
\partial_n V & = \sigma V \quad \textrm{on}~ \Gamma_{\ve_1} \,; \qquad
V  = 0 \quad \textrm{on} ~ \Gamma_{\ve_2} \cup \cdots \cup \Gamma_{\ve_N} \,,\\
\partial_n V & = 0 \quad \textrm{on} ~ \pa_0 .
\end{align}
\end{subequations}
As previously, this spectral problem is known to have a discrete
positive spectrum \cite{Levitin}, i.e., infinitely many eigenpairs
$\{\sigma_k, V_k\}$ that are enumerated by $k = 0,1,\ldots$ to form an
increasing sequence of eigenvalues: $0 < \sigma_0 \leq \sigma_1 \leq
\ldots \nearrow \infty$.  The presence of Dirichlet patches implies
that the principal eigenvalue $\sigma_0$ is strictly positive.  We aim
at determining the asymptotic behavior of the eigenvalues and
eigenfunctions of this spectral problem in the small-target limit
$\ve_j\to 0$.

In the analysis below, we treat separately two cases depending on the
integral of the Steklov eigenfunction on the Steklov patch
$\Gamma_{\ve_1}$. In particular, if
\begin{equation}  \label{eq:V_intSteklov}
\int_{\Gamma_{\ve_1}} \partial_n V \, ds = \sigma \int_{\Gamma_{\ve_1}} V \, ds \ne 0,
\end{equation}
then, from the divergence theorem, the Steklov patch produces a
logarithmic contribution to the far field.  We will mainly focus on
this generic case.  However, if the integral in
Eq. (\ref{eq:V_intSteklov}) is zero (e.g., if $V|_{\Gamma_{\ve_1}}$ is
antisymmetric, see below), there is no logarithmic contribution, and
such an eigenfunction vanishes in the far field.  This situation is
actually simpler because the decay of $V$ away from the Steklov patch
is compatible with Dirichlet patches.  In other words, we can restrict
the analysis to the inner solution near the Steklov patch as if there
were no Dirichlet patches.  We will illustrate this situation in
Sec. \ref{sec:example_two}.

\subsection{Matched asymptotic analysis}
\label{sec:SND_matched}

Expectedly, we can combine formerly derived inner solutions for
Dirichlet and Steklov patches, whereas the outer solution is still
written as the linear combination (\ref{eq:Vsolution}).  As a
consequence, we must enforce the singularity behavior
\begin{equation}  \label{eq:Vasympt_SND}
  V \sim A_j \biggl[\ln|\x-\x_j| + 1/\nu_j + \delta_{j,1} \Cmu(-\sigma \ve_1) +
  o(1)\biggr] \,,
\end{equation}
as $\x\to\x_j$ for all patches $j\in \lbrace{1,\ldots,N\rbrace}$.
Here we have used the former definition $\nu_j = -1/\ln(\ve_j/2)$ for
Dirichlet patches and $\nu_1 = -1/\ln(\ve_1)$ for the Steklov patch.
By ensuring that $V$ in Eq. (\ref{eq:Vsolution}) satisfies
Eq. (\ref{eq:Vasympt_SND}) we obtain that
\begin{equation}  \label{eq:A0}  
A_j + \nu_j \biggl[A_j R_j + \sum\limits_{i=1 \atop i\ne j}^N  G_{j,i} A_i\biggr] 
= \chi \nu_j - \nu_1 A_1 \Cmu(-\sigma \ve_1) \delta_{j,1} \,.
\end{equation}
Together with Eq. (\ref{eq:Asum}), they form a system of $N+1$ linear
equations for the unknowns $A_i$ and $\chi$.  Using the former matrix
notation in Eqs. (\ref{eq:matrices_def}), Eq.~(\ref{eq:A0}) becomes
\begin{equation}  \label{eq:A01}
  \bigl(\I + \nnu \G\bigr) \A = \chi \nnu \e - \nu_1 \Cmu(-\sigma \ve_1) \e_1
  \e_1^\dagger \A\,,
\end{equation}
with $\e^\dagger \A = 0$.  Upon left-multiplying by $\e^\dagger$, we
get
\begin{equation}  \label{eq:chi}
  \chi = \frac{1}{\bar{\nu}} \biggl[\e^\dagger \nnu \G \A + \Cmu(-\sigma \ve_1)
  \nu_1 \e^\dagger \E_1 \A\biggr]\,,
\end{equation}
where $\bar{\nu}$ was defined in Eq. (\ref{eq:barnu}), and we
introduced the matrix $\E_1 = \e_1 \e_1^\dagger$ for a shorter
notation.  Eliminating $\chi$ from Eq. (\ref{eq:A01}), we obtain that
\begin{equation}  \label{eq:system2}
\biggl(\I + \nnu \biggl(\I - \frac{\E \nnu}{\bar{\nu}} \biggr) \G \biggr) \A 
+ \Cmu(-\sigma \ve_1) \nu_1 \biggl(\I - \frac{\nnu \E}{\bar{\nu}} \biggr)
\E_1 \A = {\bf 0} \,,
\end{equation}
where we recall that $\E = \e \e^\dagger$.  Upon introducing $\B$ by
\begin{equation} 
\B  = \biggl(\I - \frac{\nnu \E}{\bar{\nu}} \biggr) \E_1 \,,
\end{equation}
and using the matrix $\M_0$ from Eq. (\ref{eq:M0}), we rewrite the
matrix system in Eq. (\ref{eq:system2}) as
\begin{equation}  \label{eq:system_A}
\bigl(\M_0 + \Cmu(-\sigma \ve_1) \nu_1 \B\bigr) \A = {\bf 0} \,.
\end{equation}
The condition, under which this matrix equation admits a
nontrivial solution is
\begin{equation}
\det\bigl(\M_0 + \Cmu(-\sigma\ve_1) \nu_1 \B\bigr) = 0\,,
\end{equation}
which is a scalar problem that determines $\sigma$.

To rewrite this problem in a more explicit form, we first observe that
$\rank(\B) = 1$ since $\E_1 \q = 0$ for any vector $\q \in \R^N$ such
that $\e_1^\dagger \q = 0$.  Since $\e^\dagger \e_1 = 1$, we can
rewrite $\B$ in a more convenient rank-one form as
\begin{equation}
  \B = \biggl(\e_1 - \frac{\nnu \e}{\bar{\nu}} \biggr) \e_1^\dagger =
  \a \b^\dagger\,, \quad \mbox{where} \quad   \a = \e_1 -
  \frac{\nnu \e}{\bar{\nu}} \,, \qquad \b =
  \Cmu(-\sigma\ve_1) \nu_1 \e_1 \,.
\end{equation}

To proceed, we need the matrix determinant lemma \cite{Ding07}.

{\bf Lemma}: {\it Let $\M = \M_0 + \a \b^\dagger$ be a perturbation of an
invertible matrix $\M_0 \in \R^{N,N}$ by a rank-one matrix $\a
\b^\dagger$.  Then 
\begin{equation}
\det(\M_0 + \a \b^\dagger) = (1 + \b^\dagger \M_0^{-1} \a)\, \det(\M_0)\,.
\end{equation}
It follows that $\det(\M) = 0$ if and only if $\b^\dagger \M_0^{-1}
\a = -1$.}

In the small-target limit, all $\nu_j \ll 1$ so that the matrix $\M_0$
is invertible since it is a small perturbation of the identity matrix
in Eq. (\ref{eq:M0}).  Applying the lemma above to our setting, we
determine the condition on $\sigma$ as
\begin{equation}  \label{eq:Cmu_N}
\Cmu(-\sigma \ve_1) = C \,,
\end{equation}
where
\begin{equation} \label{eq:Cdef}
  C = - \frac{1}{\nu_1} \biggl( \e_1^\dagger \M_0^{-1}
  \biggl[\e_1 - \frac{\nnu \e}{\bar{\nu}}\biggr]\biggr)^{-1} \,,
\end{equation}
with the vectors and matrices $\e$, $\e_1$, $\nnu$, and $\M_0$ being
defined in Eqs. (\ref{eq:matrices_def}, \ref{eq:M0}).  This is the
main result of this section that will allow us to determine the
asymptotic behavior of the Steklov eigenvalues and their dependence on
the configuration and sizes of all patches that are captured via the
constant $C$ in Eq. (\ref{eq:Cdef}).  We further emphasize that the
homogeneous matrix equation (\ref{eq:system_A}) cannot uniquely
determine the coefficients $A_i$.  In fact, an eigenfunction $V$ can
be found up to a multiplicative factor that has to be fixed by
normalization (see below).

As stated above, the matrix $\M_0$ is a small perturbation of the
identity matrix in the small-target limit, so that $\M_0^{-1} \sim \I$
to leading order, which implies that $\e_1^\dagger \M_0^{-1} [\e_1 -
\nnu \e/\bar{\nu}] \sim 1 - \nu_1/\bar{\nu} > 0$.  Since $0 < \nu_1
\ll 1$, we conclude that the constant $C$ in Eq. (\ref{eq:Cdef}) is
negative and large:
\begin{equation}  \label{eq:Cnegative}
C < 0\,, \qquad |C| \gg 1\,.
\end{equation}

\subsection{Asymptotic behavior of eigenvalues and eigenfunctions}

Denoting $\mu = -\sigma\ve_1$, we recast Eq. (\ref{eq:Cmu_N}) as
\begin{equation}  \label{eq:Cmu_C}
\Cmu(\mu) = C \,.
\end{equation}
The spectral expansion (\ref{eq:Cmu_spectral}) of the function
$\Cmu(\mu)$ allows one to solve this equation numerically for any
fixed negative value $C$ given by Eq. (\ref{eq:Cdef}).  Since the
derivative $\Cmu^{\prime}(\mu) = d\Cmu(\mu)/d\mu$ is negative, $\Cmu(\mu)$ is
a continuous and monotonically decreasing function on each interval
$(-\mu_{2j+2}, -\mu_{2j})$, with $j = 0,1,\ldots$.  Moreover, it
ranges from $+\infty$ to $-\infty$ on each interval.  As a
consequence, for any fixed value $C$, there exist infinitely many {\it
negative} solutions of Eq. (\ref{eq:Cmu_C}), denoted as
$-\hat{\mu}_{2j}$, such that
\begin{equation}
  \mu_{2j} \leq \hat{\mu}_{2j} \leq \mu_{2j+2} \quad \mbox{for} ~
  j = 0, 1, \ldots\,.
\end{equation}
This property facilitates the numerical solution, as a single zero has
to be searched on each interval.  Moreover, as the coefficients
$[\Psi_{2k}(\infty)]^2$ are small (see Table \ref{tab:mu_int}) and
decrease with $k$, whereas $C$ is negative and large, one has
$\hat{\mu}_{2j} \approx \mu_{2j}$ for $j > 0$.  This is consistent
with the intuitive picture that large Steklov eigenvalues become
insensitive to Dirichlet patches in the small-target limit, and one
retrieves the asymptotic behavior for a single Steklov patch
\cite{Grebenkov25}.

The solutions $\hat{\mu}_{2j} \approx \mu_{2j}$ determine the
leading-order term in the asymptotic behavior of the eigenvalues
$\sigma_{2j}$:
\begin{equation}   \label{eq:sigmak_asympt}
  \sigma_{2j} \approx \frac{\mu_{2j}}{\ve_1}   \qquad \mbox{for}
  ~ j = 1,2,\ldots \,.
\end{equation}
In contrast, the smallest eigenvalue $\sigma_0$ involves the solution
$\hat{\mu}_0$, which may actually be small in the small-target limit.
We discuss this case separately in Sec. \ref{eq:sigma0}.

We also mention that the analysis above provides the leading-order
approximation to the associated Steklov eigenfunction, restricted to
$\Gamma_{\ve_1}$.  We recall that the inner solution near the Steklov
patch is $V_{2j} \sim A_1 g_{-\sigma_{2j} \ve_1}(\y)$, with the
Green's function $g_\mu(\y)$ given by Eq. (\ref{eq:gR}).  As a
consequence, its restriction onto the patch reads
\begin{equation}  \label{eq:Vj_approx}
V_{2j}\biggr|_{\Gamma_{\ve_1}} \approx a_{2j} \sum\limits_{k=0}^\infty 
\frac{\Psi_{2k}(\infty)}{\mu_{2k} - \ve_1 \sigma_{2j}} \Psi_{2k}(\y)\,,
\end{equation}
where the proportionality coefficient $a_{2j}$ is fixed by the
conventional normalization of the Steklov eigenfunction:
\begin{equation}  \label{eq:Vj_norm}
  1 = \int\limits_{\Gamma_{\ve_1}} V_{2j}^2\, ds \approx \ve_1 a_{2j}^2
  \sum\limits_{k=0}^\infty 
\frac{[\Psi_{2k}(\infty)]^2}{(\mu_{2k} - \ve_1 \sigma_{2j})^2} \,,
\end{equation}
where the orthogonality of $\Psi_{2k}$ was used.  Since $\ve_1
\sigma_{2j} \approx \mu_{2j}$, the eigenfunction $\Psi_{2j}$ provides
the dominant contribution, and one gets for each $j\in
\lbrace{1,2,\ldots\rbrace}$ that
\begin{equation}  \label{eq:Vj_approx2}
V_{2j}(\x_1 + \ve_1 \Q_1^\dagger (y_1,0)^\dagger) \approx
\frac{1}{\sqrt{\ve_1}} \Psi_{2j}(y_1,0)
\end{equation}
on the Steklov patch (i.e., for $|y_1| \leq 1$).

We emphasize that the analysis above allowed us to access only half of
eigenvalues with even indices $2j$ that correspond to symmetric
eigenmodes.  In turn, the eigenvalues with odd indices correspond to
antisymmetric eigenmodes, for which the integral over the Steklov
patch is zero.  As discussed at the beginning of
Sec. \ref{sec:SteklovND}, such eigenfunctions vanish away from the
Steklov patch so that their asymptotic behavior can be determined
directly from the local solution.  As a consequence, we get a
leading-order approximation
\begin{equation}
  \sigma_{2j+1} \approx \frac{\mu_{2j+1}}{\ve_1}   \qquad \mbox{for}
  ~ j=0,1,\ldots \,,
\end{equation}
and on the Steklov patch we have
\begin{equation}  \label{eq:Vj_approx3}
  V_{2j+1}(\x_1 + \ve_1 \Q_1^\dagger (y_1,0)^\dagger) \approx
  \frac{1}{\sqrt{\ve_1}} \Psi_{2j+1}(y_1,0) \,.
\end{equation}

In summary, our analysis justifies theoretically the intuitively
appealing scaling argument that the eigenvalues $\sigma_j$ and
eigenfunctions $(V_j)|_{\Gamma_{\ve_1}}$ (on the Steklov patch) can be
approximated by the eigenvalues $\mu_j$ and eigenfunctions $\Psi_j$ of
the auxiliary problem (\ref{eq:Vkint}), as if there were no Dirichlet
patches.  Regardless of the symmetry of eigenfunctions (and parity of
its index), we can combine the former leading-order approximations as
\begin{subequations}  \label{eq:sigma_SND}
\begin{align}
  \sigma_j & \approx \frac{\mu_j}{\ve_1}   \quad \mbox{for}
             ~ j=1,2,\ldots \,, \\
  V_j(\x_1 + \ve_1 \Q_1^\dagger (y_1,0)^\dagger) & \approx
          \frac{1}{\sqrt{\ve_1}} \Psi_j(y_1,0) \,.
\end{align}
\end{subequations}
In contrast, the presence of Dirichlet patches must affect the
principal eigenvalue $\sigma_0$ and the associated eigenfunction
$V_0$, as explained below.


\subsection{The principal eigenvalue}
\label{eq:sigma0}

Since $\mu_0 = 0$, the smallest solution of Eq. (\ref{eq:Cmu_C}),
$\hat{\mu}_0$, is close to $0$.  Indeed, as the constant $C$ is large
and negative, one needs to have $|\mu| \ll 1$ to ensure that
$\Cmu(\mu)$ is also large and negative.  Under the condition $|\mu|\ll
1$, we can use the approximation (\ref{eq:Cmu_approx}), which can be
easily inverted to get the explicit result
\begin{equation}  \label{eq:mu_Cmu}
  \frac{1}{\mu} \approx \frac{2}{\pi} \bigl[\Cmu(\mu) - C_1 \bigr] =
  \frac{2}{\pi} \bigl[C - C_1\bigr] \,.
\end{equation}
As a consequence, substitution of Eq. (\ref{eq:Cdef}) here yields the
asymptotic behavior of the principal eigenvalue $\sigma_0$:
\begin{equation}  \label{eq:mu0_main}
\frac{1}{\ve_1 \sigma_0} \approx - \frac{2}{\pi} \biggl[\ln(\ve_1) 
\biggl( \e_1^\dagger \M_0^{-1} \biggl[\e_1 -
\frac{\nnu \e}{\bar{\nu}}\biggr]\biggr)^{-1} - C_1\biggr]\,.
\end{equation}
In sharp contrast to Eq. (\ref{eq:sigmak_asympt}) for $\sigma_j$ with
$j \geq 1$, the principal eigenvalue $\sigma_0$ exhibits a slower
divergence ${\mathcal O}(1/(\ve_1\ln(\ve_1)))$.  This is one of the
main results of this section.  The associated eigenfunction $V_0$ is
given by Eq. (\ref{eq:Vj_approx}) with the normalization condition
(\ref{eq:Vj_norm}).  We stress that Eq. (\ref{eq:Vj_approx}) cannot be
reduced to the approximation (\ref{eq:Vj_approx2}) in this case.

\subsection{Example of two patches}
\label{sec:example_two}

When there are two patches ($N = 2$), Eqs. (\ref{eq:Cdef},
\ref{eq:mu0_main}) can be readily solved.  Substituting
$\M_0^{-1}$ from Eq. (\ref{eq:M0_inv_two}) and
\begin{equation}
\e_1 - \frac{\nnu \e}{\bar{\nu}} = \frac{\nu_2}{\nu_1+\nu_2}
  \left(1,-1\right)^{\dagger}
\end{equation}
into Eq. (\ref{eq:Cdef}), we get after simplifications that
\begin{equation}  \label{eq:Cmu_1} 
C = \ln(\ve_1 \ve_2/2) - (R_1+R_2-2G_{1,2}) \,.
\end{equation}
Using Eq. (\ref{eq:mu_Cmu}) with $C_1={3/2}-\ln{2}$ from
Eq. (\ref{eq:C0_C1}), the asymptotic behavior of the principal
eigenvalue is
\begin{equation} \label{eq:mu_1} 
  \frac{1}{\ve_1 \sigma_0} \approx \frac{2}{\pi}\biggl( -\ln(\ve_1 \ve_2) +
  \frac32 + (R_1+R_2-2G_{1,2})\biggr) \,.
\end{equation}

\begin{figure}
\begin{center}
\includegraphics[width=40mm]{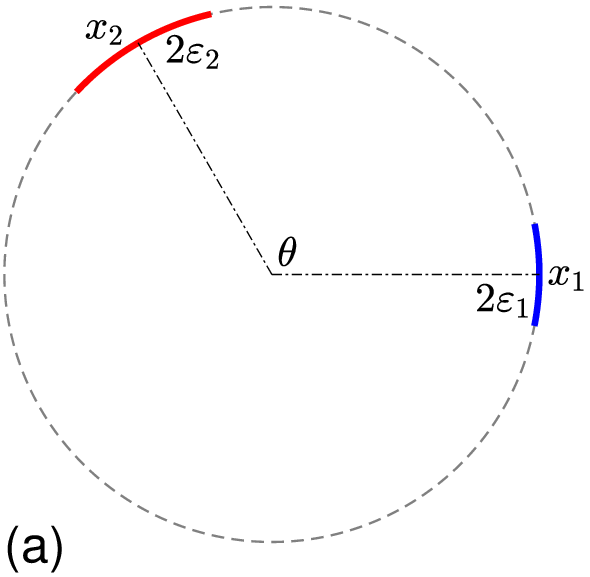} 
\includegraphics[width=40mm]{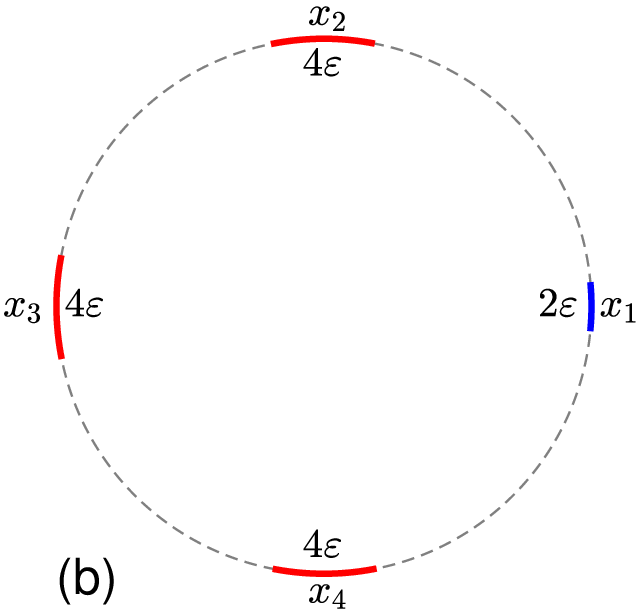} 
\end{center}
\caption{
Illustration for the unit disk with Steklov and Dirichlet patches.
{\bf (a)} One Steklov patch of length $2\ve_1 = 0.4$ at $\x_1 = (1,0)$
(blue) and one Dirichlet patch of length $2\ve_2 = 0.6$ (red), whose
center $\x_2$ is at angle $\theta = 2\pi/3$.  {\bf (b)} One Steklov
patch of length $2\ve_1 = 2\ve = 0.2$ at $\x_1 = (1,0)$ (blue)
and three Dirichlet patches of length $2\ve_j = 0.4$ (red),
whose centers $\x_j$ are equally-spaced on the boundary of the unit
disk. }
\label{fig:scheme_disk}
\end{figure}

\begin{figure}
\begin{center}
\includegraphics[width=88mm]{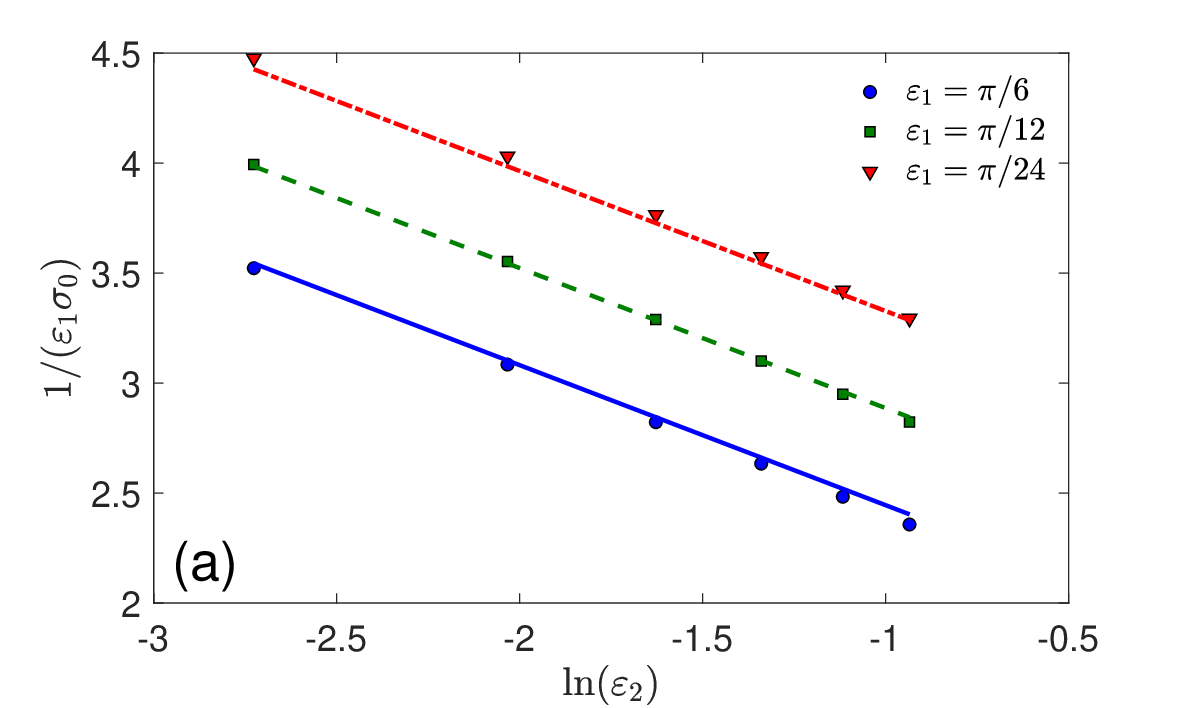} 
\includegraphics[width=88mm]{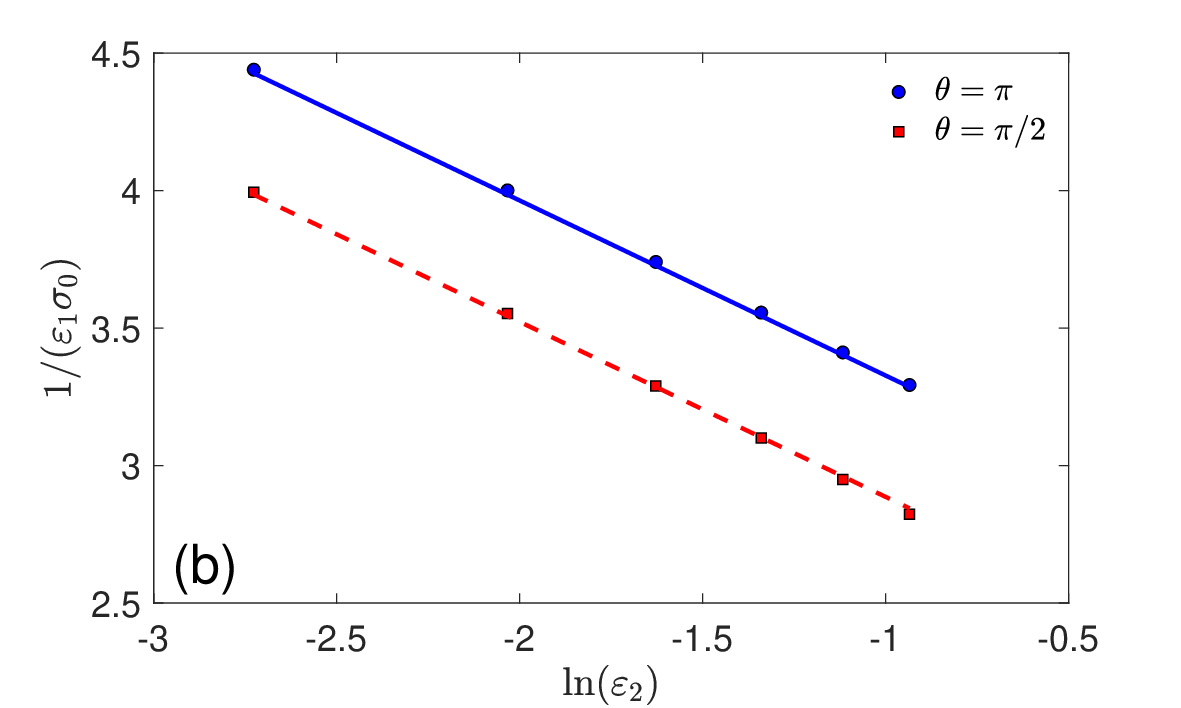} 
\end{center}
\caption{
Dependence of $1/(\ve_1\sigma_0)$ on $\ve_2$ for the unit disk with a
Steklov patch of length $2\ve_1$ (located at $\x_1 = (1,0)$), and one
Dirichlet patch of length $2\ve_2$, located at $\x_2$.  Symbols
present the numerical solution by a FEM with the maximal
meshsize $h_{\rm max} = 0.005$ and lines show Eq. (\ref{eq:mu_disk}).
{\bf (a)} $\x_2 = (0,1)$ and three values of $\ve_1$: $\ve_1 =
\pi/6$ (circles), $\ve_1 = \pi/12$ (squares), and $\ve_1 = \pi/24$
(triangles).
{\bf (b)} $\ve_1 = \pi/12$ and two locations of the Dirichlet patch:
$\x_2 = (-1,0)$ (circles, $\theta = \pi$), and $\x_2 = (0,1)$
(squares, angle $\theta = \pi/2$).  }
\label{fig:disk}
\end{figure}

For instance, when $\Omega$ is the unit disk, one can substitute
Eq. (\ref{eq:Green_disk}) into Eq. (\ref{eq:mu_1}) to get
\begin{equation}    \label{eq:mu_disk}
  \frac{1}{\ve_1 \sigma_0} \approx \frac{2}{\pi}
  \biggl( -\ln(\ve_1 \ve_2) + \frac32 + 2\ln |\x_1 - \x_2| \biggr) \,.
\end{equation}
Figure \ref{fig:disk} illustrates the remarkable accuracy of this
asymptotic relation.

\begin{figure}
\begin{center}
\includegraphics[width=42mm]{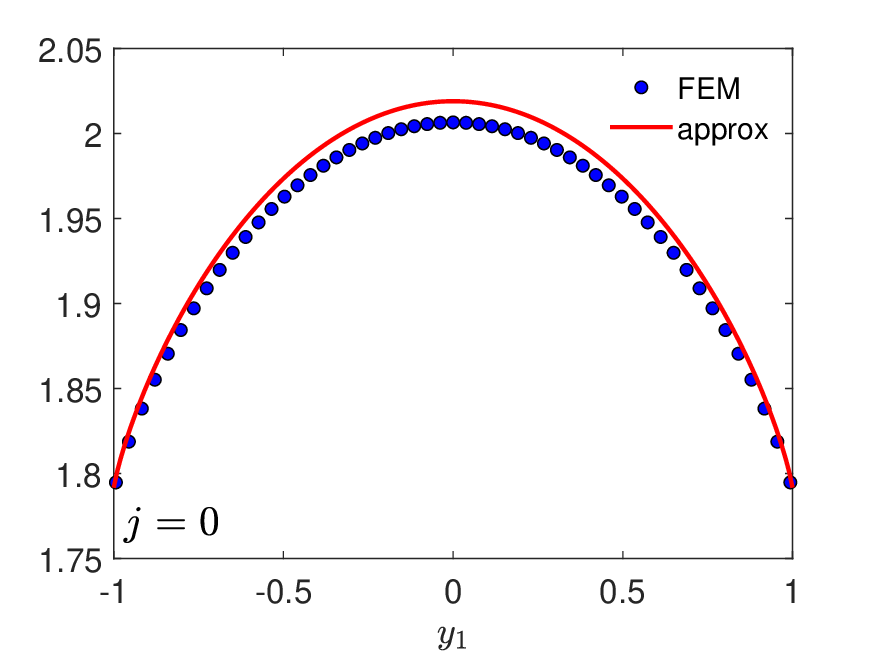} 
\includegraphics[width=42mm]{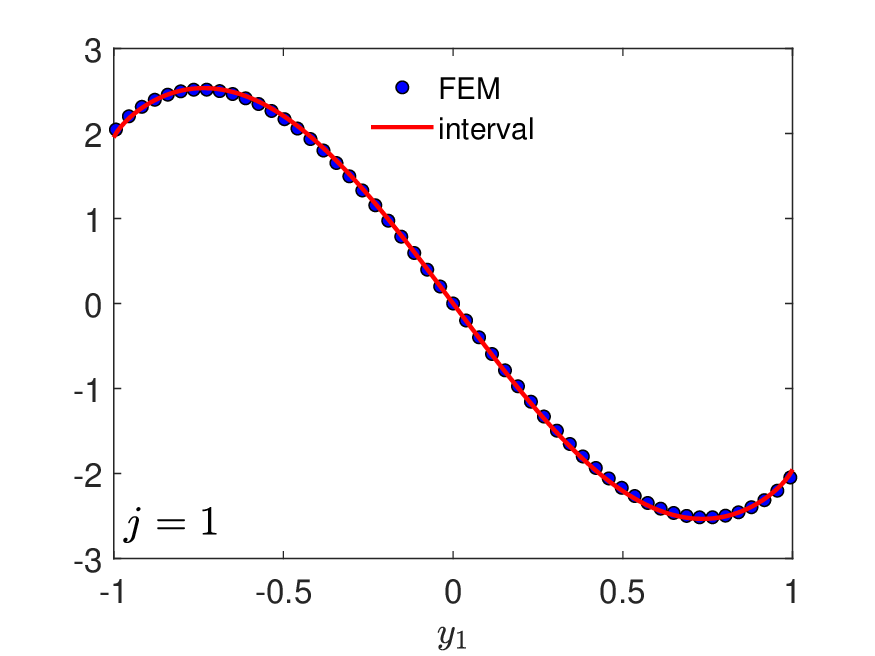} 
\includegraphics[width=42mm]{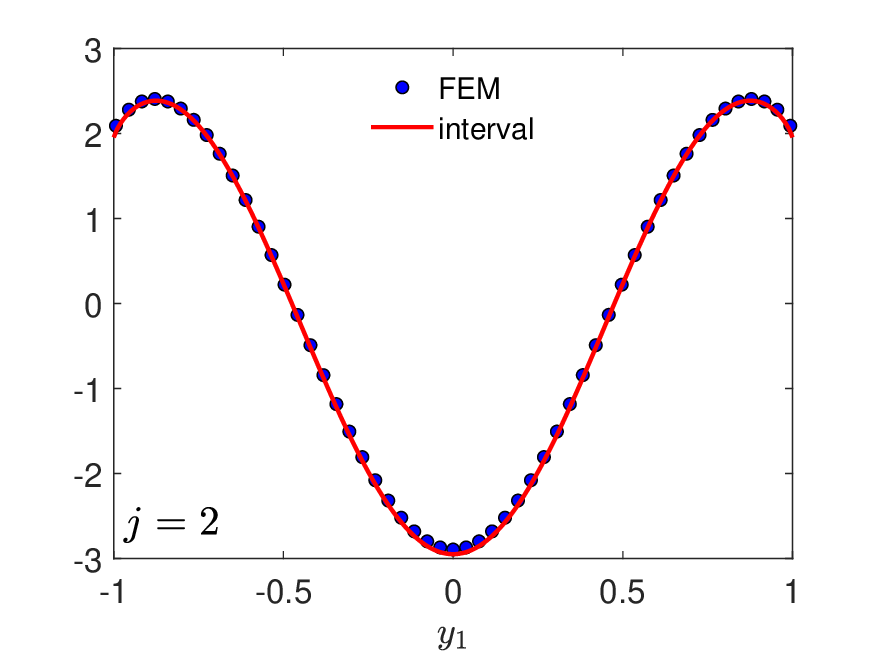} 
\includegraphics[width=42mm]{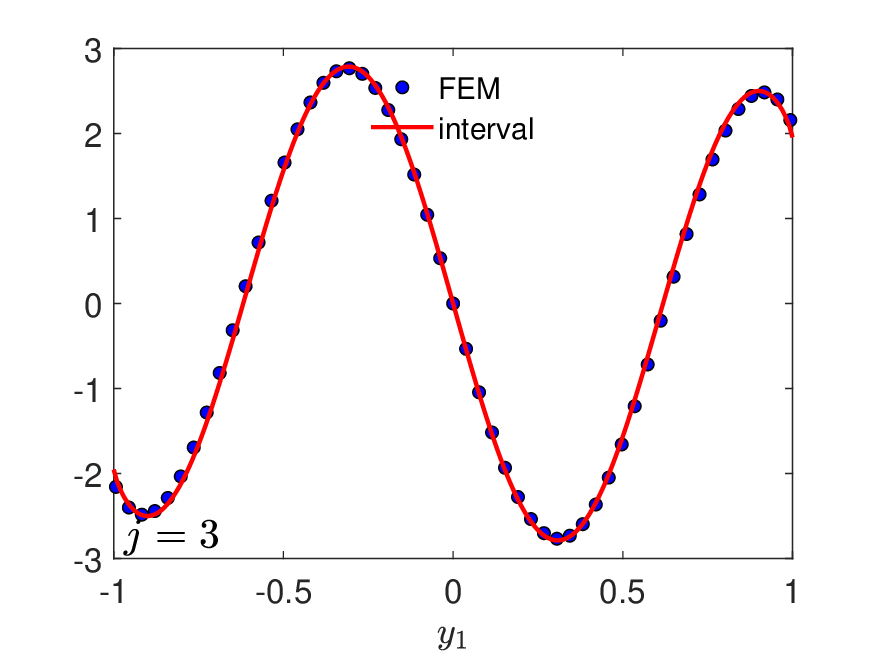} 
\end{center}
\caption{
The eigenfunctions $V_j$ restricted on the Steklov patch
$\Gamma_{\ve_1}$, for the unit disk with a Steklov patch of length
$2\ve_1 = \pi/12 \approx 0.26$ (located at $\x_1 = (1,0)$), and one
Dirichlet patch of length $2\ve_2 = \pi/6 \approx 0.52$, located at
$\x_2 = (-1,0)$.  Filled circles present the numerical solution by a
FEM with the maximal meshsize $h_{\rm max} = 0.005$, while solid lines
show Eq. (\ref{eq:Vj_approx}) for $j = 0$ and Eq. (\ref{eq:sigma_SND})
for $j > 0$.  Four panels present the cases $j = 0,1,2,3$.}
\label{fig:disk_vk}
\end{figure}

Figure \ref{fig:disk_vk} shows the behavior of the Steklov
eigenfunctions $V_j$ restricted onto the Steklov patch.  For the
principal eigenmode with $j = 0$, this restriction is positive, as
expected.  The asymptotic formula (\ref{eq:Vj_approx}) yields an
accurate approximation.  Let us now look at other eigenmodes with $j =
1,2,3$, for which $\sigma_j \approx \mu_j/\ve_1$.  We see that the
restriction of $V_j$ and its approximation (\ref{eq:sigma_SND}) are in
excellent agreement, for both symmetric and antisymmetric
eigenfunctions, even though both considered patches are not
small.

\subsection{Example of equally-spaced patches on the unit disk}

We now consider another setting where the general formulas
(\ref{eq:Cmu_N}) and (\ref{eq:Cdef}) can be simplified.  We suppose
that all Dirichlet patches with $j = 2,\ldots,N$ are of the same
length, so that $\ve_j = \ve$ for $j\in\lbrace{2,\ldots,N\rbrace}$,
whereas for the Steklov patch we have $\ve_1 = \ell_1 \ve/2$, for some
$\ell_1>0$. To treat this case, we can impose in our general formula
(\ref{eq:Cdef}) for $C$ that
\begin{equation}
\nu=\nu_1 = \ldots = \nu_N = -{1/\ln(\ve/2)}, \qquad \nnu=\nu \I ,
\end{equation}
provided that we shift $\Cmu(-\sigma\ve_1)$ appropriately in the
relation (\ref{eq:Cmu_N}). To determine this shift, we write the
singularity condition (\ref{eq:Vasympt_SND}) as $\x\to\x_1$ associated
with the Steklov patch $j=1$ as
\begin{align*}  
  V & \sim A_1 \biggl[\ln|\x-\x_1| - \ln\ve_1 + \Cmu(-\sigma \ve_1) +
      o(1)\biggr]  \\
  & \sim A_1 \biggl[\ln|\x-\x_1| + 1/\nu + \tilde{\Cmu}(-\sigma \ve_1)
    + o(1)\biggr]\,,
\end{align*}
where we have defined $\tilde{\Cmu}(-\sigma \ve_1)$ by
\begin{equation}\label{eq:new_C}
  \tilde{\Cmu}(-\sigma \ve_1)=\Cmu(-\sigma\ve_1) -
  \ln\left(\frac{2\ve_1}{\ve}\right) \,.
\end{equation}
As a result, by repeating the steps of Sec. \ref{sec:SND_matched}, we need
only replace Eq. (\ref{eq:Cmu_N}) by
\begin{equation}  \label{eq:Cmu_Nbis}
\Cmu(-\sigma \ve_1) - \ln\left( \frac{2\ve_1}{\ve} \right) = C \,,
\end{equation}
where $C$ is given by Eq. (\ref{eq:Cdef}) with $\nnu = \nu \I$,
$\bar{\nu} = N\nu$, and $\nu={-1/\ln\left({\ve/2}\right)}$, which yields
\begin{equation}  \label{eq:Cmu_multiple0}
  C = \frac{1}{\nu} \biggl(\frac{\e_1^\dagger \M_0^{-1} \e}{N}
  - \e_1^\dagger \M_0^{-1} \e_1\biggr)^{-1} \,.
\end{equation}

If all $\x_j$ are equally-spaced on the boundary of the unit disk
(such as shown in Fig. \ref{fig:scheme_disk}(b)), then the matrix $\G$ is
circulant and symmetric, so that its eigenvectors and eigenvalues are
known exactly (see Sec. \ref{sec:Dirichlet_N}).  Moreover, the matrix
$\M_0$ admits a spectral representation and thus can be inverted
explicitly.  Using Eq. (\ref{eq:M0_inv_N}), we get
\begin{equation}
  \e_1^\dagger \M_0^{-1} \e  = 1 \,, \qquad \e_1^\dagger \M_0^{-1} \e_1  =
  (\e_1^\dagger \q_N)  (\q_N^\dagger \e_1) + \sum\limits_{j=1}^{N-1}
  \frac{(\e_1^\dagger \q_j) \, (\q_j^\dagger \e_1)}{1 + \nu \kappa_j} \,,
  \label{eq:iprod_1}
\end{equation}
where $\q_j$ and $\kappa_j$ were defined by Eqs. (\ref{eq:qj_def},
\ref{eq:kappaj}).  However, since $(\e_1^\dagger \q_j) =
\omega^j/\sqrt{N}$ and $(\q_1^\dagger \e_1) = \omega^{-j}/\sqrt{N}$,
we get
\begin{equation}
  \e_1^\dagger \M_0^{-1} \e_1 = \frac{1}{N}\biggl[1 + \sum\limits_{j=1}^{N-1}
  (1 + \nu \kappa_j)^{-1} \biggr]\,.
\end{equation}
Substituting this expression together with Eq. (\ref{eq:iprod_1}) into
Eq. (\ref{eq:Cmu_multiple0}), we find
\begin{equation}  \label{eq:Cmu_multiple}
  C = - \frac{N}{\nu} \biggl(\sum\limits_{j=1}^{N-1} \frac{1}{1 + \nu \kappa_j}
  \biggr)^{-1} \,.
\end{equation}
As a result, Eqs. (\ref{eq:mu_Cmu}, \ref{eq:Cmu_Nbis}) with
$C_1={3/2}-\ln{2}$ imply
\begin{align}  \nonumber
  \frac{1}{\ve_1 \sigma_0} & \approx - \frac{2}{\pi} \bigl[\Cmu(-\sigma\ve_1) -
                             C_1 \bigr]
= - \frac{2}{\pi} \bigl[C + \ln\left(\frac{2\ve_1}{\ve}\right) - C_1\bigr] \\
  \label{eq:mu_multiple}
   & \approx \frac{2}{\pi} \biggl[\frac{N}{\sum\nolimits_{j=1}^{N-1}
                             \bigl(\ln(2/\ve) + \kappa_j\bigr)^{-1}} 
+ C_1 - \ln\left(\frac{2\ve_1}{\ve}\right)\biggr]\,.
\end{align}

\begin{figure}
\begin{center}
\includegraphics[width=88mm]{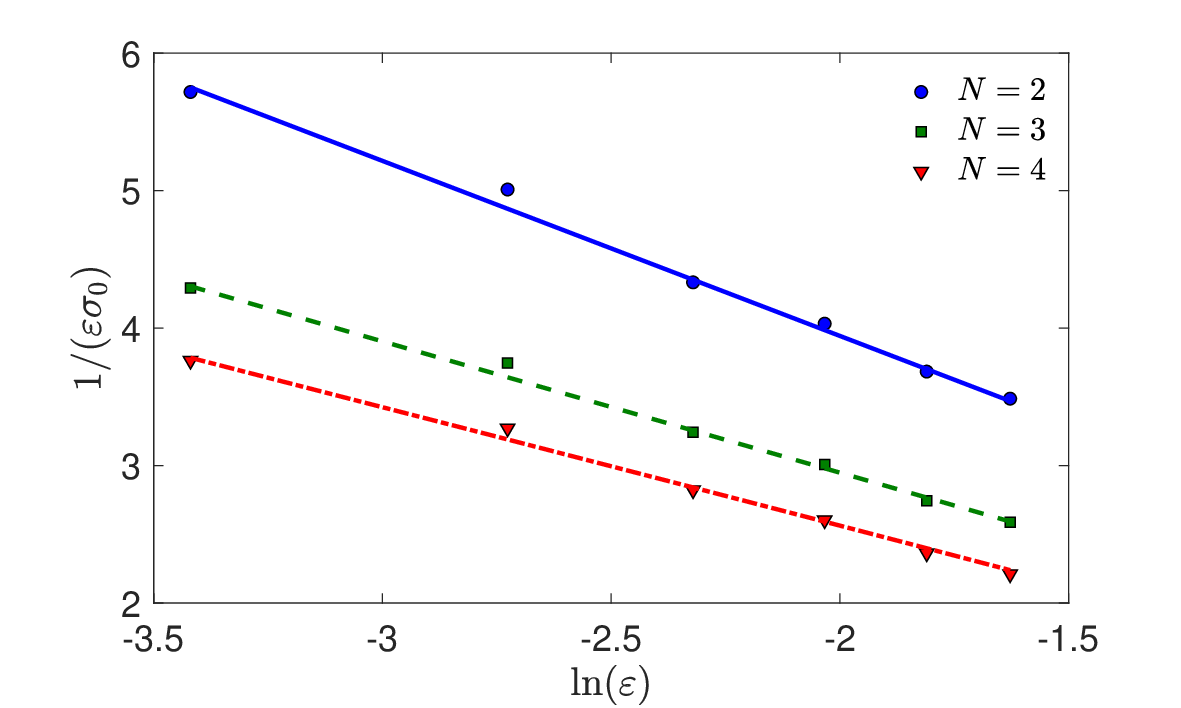} 
\end{center}
\caption{ 
Dependence of $1/(\ve\sigma_0)$ on $\ve$ for the unit disk with one
Steklov patch of length $2\ve$ (located at $\x_0 = (1,0)$), and $N-1$
Dirichlet patches of length $4\ve$, equally-spaced on the boundary of
the unit disk.  Symbols present the numerical solution by a FEM with
the maximal meshsize $h_{\rm max} = 0.005$ and lines show
Eq. (\ref{eq:mu_multiple}).}
\label{fig:disk_multiple}
\end{figure}

Figure \ref{fig:disk_multiple} illustrates the behavior of
$1/(\ve_1\sigma_0)$ as a function of $\ln(\ve)$ for the unit disk with
one Steklov patch of length $2\ve$, and several Dirichlet patches of
length $4\ve$, which are equally-spaced on the boundary of the unit
disk.  We observe an excellent agreement between the asymptotic
formula (\ref{eq:mu_multiple}) and numerical results.

Although the eigenvalues $\kappa_j$ are known explicitly via
Eq. (\ref{eq:kappaj}), it is instructive to inspect their asymptotic
behavior for large $N$ (see Appendix \ref{sec:kappa}).  We aim at
approximating the sum in the denominator of
Eq. (\ref{eq:mu_multiple}):
\begin{equation}  \label{eq:sum_def}
S = \frac{1}{N} \sum\limits_{j=1}^{N-1} \frac{1}{\ln(2/\ve) + \kappa_j} \,.
\end{equation}
The degeneracy $\kappa_{N-j} = \kappa_j$ allows us to limit this sum
to $j \leq {N/2}$ when $N$ is even.  By using the asymptotic result in
Eq. (\ref{eqb:kappa_final4}) for $\kappa_j$ when $N\gg 1$, we find
that
\begin{equation}\label{eq:sum_kappa_xi}
  \ln\left({2/\ve}\right) + \kappa_j \approx \frac{1}{2\xi} + \zeta +
  a \frac{\pi^2\xi^2}{9} \,,
\end{equation}
where we have defined $\xi={j/N}$, 
\begin{equation}\label{eq:sum_param}
  \zeta = -\ln\left( \frac{N\ve}{b}\right) \,, \quad \mbox{and}
  \quad b = 4\pi e^{-11/6}\approx 2.009\,.
\end{equation}
Here the coefficient $a = 1.25$ was empirically introduced to
improve the accuracy of the approximation of $\kappa_j$ (see Appendix
\ref{sec:kappa} for details).  In terms of $\zeta$, and substituting
${j/N} = \xi$, we view $S$ as a Riemannian approximation of the
integral defined by
\begin{equation} \label{eq:sum_approx} 
  S \approx S(\zeta) = 4
  \int\limits_0^{1/2} \frac{\xi \, d\xi}{1 + 2\xi
    \left(\zeta + a \pi^2\xi^2/9\right)} \,.
\end{equation}
As a result, Eq. (\ref{eq:mu_multiple}) is approximated for
$N\gg 1$ by
\begin{align}  \label{eq:mu_multiple2}
\frac{1}{\ve_1 \sigma_0} & \approx \frac{2}{\pi} 
\biggl[ \frac{1}{S(\zeta)} + C_1 - \ln\left(\frac{2\ve_1}{\ve}\right)\biggr],
\end{align}
where $\zeta$ is defined in Eq. (\ref{eq:sum_param}). In this way, we
have reduced the problem of estimating ${1/(\ve_1 \sigma_0)}$ to a
simple numerical quadrature of the function $S(\zeta)$ in
Eq.~(\ref{eq:sum_approx}).  To obtain a more explicit, but less
accurate, approximation, we neglect the ${a \pi^2\xi^2/9}$ term in
Eq. (\ref{eq:sum_approx}), which corresponds to using the result
(\ref{eqb:kappa_approx}) for $\kappa_j$, and then evaluate the
resulting integral to get $S= {\left(\zeta -
\ln(1+\zeta)\right)/\zeta^2}$.  In this way,
Eq. (\ref{eq:mu_multiple}) can be approximated more explicitly as
\begin{equation}  \label{eq:mu_loworder}
\frac{1}{\ve_1 \sigma_0}  \approx \frac{2}{\pi} 
 \left[\frac{\zeta^2}{\zeta - \ln(1+\zeta)} + C_1 - \ln\left(\frac{2\ve_1}{\ve}
                           \right)\right]\,.
\end{equation}

\begin{figure}
\begin{center}
\includegraphics[width=88mm]{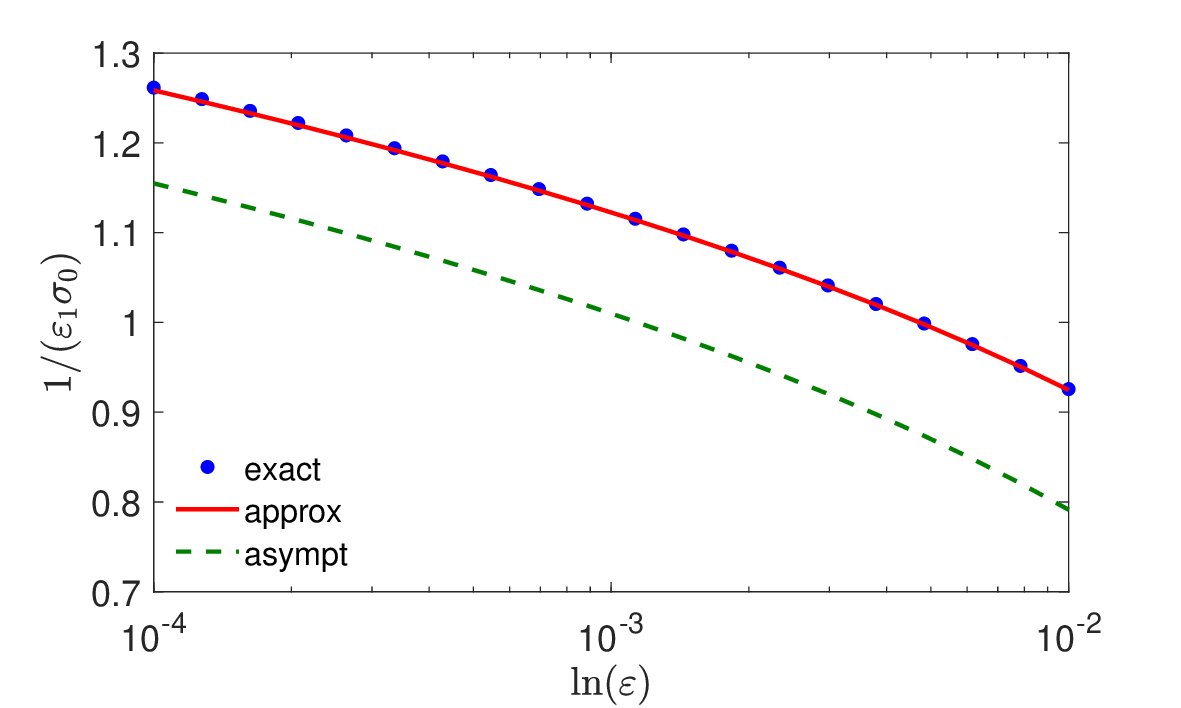} 
\end{center}
\caption{
Dependence of $1/(\ve_1 \sigma)$ on $\ve$ for the unit disk with the
Steklov patch of length $2\ve_1 = 0.2$ and $63$ Dirichlet patches
(each of length $2\ve$) that are equally-spaced on the boundary of the
unit disk.  Filled circles correspond to $\kappa_j$ obtained via the
discrete sum (\ref{eq:mu_multiple}), the solid line indicates its
large-$N$ approximation (\ref{eq:mu_multiple2}), and the dashed line
is the low-order approximation (\ref{eq:mu_loworder}).}
\label{fig:kappaj1}
\end{figure}

For $63$ Dirichlet patches (i.e.~$N=64$) and with $\ve_1=0.1$,
Fig. \ref{fig:kappaj1} compares the asymptotic results obtained by
using the discrete sum (\ref{eq:mu_multiple}) with its large-$N$
approximation (\ref{eq:mu_multiple2}) and with the simpler, more
explicit, result (\ref{eq:mu_loworder}).  We observe that
Eq. (\ref{eq:mu_multiple2}) provides an excellent approximation, while
Eq. (\ref{eq:mu_loworder}) has a small systematic underestimate.

\section{Further extensions}
\label{sec:extensions}

In the previous four sections, we progressively increased the
complexity of the problem: (i) $N$ perfectly reactive (Dirichlet)
patches; (ii) $N$ partially reactive (Robin) patches; (iii) $N$
imperfect (Steklov) patches; and (iv) one imperfect patch with $N-1$
Dirichlet patches.  We showed that the asymptotic analysis required to
study these settings is similar, although the resulting formulas
became progressively more intricate.  In the same vein, we can treat
any combination of perfectly reactive, partially reactive and
imperfect patches.

In this section, we briefly discuss two other extensions that are
relevant for applications: the case of interior targets
(Sec. \ref{sec:interior}) and the exterior problem
(Sec. \ref{sec:exterior}).

\subsection{Interior targets}
\label{sec:interior}

\begin{figure}
\begin{center}
\includegraphics[width=60mm]{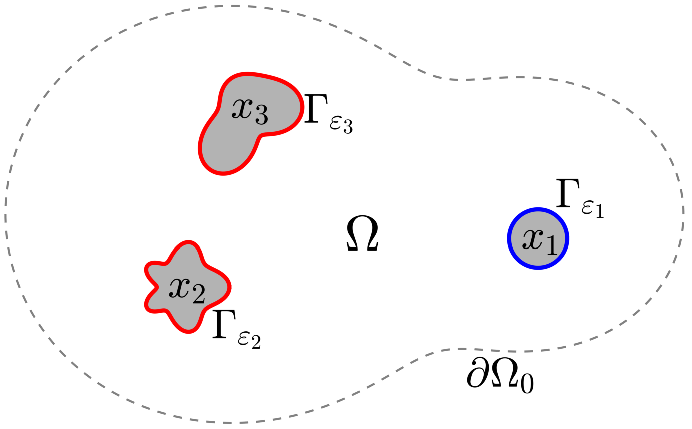} 
\end{center}
\caption{
Illustration of a bounded domain $\Omega \subset\R^2$ with a smooth
reflecting boundary $\pa_0$ (gray dashed line) and three interior
targets $\Omega_{\ve_j}$ (filled in gray), centered at $\x_j$, with
reactive boundaries $\Gamma_{\ve_j}$ (in red and blue).  For a
particle started from a point $\x\in\Omega$, the splitting probability
$S_1(\x)$ is the probability of hitting the boundary $\Gamma_{\ve_1}$
first.}
\label{fig:scheme2}
\end{figure}

Throughout this paper, we focused on reactive patches on the boundary
of a bounded domain.  In many applications, however, absorbing sinks,
traps and/or reactive targets can be hidden {\it inside} a bounded
domain $\Omega_0$, surrounded by a reflecting boundary $\pa_0$.  Let
us consider the problem with $N$ interior targets, where each target
is a compact set $\Omega_{\ve_j} \subset \Omega_0$ of size $\ve_j$,
centered at a point $\x_j$ (Fig. \ref{fig:scheme2}).  Since the
targets are impenetrable for a diffusing particle, we still consider
surface reactions on their boundaries, denoted as $\Gamma_{\ve_j} =
\pa_{\ve_j}$.
As before, the targets are small ($\ve_j \sim o(1)$), comparable in
size, and well-separated from each other: $|\x_j - \x_k| \sim
{\mathcal O}(1)$ for $j\ne k$, and from the domain boundary $\pa_0$:
$|\x_j - \x| \sim {\mathcal O}(1)$ for any $\x\in\pa_0$.  This setting
represents diffusion in a perforated domain $\Omega = \Omega_0
\backslash (\Omega_{\ve_1} \cup \cdots \cup \Omega_{\ve_N})$ with the
boundary $\pa = \pa_0 \cup \Gamma_{\ve_1} \cup \cdots \cup
\Gamma_{\ve_N}$.  These notations allow us to make an equivalence with
the earlier setting introduced at the beginning of
Sec. \ref{sec:Dirichlet}.  In particular, we can retain the same
formulations of the four considered problems.  As expected, their
solutions can be constructed analogously, but with some modifications
in the ``building blocks''. In this section, we briefly describe these
modifications and illustrate a few results.

\subsubsection*{Perfectly reactive targets: Dirichlet boundary condition}

If the $j$-th target is perfectly reactive (with Dirichlet condition),
the inner solution around this target is proportional to the exterior
Dirichlet Green's function
\begin{equation} 
\Delta g_\infty  = 0  \quad \textrm{in}~ \R^2 \backslash \Omega_{j}\,; \qquad
g_\infty  = 0 \quad \textrm{on}~ \pa_{j} \,; \qquad
g_\infty  \sim \ln |\y| + {\mathcal O}(1) \quad \textrm{as}~ |\y| \to \infty\,,
\end{equation}
where $\Omega_j = \ve_j^{-1} \Omega_{\ve_j}$ is the rescaled target
$\Omega_{\ve_j}$.  In contrast to the earlier studied case of a
Dirichlet patch, the Green's function $g_\infty(\y)$ depends on the shape
of the target $\Omega_j$ and thus is not universal.  In particular,
its asymptotic behavior is $g_\infty(\y) \sim \ln |\y| - \ln(d_j) + o(1)$
as $|\y|\to\infty$, where $d_j$ is the logarithmic capacity of
$\Omega_j$.  For instance, if $\Omega_j$ is the unit disk, its
logarithmic capacity is $1$.  Numerical values for $d_j$ for various
shapes of $\Omega_j$ are given in Table 1 of \cite{Kurella15}.

In addition, the outer solution involves the Neumann Green's function
$G_{\rm b}(\x,\xxi)$ (also known as pseudo-Green's function) that
satisfies
\begin{subequations} \label{eq:Neumann-Green}
\begin{align}  
  \Delta G_{\rm b} & = \frac{1}{|\Omega_0|}-\delta(\x-\xxi)  \quad \textrm{in}~
                     \Omega_0\,; \qquad
\partial_n G_{\rm b}  = 0 \quad \textrm{on}~ \pa_0 \,; \qquad
                     \int\limits_{\Omega_0} G_{\rm b}(\x,\xxi)\, d\x  = 0 \,,
                     \label{eq:Neumann-Green_eq}\\
  G_{\rm b}(\x,\xxi) & \sim - \frac{1}{2\pi} \ln|\x-\xxi| + R_{\rm b}(\xxi) + o(1)
  \quad \textrm{as}~ \x \to \xxi\in \Omega_0\,,\label{eq:G_asympt_bulk}
\end{align}
\end{subequations}
where we added the subscript $b$ to distinguish it from the surface
Neumann Green's function $G(\x,\xxi)$. The main difference between
this Green's function and the surface Neumann Green's function defined
by Eqs.  (\ref{eq:surfNeumann}) is that the singularity at $\xxi$ is
located in the bulk and not on the boundary (accordingly, there
is the factor $1/(2\pi)$ in Eq. (\ref{eq:G_asympt_bulk}) instead of
$1/\pi$).  For instance, when $\Omega_0$ is the unit disk, the Neumann
Green's function is well-known \cite{Kolokolnikov05}:
\begin{subequations}  \label{eq:NGreen_disk}
\begin{align} 
  G_{\rm b}(\x,\xxi) & = -\frac{1}{2\pi}\ln|\x-\xxi| - \frac{1}{4\pi}
                       \ln\left(|\x|^2|\xxi|^2 + 1 -2\x {\bf \cdot}\xxi \right)
                    + \frac{|\x|^2 + |\xxi|^2}{4\pi} - \frac{3}{8\pi}\,, \\
  R_{\rm b}(\xxi) & = -\frac{1}{2\pi}\ln\left(1-|\xxi|^2\right) +
              \frac{|\xxi|^2}{2\pi} - \frac{3}{8\pi} \,.
\end{align}
\end{subequations}
We remark that rapidly converging infinite series representations for
$G_{\rm b}$ and $R_{\rm b}$ are also known explicitly for ellipses
(see Eqs. (4.6, 4.7) from \cite{Iyaniwura21}) and for rectangles (see
Eq. (4.13) of Section 4.2 of \cite{kolok_split}, as well as
\cite{McCann01}).

With these minor changes in the ``building blocks'', we can repeat the
steps in Secs. \ref{sec:outer_solution} and \ref{sec:matrix} to obtain
that the splitting probability $S_k(\x)$ is
\begin{equation}  \label{eq:solution_bulk}
S_k(\x) = \chi_k - \sum\limits_{i=1}^N 2\pi A_i G_{\rm b}(\x,\x_i) \,,
\end{equation}
where the coefficients $\chi_k$ and $A_i$ are still determined
via Eqs. (\ref{eq:chiS}, \ref{eq:A_splitting}), but now with $\nu_j =
-1/\ln(\ve_jd_j)$.  Note also that the matrix $\G$ from
Eq. (\ref{eq:matrices_def}) is now based on the Neumann Green's
function $G_{\rm b}(\x,\xxi)$ from Eqs. (\ref{eq:Neumann-Green}) and
its regular part $R_{\rm b}(\xxi)$.  In addition, the factor $\pi$ in
the definition (\ref{eq:Gij}) of the elements of the matrix $\G$
should be replaced by $2\pi$.  The example of two targets can be
worked out explicitly.  Moreover, equally-spaced targets located on a
circular ring that is concentric within a unit disk can also be
treated explicitly.

\subsubsection*{Partially reactive targets: Robin boundary condition}

If the $j$-th target is {\em partially reactive}, we should replace the
Dirichlet boundary condition by a Robin condition with the reactivity
parameter $q_j$.  In the same vein, the former Robin Green's function
$g_\mu(\y)$ now satisfies
\begin{equation}
  \Delta g_{\mu}  = 0  \quad \textrm{in}~ \R^2 \backslash \Omega_{j}\,; \qquad
  \partial_n g_{\mu} + \mu g_{\mu} = 0 \quad \textrm{on}~ \pa_{j} \,; \qquad
  g_{\mu}  \sim \ln |\y| + \Cmu_j(\ve_j q_j) + o(1) \quad \textrm{as}~ |\y|
  \to \infty\,,
\end{equation}
with $\mu = \ve_j q_j$.  In particular, the constant term
$\Cmu_j(\mu)$ of its asymptotic behavior at infinity is not universal
and depends on the shape of $\Omega_j$.  We can still employ the
eigenmodes of the exterior Steklov problem in $\R^2 \backslash
\Omega_j$ to construct $g_\mu(\y)$ and to determine the spectral
expansion for $\Cmu_j(\mu)$:
\begin{equation}  \label{eq:Cj_bulk}
\Cmu_j(\mu) = -\ln(d_j) + \frac{2\pi}{\mu |\pa_j|} + 2\pi \sum\limits_{k=1}^\infty \frac{[\Psi_k^{j}(\infty)]^2}{\mu_k^{j} + \mu} \,,
\end{equation}
where $\mu_k^{j}$ and $\Psi_k^{j}$ are the eigenvalues and
eigenfunctions of the auxiliary exterior Steklov problem,
\begin{equation}
  \Delta \Psi_k^j  = 0 \quad \textrm{in}~ \R^2 \backslash \Omega_j\,; \qquad
\partial_n \Psi_k^j  = \mu_k^j \Psi_k^j \quad \textrm{on}~ \pa_j \,; \qquad
\Psi_k^j  \sim {\mathcal O}(1) \quad \textrm{as}~ |\y| \to \infty\,.
\end{equation}
We also used the normalization $\Psi_0^j = 1/\sqrt{|\pa_j|}$ of the
principal eigenfunction associated to $\mu_0^j = 0$.  As earlier, a
partially reactive target can be treated a perfect one, but with the
reduced size:
\begin{equation}
\ve_j^{\rm eff} = \ve_j \,\exp\bigl(-\ln(d_j) - \Cmu_j(\ve_j q_j)\bigr)\,,
\end{equation}
(see further discussion and examples in Sec. \ref{sec:Robin}).

For the target $\Omega_j$ of an arbitrary shape, the computation of
the Steklov eigenmodes $\mu_k^j$ and $\Psi_k^j$ requires numerical
techniques (e.g., a finite-element method, see
\cite{Chaigneau24,Grebenkov25b} and references therein).  However, if
the rescaled target $\Omega_j$ is the unit disk, the eigenmodes are
known explicitly and they all vanish at infinity, except $\Psi_0^j$.
As the logarithmic capacity of the unit disk is equal to $1$, we get a
particularly simple {\it exact} expression:
\begin{equation}  \label{eq:Cmu_disk}
\Cmu_{\rm disk}(\mu) = \frac{1}{\mu} \,.
\end{equation}
This is not surprising given that the exterior Green's function
for the unit disk is simply $g_\mu(\y) = 1/\mu + \ln |\y|$.

\subsubsection*{Imperfect targets}

In a similar way, we can easily reproduce the derivations of
Secs. \ref{sec:SteklovN} and \ref{sec:SteklovND} for the mixed
Steklov-Neumann and Steklov-Neumann-Dirichlet problems.  For an
imperfect target $\Omega_{\ve_j}$ of arbitrary shape, the main
difficulty is the lack of knowledge of the function $\Cmu_j(\mu)$,
which is formally accessible via the spectral expansion
(\ref{eq:Cj_bulk}) but its ``ingredients'' require numerical
computations.  Moreover, as the coefficients of the Taylor expansion
of $\Cmu_j(\mu)$ are unknown, we cannot rely on the approximation
(\ref{eq:Cmu_approx}).  As a consequence, many numerical steps would
be involved, and the analytical, almost explicit form of the
asymptotic results would in general be lost.  An interesting 
extension of this work consists of a systematic study of the
function $\Cmu_j(\mu)$ for targets of various shapes.

A drastic simplification appears when the imperfect targets are disks
due the explicit form (\ref{eq:Cmu_disk}) of the function $\Cmu(\mu)$.
In this case, the analysis of Secs. \ref{sec:SteklovN} and
\ref{sec:SteklovND} can be reproduced and will actually be even
simpler.  For instance, for the mixed Steklov-Neumann-Dirichlet
problem with a single disk-shaped target $\Omega_{\ve_1}$ and one
perfectly reactive target $\Omega_{\ve_2}$ (of arbitrary
shape), one can rewrite Eq. (\ref{eq:Cmu_1}) as
\begin{equation}  
  C = \ln(d_2 \ve_1 \ve_2) - 2\pi \bigl[R_{\rm b}(\x_1) + R_{\rm b}(\x_2)
  - 2G_{\rm b}(\x_1,\x_2)\bigr]\,,
\end{equation}
from which the principal eigenvalue reads
\begin{equation}  \label{eq:sigma0_disk_interior}
\frac{1}{\ve_1 \sigma_0} \approx - \ln(d_2 \ve_1 \ve_2) + 2\pi \bigl[R_{\rm b}(\x_1) + R_{\rm b}(\x_2) - 2G_{\rm b}(\x_1,\x_2)\bigr]\,.
\end{equation}

\begin{figure}
\begin{center}
\includegraphics[width=88mm]{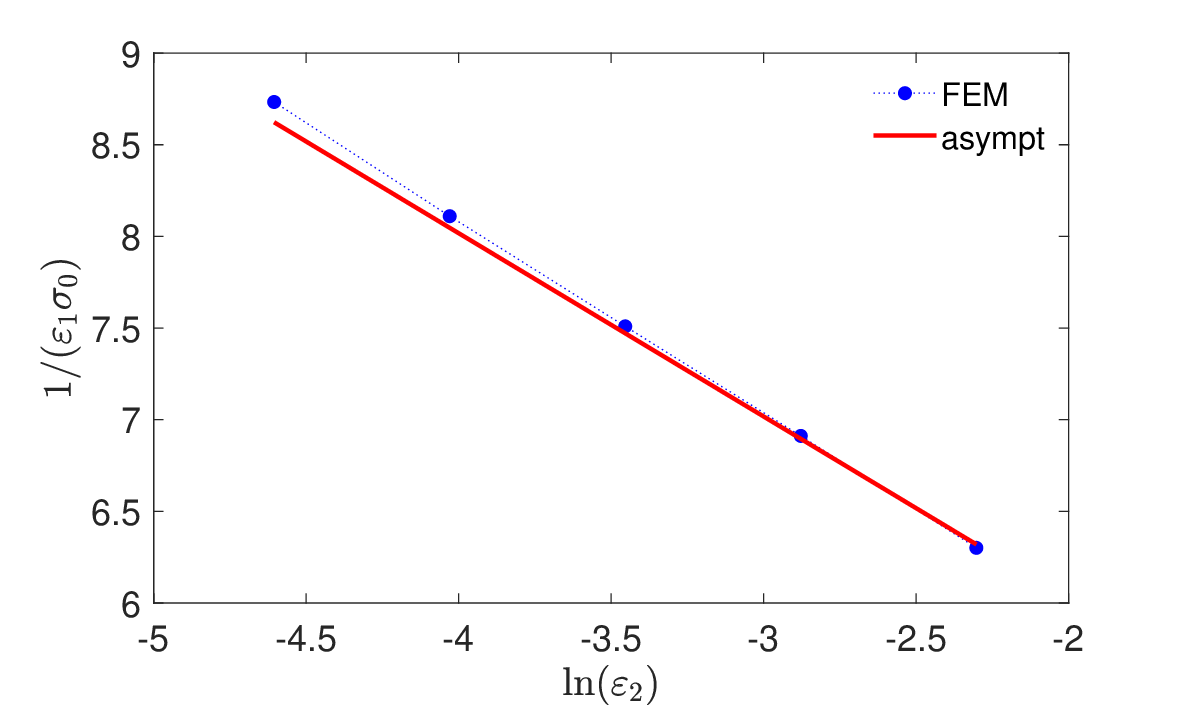} 
\end{center}
\caption{
The asymptotic behavior of the principal eigenvalue of the mixed
Steklov-Neumann-Dirichlet problem, plotted as $1/(\ve_1\sigma_0)$
versus $\ve_2$, for the unit disk with two interior circular targets
of radii $\ve_1 = 0.05$ and $\ve_2$ (variable from $0.01$ to $0.1$),
located at $\x_1 = (-0.5,0)$ and $\x_2 = (0.5,0)$.  Filled circles
present the numerical solution by a FEM with the maximal meshsize of
$0.005$, solid line indicates Eq. (\ref{eq:sigma0_disk_interior}).}
\label{fig:disk_interior}
\end{figure}

Figure \ref{fig:disk_interior} illustrates the accuracy of the
asymptotic relation (\ref{eq:sigma0_disk_interior}) for the unit disk
with two interior circular targets.  We observe a close agreement
between a numerical solution and the asymptotic formula.  One can
notice a small deviation between two lines that slightly increases as
the target radius $\ve_2$ decreases.  This minor discrepancy seems to
be a numerical artefact due to the available meshsize $0.005$, which
becomes comparable to the target at small $\ve_2$.  To check this
point, we computed the eigenvalue of the mixed Steklov-Dirichlet
problem for a circular annulus with radii $\ve_2$ and $R = 1$, with
Steklov condition on the outer circle and Dirichlet condition on the
inner circle.  As the exact solution of this problem is known,
$\sigma_0 = 1/(R\ln(R/\ve_2))$, we could compare it with the numerical
results, and found the same minor discrepancy.

In summary, we conclude that interior targets can be handled in
essentially the same way as boundary patches, even though the
asymptotic analysis becomes sensitive to the shapes of the targets.
Moreover, one can combine interior targets with boundary patches that
opens a way to access a broad variety of various geometric settings.

\subsection{Exterior problems}
\label{sec:exterior}

Another extension of the present approach is related to
exterior problems in $\Omega = \R^2 \backslash \Omega_0$, where
$\Omega_0$ is a simply-connected compact domain.  While the inner 
solutions remain unchanged, the outer solution is now constructed using the
exterior surface Neumann Green's function, labeled by $G_{\rm e}$, which
satisfies
\begin{subequations} \label{eq:Neumann-Green_ext}
\begin{align}
\Delta G_{\rm e} & = 0  \quad \textrm{in}~ \R^2\backslash \Omega_0\,,  \\  
  G_{\rm e}(\x,\xxi) & \sim - \frac{1}{\pi} \ln|\x-\xxi| + R_{\rm e}(\xxi) + o(1)
                       \quad \textrm{as}~  \x \to \xxi\in \pa_0\,;  \qquad
 \partial_n G_{\rm e}  = 0 \quad \textrm{on}~ \pa_0 \backslash \{\xxi\} \,,
                       \label{eq:ext_bc}\\
  G_{\rm e} (\x,\xxi) & \sim -\frac{1}{2\pi}\ln |\x| + o(1) \quad
                        \textrm{as}~ |\x|\to \infty\,, \label{eq:ext_ff}
\end{align}
\end{subequations}
where $R_{\rm e}(\xxi)$ is the regular part of $G_{\rm
e}(\x,\xxi)$ at $\xxi$.  The condition that $G_{\rm e}(\x,\xxi)
+(2\pi)^{-1}\ln|\x|\to 0$ as $|\x|\to \infty$ determines $G_{\rm e}$
uniquely.

In the case when $\Omega_0$ is the unit disk and $|\xxi|=1$, we
claim that
\begin{equation}     \label{eq:Green_disk_ext1}
  G_{\rm e}(\x,\xxi) = -\frac{1}{\pi}\ln|\x-\xxi| + \frac{1}{2\pi} \ln|\x| \,,
  \qquad   R_{\rm e}(\xxi) = 0 \,.
\end{equation}
Clearly Eq. (\ref{eq:Green_disk_ext1}) satisfies the behavior
(\ref{eq:ext_ff}) as $|\x|\to\infty$ as well as Eq. (\ref{eq:ext_bc})
as $\x\to\xxi$, where we identify that $R_{\rm e}(\xxi) = 0$ (see more
details in Appendix \ref{sec:Gext_disk}).

For instance, if there are two patches on the unit circle, the
solution of any of four earlier considered problems involves
\begin{equation}
  R_{\rm e}(\x_1) + R_{\rm e}(\x_2) - 2G_{\rm e}(\x_1,\x_2) =
  \frac{2}{\pi} \ln|\x_1-\x_2|\,,
\end{equation}
which is identical for both interior and exterior domains.  This
property is consistent with the fact that the eigenvalues of interior
and exterior mixed Steklov problems for the unit disk are identical.
In contrast, the associated eigenfunctions behave differently.

\section{Discussion}\label{sec:conclusion}

In this paper, we have established a general mathematical framework
for studying the competition of small targets for a diffusing particle
in planar domains.  Using the method of matched asymptotic expansions,
as tailored for problems with localized defects \cite{Ward93} and with
logarithmic gauge functions \cite{Ward93b}, we solved four different
problems of increasing complexity in the boundary conditions: (i)
splitting probabilities for perfectly reactive patches with Dirichlet
condition; (ii) their extension to partially reactive patches
with Robin condition; (iii) mixed Steklov-Neumann problem describing
imperfect patches; and (iv) mixed Steklov-Neumann-Dirichlet problem
describing the escape of a particle through Dirichlet patches in the
presence of an imperfect patch.  Although the first problem was
thoroughly studied in the past, we have improved some former results.
To our knowledge, the asymptotic behavior for the three other problems
in the small-target limit has not been reported previously.  Moreover,
we discussed two further extensions of our results to the case of
interior targets and to exterior problems.  The established asymptotic
formalism in our 2-D setting can be applied to a broad variety of
natural and industrial phenomena such as diffusion-controlled
reactions in chemistry and biology.

It would be worthwhile to extend our analytical framework to
determine high-order asymptotic expansions to treat analogous 3-D
problems with many either partially reactive or imperfect (Steklov)
patches on the domain boundary. For a locally circular partially
reactive patch on the boundary of a 3-D domain, a leading-order
asymptotic theory was derived in \cite{Cengiz24} to determine the mean
first-passage time for small, intermediate, and large patch
reactivities. In \cite{Guerin23}, the large reactivity limit was
analyzed in detail. However, it is an open problem to derive
high-order asymptotic expansions allowing for multiple partially
reactive or imperfect patches, as the local geometry of the domain
boundary will play a key role in the analysis.

\begin{acknowledgments}
The authors thank professors I. Polterovich and M. Levitin for
fruitful discussions, and A. Chaigneau for his implementation of
the FEM code.  D.S.G. acknowledges the Simons Foundation for
supporting his sabbatical sojourn in 2024 at the CRM, University of
Montr\'eal, Canada, and the Alexander von Humboldt Foundation for
support within a Bessel Prize award.  M.J.W. was supported by the
NSERC Discovery grant program.
\end{acknowledgments}

\subsection*{Competing interests declaration}

Competing interests: The author(s) declare none

\appendix
\section{Green's function for the Dirichlet patch in the half-plane}
\label{sec:wc}

The Green's function satisfying Eq. (\ref{eq:w_problem}) can be found
exactly.  Even though this solution is classical \cite{Saff} we
reproduce it here for completeness.
For this purpose, we use the elliptic coordinates for an ellipse with
semi-axes $a > b$:
\begin{equation}  \label{eq:ellip_coord}
y_1 = a_E \cosh\alpha \cos\theta\,, \qquad y_2 = a_E \sinh\alpha \sin\theta \,,
\end{equation}
where $a_E = \sqrt{a^2-b^2}$, $\alpha\ge 0$ and $-\pi < \theta \leq
\pi$.  It is worth noting that all points on the horizontal
interval $(-a_E,a_E)\times \{0\}$ correspond to $\alpha = \theta = 0$
and are thus indistinguishable.

In our setting, we fix $a = 1$ and $b = 0$.  We search for
$g_\infty(\y)$ in the form
\begin{equation}
  g_\infty(\y) = \ln|\y| - \ln(d) - \sum\limits_{n=1}^\infty c_n
  \cos(n\theta) e^{-n\alpha} \,,
\end{equation}
with unknown coefficients $c_n$.  This is a general form of a harmonic
function, which behaves as $\ln|\y| - \ln(d)$ at infinity and
satisfies the condition $\partial_n g_\infty = 0$ on $|y_1|>1$, $y_2 =0$.
The coefficients $c_n$ are determined by the condition $g_\infty = 0$ on
$|y_1|\leq 1$, $y_2=0$, which yields
\begin{equation} \label{eq:wc_BC}
  0 = \ln |\cos\theta| - \ln(d) - \sum\limits_{n=1}^\infty c_n \cos(n\theta)\,,
  \qquad (0 \leq \theta \leq \pi)
\end{equation}
(here we restrict the analysis to the upper half-plane, with $\theta
\geq 0$).  Using the expansion
\begin{equation}
\ln(2|z-z_0|) = - \sum\limits_{n=1}^\infty \frac{2}{n} T_n(z) T_n(z_0)\,,
\end{equation}
where $T_n(\cos z) = \cos(nz)$ are the Chebyshev polynomials and
$z_0=0$, we immediately see that the boundary condition
(\ref{eq:wc_BC}) implies
\begin{equation}
\ln(d) = - \ln(2)\,, \qquad c_n = \frac{2}{n} \,.
\end{equation}
The solution then reads 
\begin{equation}  \label{eq:gD_interval}
  g_\infty(\y) = \ln|\y| + \ln(2) - \sum\limits_{n=1}^\infty \frac{2}{n}
  \cos(n\theta) e^{-n\alpha} \,.
\end{equation}
By summing this series in terms of the logarithm we get
\begin{equation} \label{eq:wc_explicit}
  g_\infty(\y) = \ln|\y| + \ln(2) - \ln\bigl(1 - 2\cos\theta e^{-\alpha} +
  e^{-2\alpha} \bigr) \,.
\end{equation}
Note that the elliptic coordinates $\alpha$ and $\theta$ can be easily
expressed in terms of $\y = (y_1,y_2)$ by setting
\begin{equation}  \label{eq:rpm}
r_{\pm} = \sqrt{(y_1\pm a_E)^2 + y_2^2} = a_E (\cosh\alpha \pm \cos\theta)\,,
\end{equation}
from which
\begin{equation}  \label{eq:elliptic_coord}
\cosh\alpha = \frac{r_+ + r_-}{2a_E} \,,   \qquad
\cos\theta =  \frac{r_+ - r_-}{2a_E} \,.
\end{equation} 

Using Eq. (\ref{eq:gD_interval}), we calculate that
\begin{equation}  \label{eq:dnwc_explicit}
  - \partial_n g_\infty|_{(-1,1)} = \biggl(\frac{1}{h_{\alpha}}
  \partial_\alpha g_\infty\biggr)_{\alpha=0}
= \frac{1}{|\sin\theta|} = \frac{1}{\sqrt{1-y_1^2}} \,,
\end{equation}
where we used
$h_\alpha = a_E \sqrt{\cosh^2\alpha - \cos^2\theta} = |\sin\theta|$ at
$\alpha = 0$ for the scale factor.  In particular, the integral of
this expression over the interval $(-1,1)$ is equal to $\pi$, as
expected from the divergence theorem.

\section{Limit of many small targets}
\label{sec:kappa}

In this Appendix, we study the large-$N$ behavior of the eigenvalues
$\kappa_j$, given in Eq. (\ref{eq:kappaj_new}), of the matrix $\G$ for
$N$ identical equally-spaced patches on the boundary of the unit disk.
Since the $\kappa_j$ are the eigenvalues of the symmetric and
circulant matrix $\G$, we have $\kappa_j=\kappa_{N-j}$ for
$j=1,\ldots,{N/2}$ when $N$ is even. As such, we need only 
estimate $\kappa_j$ for $j=1,\ldots,{N/2}$ when $N\gg 1$ is even.

To do so, we use the Euler-Maclaurin expansion for a $C^{\infty}$
function $f(\theta)$ on $1\leq \theta\leq N-1$, which is given by
\begin{equation}
  \sum_{m=1}^{N-1} f(m) = \int_{1}^{N-1} f(\theta)\, d\theta  +
  \frac{1}{2}\left[ f(N-1)+f(1)\right] + \frac{1}{12}\left[
    f^{\prime}(N-1)-f^{\prime}(1)\right] + \cdots\,, \label{eqb:em}
\end{equation}
where from Eq. (\ref{eq:kappaj_new}) we define $f(\theta)$ by
\begin{equation}\label{eqb:fdef}
  f(\theta)=\cos(2j a\theta) \ln\left[\sin(a\theta)\right] \,, \quad
   a=\frac{\pi}{N}\,.
\end{equation}
From Eq. (\ref{eq:kappaj_new}) we identify
\begin{equation}\label{eqb:em_kappa}
  \kappa_j = \ln{2} - \sum_{m=1}^{N-1} f(m) \,, \quad j\in \lbrace{1,\ldots,
    {N/2}\rbrace} .
\end{equation}
We first estimate the integral in Eq. (\ref{eqb:em}), labeled by
$I=\int_{1}^{N-1} f(\theta)\, d\theta$.  Upon substituting
$x=a\theta$, we get
\begin{equation}
  I = \frac{N}{\pi} \int_{0}^{\pi}
      \cos(2j x)\ln(\sin{x})\, dx - \frac{2N}{\pi}\int_{0}^{\pi/N}
  \cos(2j x)\ln(\sin{x})\, dx\,. \label{eqb:int_1}
\end{equation}
The first integral on the right-side of Eq. (\ref{eqb:int_1}) can be
evaluated explicitly as ${-\pi/(2j)}$, whereas in the second integral
we use $\sin(x)\approx x$ on the range $0<x<{\pi/N}$, which is valid
for $N\gg 1$.  In this way, for $N\gg 1$ we obtain
\begin{equation}\label{eqb:int_2}
  I \sim -\frac{N}{2j} - \frac{2N}{\pi}
  \int_{0}^{\pi/N} \cos(2jx)\ln{x}\, dx \,.
\end{equation}
Upon integrating by parts in Eq. (\ref{eqb:int_2}), we find that
\begin{equation}\label{eqb:int_3}
  I \sim -\frac{N}{2j} - \frac{N}{\pi j}
  \sin\left(\frac{2\pi j}{N}\right) \ln\left(\frac{\pi}{N}\right)
  +\frac{N}{\pi j} \mbox{Si}\left(\frac{2\pi j}{N}\right)\,,
\end{equation}
where $\mbox{Si}(x)=\int_{0}^{x}\xi^{-1}\sin\xi \,d\xi$ is the 
sine integral function.  Moreover, we readily calculate for $N\gg 1$
that
\begin{subequations}\label{eqb:endp}
\begin{align}
  f(1)=f(N-1) &\sim \cos\left(\frac{2\pi j}{N}
                \right) \ln\left(\frac{\pi}{N}\right) \,,  \\
  f^{\prime}(1)=-f^{\prime}(N-1) &\sim -\frac{2\pi j}{N}
     \sin\left(\frac{2\pi j}{N}\right)\ln\left(\frac{\pi}{N}\right) +
    \cos\left(\frac{2\pi j}{N}\right) \,.
\end{align}
\end{subequations}
Upon substituting Eqs. (\ref{eqb:int_3}, \ref{eqb:endp}) into Eq.
(\ref{eqb:em}), and recalling Eq. (\ref{eqb:em_kappa}), we conclude
for $N\gg 1$ and for $j=1,\ldots,{N/2}$ that
\begin{equation} \label{eqb:kappa_final}
  \kappa_j \sim \frac{N}{2j} + \ln\left(\frac{\pi}{N}\right) {\mathcal
             B}\left(\frac{j}{N}\right) + \ln{2} + \frac{1}{6}
           \cos\left(\frac{2\pi j}{N}\right) - \frac{N}{\pi j}
             \mbox{Si}\left(\frac{2\pi j}{N}\right)\,,
\end{equation}
where ${\mathcal B}(\xi)$, with $\xi={j/N}$, is defined by
\begin{equation}
  {\mathcal B}(\xi) = \frac{1}{\pi \xi} \sin(2\pi \xi) - \cos(2\pi \xi)
  - \frac{\pi \xi}{3} \sin(2\pi\xi)\,.
\end{equation}
We calculate from a Maclaurin series that ${\mathcal B}(\xi)= 1+
{(2\pi \xi)^4/360}+{\mathcal O}(\xi^6)$, and so we will approximate
${\mathcal B}(\xi)\approx 1$ on $0<\xi<{1/2}$.  We then write
Eq. (\ref{eqb:kappa_final}) as
\begin{subequations}\label{eqb:kappa_final2}
\begin{equation}
  \kappa_j \sim \frac{N}{2j} + \ln\left( \frac{2\pi e^{-11/6}}{N}\right)
  + {\mathcal D}\left( \frac{j}{N}\right) \,,
\end{equation}
where ${\mathcal D}(\xi)$, with ${\mathcal D}(0)=0$, is defined by
\begin{equation}
  {\mathcal D}(\xi) = \frac{11}{6} + \frac{1}{6}\cos(2\pi \xi) -
  \frac{1}{\pi \xi} \mbox{Si}\left(2\pi \xi\right)\,.
\end{equation}
\end{subequations}
By using $\cos{z}\sim 1-{z^2/2}$ and $\mbox{Si}(z)\sim z-{z^3/18}$,
the Maclaurin series for ${\mathcal D}(\xi)$ is ${\mathcal
D}(\xi)={\pi^2\xi^2/9} + {\mathcal O}(\xi^4)$.  From this lowest-order
approximation, Eq.~(\ref{eqb:kappa_final2}) becomes
\begin{equation}\label{eqb:kappa_final3}
  \kappa_j \sim \frac{N}{2j} + \ln\left( \frac{2\pi e^{-11/6}}{N}\right)
  + \frac{\pi^2 j^2}{9N^2} \,,
\end{equation}
which should be rather accurate if $j\ll {N/2}$.  Neglecting the
correction term in Eq.~(\ref{eqb:kappa_final3}) gives a simpler, but
less accurate, approximation
\begin{equation}\label{eqb:kappa_approx}
  \kappa_j \sim \frac{N}{2j} + \ln\left( \frac{2\pi e^{-11/6}}{N}\right)\,.
\end{equation}
Table \ref{tab:kappaj} illustrates the accuracy of the three
approximate relations (\ref{eqb:kappa_final}, \ref{eqb:kappa_final3},
\ref{eqb:kappa_approx}) for two cases: $N = 16$ and $N = 64$.  Even
for a moderate number of patches ($N = 16$), these three relations
approximate $\kappa_1$ very accurately.  As the index $j$ increases,
the accuracy of both relations expectedly reduces but remains
good. The accuracy is even higher when $N = 64$.


\begin{table}
\begin{center}
\begin{tabular}{|c|c|c|c|c|c|c|c|c|c|}  \hline
&  $j$  & 1  &  2  &  3  & 4 & 5 & 6 & 7 & 8  \\  \hline  
\multirow{3}{3mm}{\begin{turn}{90} $N = 16$~ \end{turn}}
& Eq. (\ref{eq:kappaj_new})    & 5.2321 & 1.2465 & -0.0623 & -0.6931 & -1.0443 & -1.2465 & -1.3529 & -1.3863 \\
  & Eq. (\ref{eqb:kappa_final})  & 5.2362 & 1.2489 & -0.0634 & -0.6986 & -1.0514 & -1.2442 & -1.3179 & -1.2804 \\
 & Eq. (\ref{eqb:kappa_final3}) & 5.2320 & 1.2491 & -0.0628 & -0.6995 & -1.0610 & -1.2805 & -1.4153 & -1.2740 \\
  & Eq. (\ref{eqb:kappa_approx}) & 5.2362 & 1.2320 & -0.1014 & -0.7680 & -1.1680 & -1.4347 & -1.4939 & -1.7680 \\ \hline \hline
\multirow{3}{3mm}{\begin{turn}{90} $N = 64$~ \end{turn}}
& Eq. (\ref{eq:kappaj_new})    & 27.8414 & 11.8423 & 6.5104 & 3.8458 & 2.2485  & 1.1851  &  0.4271 &  -0.1398 \\
  & Eq. (\ref{eqb:kappa_final})  & 27.8459 & 11.8467 & 6.5147 & 3.8498 & 2.2521  & 1.1881  &  0.4293 &  -0.1388 \\
 & Eq. (\ref{eqb:kappa_final3}) & 27.846 & 11.8470 & 6.5147 & 3.8499 & 2.2524 & 1.1886 & 0.43021 & -0.1372 \\
& Eq. (\ref{eqb:kappa_approx}) & 27.8457 & 11.8457 & 6.5123 & 3.8457 & 2.2457  & 1.1790  &  0.4171 &  -0.1543 \\  \hline
\end{tabular}
\end{center}
\caption{
Comparison between the exact values of $\kappa_j$ from
Eq. (\ref{eq:kappaj_new}), their approximation
(\ref{eqb:kappa_final}), and its simpler asymptotic forms
(\ref{eqb:kappa_final3}) and (\ref{eqb:kappa_approx}).  }
\label{tab:kappaj}
\end{table}

\begin{figure}
\begin{center}
\includegraphics[width=88mm]{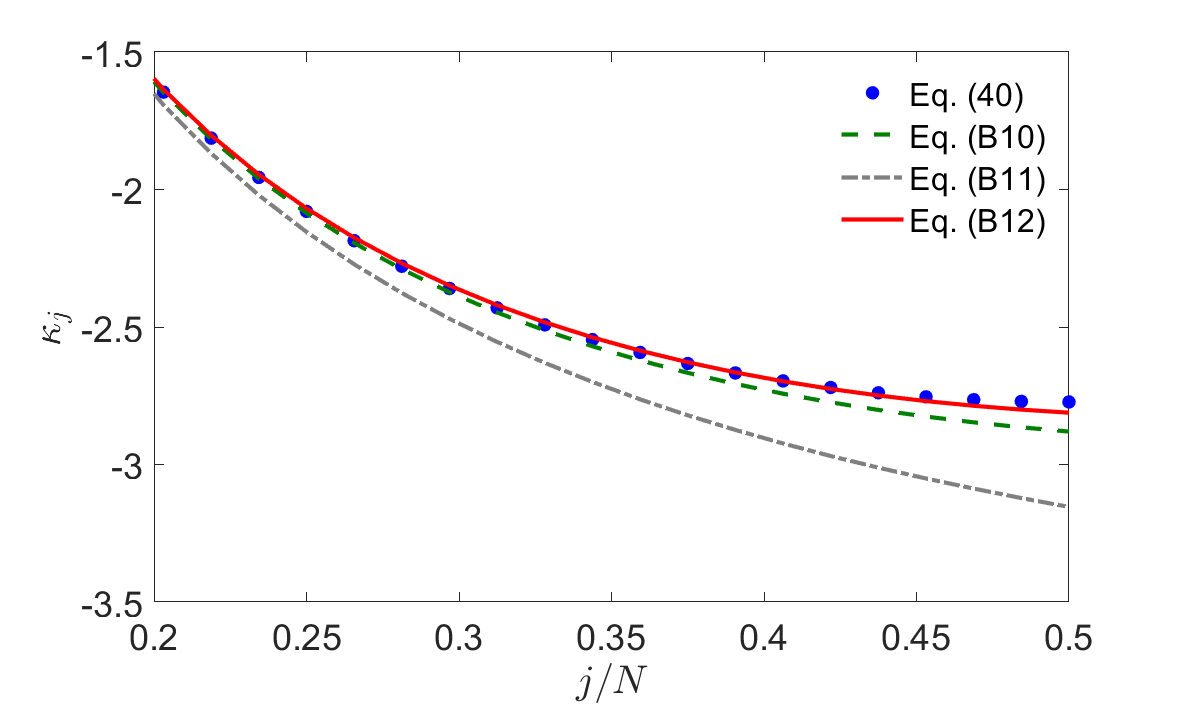} 
\end{center}
\caption{
Exact values of $\kappa_j$ versus ${j/N}$ on $0.2 <{j/N}< 0.5$ for
$N=64$ from Eq. (\ref{eq:kappaj_new}), shown by filled circles, and
their approximations (\ref{eqb:kappa_final3}, \ref{eqb:kappa_approx},
\ref{eqb:kappa_final4}) shown by lines.}  
\label{fig:kappa_approx}
\end{figure}

In turn, Fig. \ref{fig:kappa_approx} illustrates the accuracy of the
approximations (\ref{eqb:kappa_final3}) and (\ref{eqb:kappa_approx})
on a broader range ${1/5}< {j/N}<{1/2}$, for $N = 64$.  While both
approximations are accurate at small $j/N$, one can still observe
deviations for $j/N$ around $1/2$.  These deviations have (at least)
two origins: (i) neglection of higher-order terms
${\mathcal O}(\xi^4)$, and (ii) omission of the higher-order
derivatives in the Euler-Maclaurin expansion (\ref{eqb:em}).  A
careful examination of the next-order term in this expansion,
$-\tfrac{1}{720} (f^{\prime\prime\prime}(N-1)-
f^{\prime\prime\prime}(1))$, reveals that it yields the contribution
$\tfrac{1}{90} \pi^2 j^2/N^2$ to $\kappa_j$ that increases by $10\%$
the last term in Eq. (\ref{eqb:kappa_final3}).  Skipping a systematic
analysis of these higher-order contributions, we adjust the
coefficient in front of this term to get an empirical approximation:
%
%
\begin{equation}\label{eqb:kappa_final4}
  \kappa_j \sim \frac{N}{2j} + \ln\left( \frac{2\pi e^{-11/6}}{N}\right)
  + 1.25 \frac{\pi^2 j^2}{9N^2} \,.
\end{equation}
In this way, we achieve a decent approximation over the whole range of
$j/N$, as illustrated by the solid line in
Fig. \ref{fig:kappa_approx}.

\section{Green's function for the Robin patch in the half-plane}
\label{sec:Cmu_asympt}

In this Appendix, we obtain the exact form of the  Robin Green's
function satisfying Eq. (\ref{eq:Steklov_inner}) in the upper
half-plane.  We will determine $g_\mu(\y)$ in the form
\begin{equation} \label{eq:ws_spectral}
g_\mu(\y) = g_\infty(\y) + \sum\limits_{k=0}^\infty c_k \, \Psi_k(\y)\,,
\end{equation}
where $g_\infty(\y)$ is the Dirichlet Green's function satisfying
(\ref{eq:w_problem}) with Dirichlet condition on the interval
$(-1,1)$, $c_k$ are unknown coefficients, and $\Psi_k(\y)$ are the
eigenfunctions of the auxiliary Steklov-Neumann problem
(\ref{eq:Vkint}) in the upper-half plane $\H_2$, with $\mu_k$ being
the associated eigenvalues.  We recall that this spectral problem has
infinitely many solutions, enumerated by the index $k = 0,1,2,\ldots$,
with the nonnegative eigenvalues increasing up to infinity, $0 = \mu_0
\leq \mu_1 \leq \mu_2 \leq \ldots \nearrow +\infty$, while the
restrictions of Steklov eigenfunctions $\Psi_k(\y)$ onto the interval
$(-1,1)$ of the horizontal axis form a complete orthonormal basis of
$L^2(-1,1)$:
\begin{equation} \label{eq:Vk_orthonormal}
\int\limits_{-1}^1 \Psi_j(y_1,0) \, \Psi_k(y_1,0) \, dy_1 = \delta_{j,k} \,.
\end{equation}
Note that this condition fixes the normalization of the Steklov
eigenfunctions $\Psi_k$.  A rigorous formulation of the exterior
Steklov problem in the plane is discussed in \cite{Bundrock25} 
(see also \cite{Christiansen23}), whereas a numerical construction of
the eigenfunctions in elliptic coordinates, described in
\cite{Grebenkov25}, is summarized in Appendix
\ref{sec:Steklov_interval}.

We substitute Eq. (\ref{eq:ws_spectral}) into
Eq. (\ref{eq:Steklov_inner_BC}) to get
\begin{equation*}
  (\partial_n g_\infty)\bigr|_{y_2=0} + \sum\limits_{k=0}^\infty
  c_k (\mu_k + \mu) \Psi_k(y_1,0) = 0\,,  \qquad |y_1|< 1\,,
\end{equation*}
where we used $g_\infty|_{(-1,1)} = 0$.  Multiplying this equation by
$\Psi_j$, integrating over the interval $(-1,1)$ on the horizontal
axis, and using the orthonormality (\ref{eq:Vk_orthonormal}), we
identify the coefficients $c_k$ and thus the Green's function as
\begin{equation}  \label{eq:ws_wc}
  g_\mu(\y) = g_\infty(\y) + \sum\limits_{k=0}^\infty \frac{b_k}{\mu + \mu_k}
  \Psi_k(\y)\,,
\end{equation}
where
\begin{equation}  \label{eq:bk_def}
b_k = \int\limits_{-1}^1 \Psi_k(y_1,0) \, (-\partial_n g_\infty) \, dy_1 \,.
\end{equation}
Note that explicit formulas for $g_\infty(\y)$ and $\partial_n
g_\infty$ are given by Eqs. (\ref{eq:wc_explicit},
\ref{eq:dnwc_explicit}).  However, we can actually get the
coefficients $b_k$ without computing the integral in
Eq. (\ref{eq:bk_def}).  For this purpose, Eq. (\ref{eq:Vkint_Laplace})
is multiplied by $g_\infty(\y)$, Eq. (\ref{eq:wc_Laplace}) is
multiplied by $\Psi_k(\y)$, they are subtracted from each other, and
integrated over the upper half-plane, to get using Green's second
identity that
\begin{equation*}
  0  = \int\limits_{\H_2} \bigl(g_\infty \Delta \Psi_k - \Psi_k \Delta
  g_\infty\bigr) \, d\y = - \pi \Psi_k(\infty)  - \int\limits_{-1}^1 \Psi_k(y_1,0)
  (\partial_n g_\infty) \, dy_1\,,
\end{equation*}
from which we identify that
\begin{equation}
b_k = \pi\, \Psi_k(\infty)\,.
\end{equation}
Moreover, the symmetry of the problem (\ref{eq:Vkint}) with respect to
the vertical axis implies that the eigenfunctions $\Psi_k$ should be
either symmetric or antisymmetric:
\begin{subequations}
\begin{align}
\Psi_{2k}(-y_1,y_2) & = \Psi_{2k}(y_1,y_2)\,,   \\
\Psi_{2k+1}(-y_1,y_2) & = - \Psi_{2k+1}(y_1,y_2)\,,
\end{align}
\end{subequations}
(even and odd indices are used to distinguish them).  Since the
function $(\partial_n g_\infty)$ is symmetric, the integrals in
Eq. (\ref{eq:bk_def}) are zero for odd indices, implying
\begin{equation}
b_{2k+1} = 0\,.
\end{equation}
We conclude that
\begin{equation}  \label{eq:gR_bis}
  g_\mu(\y) = g_\infty(\y) + \pi \sum\limits_{k=0}^\infty
  \frac{\Psi_{2k}(\infty)}{\mu + \mu_{2k}} \Psi_{2k}(\y)\,.
\end{equation}

As $g_\mu(\y)$ is the Robin Green's function of the Laplace equation,
it is necessarily positive for any $\y$ and $\mu > 0$ \cite[Chapter
V.1]{Bergman}.  Moreover, we observed numerically the following
property: there exists $\hat{\mu} < 0$ such that, for any $\hat{\mu} <
\mu < 0$, the restriction of $g_\mu(\y)$ onto the interval $(-1,1)$ is
negative:
\begin{equation}  \label{eq:gR_negative}
g_\mu(y_1,0) < 0    \quad \textrm{for any} ~ -1 \leq y_1 \leq 1\,.
\end{equation}  
We obtained numerically that $\hat{\mu} \approx -2.006$, which is very
close to and possibly identical with $-\mu_1$.  Qualitatively, when
$\mu$ is negative but small, the first term of the sum, $\pi/(2\mu)$,
provides the dominant (negative) contribution to $g_\mu(y_1,0)$, as
compared to the remaining terms whose sum is expected to be bounded by
a constant.  However, we are not aware of the proof of this statement.

According to the spectral expansion (\ref{eq:gR_bis}), the constant
term $\Cmu(\mu)$ of the Robin Green's function $g_\mu(\y)$ at
infinity, as defined in Eq. (\ref{eq:ws_asympt}), is
\begin{equation}  \label{eq:Cmu_spectral_bis}
  \Cmu(\mu) = \ln (2) + \frac{\pi}{2\mu} + \pi \sum\limits_{k=1}^\infty
  \frac{[\Psi_{2k}(\infty)]^2}{\mu_{2k} + \mu} \,,
\end{equation}
where we used Eq. (\ref{eq:wc_infinity}) with $d = 1/2$ and wrote
explicitly the term with $k = 0$, for which $\mu_0 = 0$ and $\Psi_0 =
1/\sqrt{2}$ that yielded $\pi/(2\mu)$.  The numerical eigenvalues
$\mu_{2k}$ and the coefficients $[\Psi_{2k}(\infty)]^2$ for the first
ten terms are reported in Table \ref{tab:mu_int}.  The asymptotic
behavior of the eigenvalues is well known (see
\cite{Grebenkov25,Polosin22} and references therein):
\begin{equation}  \label{eq:muk_asympt}
\mu_k \sim \frac{\pi}{2} k \qquad (k \gg 1)\,.
\end{equation}
In turn, the oscillating eigenfunction $\Psi_{2k}(y_1)$ can be roughly
approximated as $\cos(\pi ky_1)$ at large $k$ (see Appendix
\ref{sec:Steklov_interval}).  As a consequence, we get
\begin{equation} \label{eq:Vkinf_asympt}
  \Psi_{2k}(\infty)  = \frac{1}{\pi} \int\limits_{-1}^1 \frac{\Psi_{2k}(y_1)}
  {\sqrt{1-y_1^2}} \, dy_1  
  \approx \frac{1}{\pi} \int\limits_{-1}^1 \frac{\cos (k\pi y_1)}{\sqrt{1-y_1^2}}
  \, dy_1 = J_0(\pi k) \simeq \frac{(-1)^k}{\pi \sqrt{k}} \,, \qquad (k\gg 1)\,.
\end{equation}
The decay of $[\Psi_{2k}(\infty)]^2/\mu_{2k} \propto 1/k^2$ is rapid
enough to ensure that the reported ten coefficients are sufficient for
an accurate approximation of $\Cmu(\mu)$, at least for small $\mu$.

To get the small-$\mu$ approximation, we expand the last term
of Eq. (\ref{eq:Cmu_spectral_bis}) into a Taylor series in powers of
$\mu$ as
\begin{equation} 
\Cmu(\mu) = \frac{\pi}{2\mu} + \sum\limits_{n=0}^\infty (-\mu)^n C_{n+1} \,,
\end{equation}
with
\begin{equation}  \label{eq:Cn_sum}
  C_n = \delta_{n,1} \ln 2 + \pi \sum\limits_{k=1}^\infty
  \frac{[\Psi_{2k}(\infty)]^2}{[\mu_{2k}]^{n}}\,, \qquad (n = 1,2,\ldots)\,.
\end{equation}
Substituting the first ten contributing terms from Table
\ref{tab:mu_int}, we get $C_1 \approx 0.7976$ and $C_2 \approx
0.0222$.  In Appendix \ref{sec:Cmu}, we provide an exact computation
of these coefficients that yields
\begin{equation}
  C_1 = 3/2 - \ln 2 \approx 0.8069\,, \qquad
  C_2 = \frac{21-2\pi^2}{18\pi} \approx 0.0223\,.
\end{equation} 
One sees that the numerically computed values are very close to the
exact ones.  Most importantly, the coefficient $C_2$, as well as
higher-order coefficients, are small and can thus be neglected when
$\mu \ll 1$.

\begin{table}
\begin{center}
\begin{tabular}{|c|c|c|c|c|c|} \hline
   $k$                  &    1   &   2    &   3    &    4    &    5    \\    \hline  
$\mu_{2k-1}$            & 2.0061 & 5.1253 & 8.2600 & 11.3982 & 14.5378 \\              
$\mu_{2k}$              & 3.4533 & 6.6286 & 9.7839 & 12.9330 & 16.0794 \\ 
$\pi k$                 & 3.1416 & 6.2832 & 9.4248 & 12.5664 & 15.7080 \\ 
$[\Psi_{2k}(\infty)]^2$ & 0.0664 & 0.0391 & 0.0279 & 0.0218  &  0.0178 \\  
$1/(\pi^2 k)$           & 0.1013 & 0.0507 & 0.0338 & 0.0253  &  0.0203 \\  \hline \hline
   $k$                  &    6    &   7     &     8   &     9   &    10   \\  \hline
$\mu_{2k-1}$            & 17.6780 &  20.8187& 23.9596 & 27.1006 & 30.2418 \\          
$\mu_{2k}$              & 19.2242 & 22.3682 & 25.5116 & 28.6547 & 31.7974 \\
$\pi k$                 & 18.8496 & 21.9911 & 25.1327 & 28.2743 & 31.4159 \\
$[\Psi_{2k}(\infty)]^2$ & 0.0151  &  0.0131 & 0.0116  & 0.0104  &  0.0094 \\  
$1/(\pi^2 k)$           & 0.0169  &  0.0145 & 0.0127  & 0.0113  &  0.0101 \\  \hline
\end{tabular}
\end{center}
\caption{
List of eigenvalues $\mu_{2k}$ and coefficients
$[\Psi_{2k}(\infty)]^2$ of the first 10 contributing terms in the
spectral expansion (\ref{eq:Cmu_spectral_bis}) for $\Cmu(\mu)$ 
(in addition, one has $\mu_0 = 0$ and $\Psi_0(\infty) = 1/\sqrt{2}$).
The reported values were obtained numerically by using a matrix
representation of the Steklov problem in elliptic coordinates (see
Appendix \ref{sec:Steklov_interval}).  The matrix was truncated to the
size $100 \times 100$ and then diagonalized numerically.  The shown
values did not change when the truncation order was increased to $500
\times 500$.  For comparison, the large-$k$ asymptotic approximations
of $\mu_{2k}$ and $[\Psi_{2k}(\infty)]^2$ from
Eqs. (\ref{eq:muk_asympt}, \ref{eq:Vkinf_asympt}) are also present.
For completeness, we also present the first ten eigenvalues
$\mu_{2k-1}$ that correspond to antisymmetric eigenfunctions
$\Psi_{2k-1}$ that vanish at infinity.}
\label{tab:mu_int}
\end{table}

\section{Steklov eigenmodes}
\label{sec:Steklov_interval}

In this Appendix, we recall a numerical computation of the Steklov
eigenfunctions $\Psi_k$ satisfying Eqs. (\ref{eq:Vkint}).  The details
of this computation are provided in Appendix D of
Ref. \cite{Grebenkov25}.  Since $\Psi_0 = 1/\sqrt{2}$ is known, we
focus on the other eigenfunctions with $k = 1,2,\ldots$.

In elliptic coordinates $(\alpha,\theta)$, one has
\begin{equation}
y_1 = \cosh \alpha \cos \theta\,, \qquad y_2 = \sinh\alpha \sin \theta \,,
\end{equation}
with $0 \leq \alpha < +\infty$ and $0 \leq \theta \leq \pi$.  Note
that Eqs. (\ref{eq:rpm}, \ref{eq:elliptic_coord}) with $a_E = 1$ allow
one to express $\alpha$ and $\theta$ in terms of $y_1$ and $y_2$.

The Steklov eigenfunctions can be written as
\begin{equation}  \label{eq:Psik_ck}
  \Psi_k(\alpha,\theta) = \sum\limits_{n=0}^\infty c_{k,n}
  \cos(n\theta) e^{-n \alpha} \,,
\end{equation}
with unknown coefficients $c_{k,n}$.  Imposing the Steklov condition
yields the infinite system of linear equations:
\begin{equation}  \label{eq:eigen_ellipse}
  \sum\limits_{n=1}^\infty c_{k,n} \MM_{n,m} = \frac{1}{\mu_k} c_{k,m} \,
  \quad \mbox{for}\,\, m=1,2,\ldots\,,
\end{equation}
where 
\begin{equation}  \label{eq:M_ellipse}
  \MM_{n,m} = \frac{1}{m} \biggl[\AA_{n,m} -
  \frac{\AA_{n,0} \AA_{0,m}}{\AA_{0,0}}\biggr]\,,
\end{equation} 
and
\begin{equation} \nonumber
  \AA_{n,m}  = \frac{1+(-1)^{m+n}}{\pi}
  \biggl(\frac{1}{1-(m-n)^2} + \frac{1}{1-(m+n)^2}\biggr)\,.
\end{equation}
This matrix equation determines the coefficients $c_{k,n}$ up to a
multiplicative factor that has to be fixed by the normalization
(\ref{eq:Vk_orthonormal}) of Steklov eigenfunctions.  Using the
following relation derived in \cite{Grebenkov25} for $k = 1,2,\ldots$
\begin{equation}
  1 = \int\limits_{-1}^1 |\Psi_k(y_1,0)|^2 \, dy_1 =
  \frac{\pi}{2\mu_{k}} \sum\limits_{m=1}^{\infty} m |c_{k,m}|^2 \,,
\end{equation}
one can ensure the required normalization to the coefficients
$c_{k,m}$.  Once the coefficients $c_{k,n}$ with $n=1,2,\ldots$ are
found, one also gets
\begin{equation}
c_{k,0} = - \frac{1}{\AA_{0,0}} \sum\limits_{n=1}^\infty c_{k,n} \AA_{n,0}\,.
\end{equation}
Note that $c_{k,0}$ is actually the value of $\Psi_k$ at infinity.  

\begin{figure}
\begin{center}
\includegraphics[width=88mm]{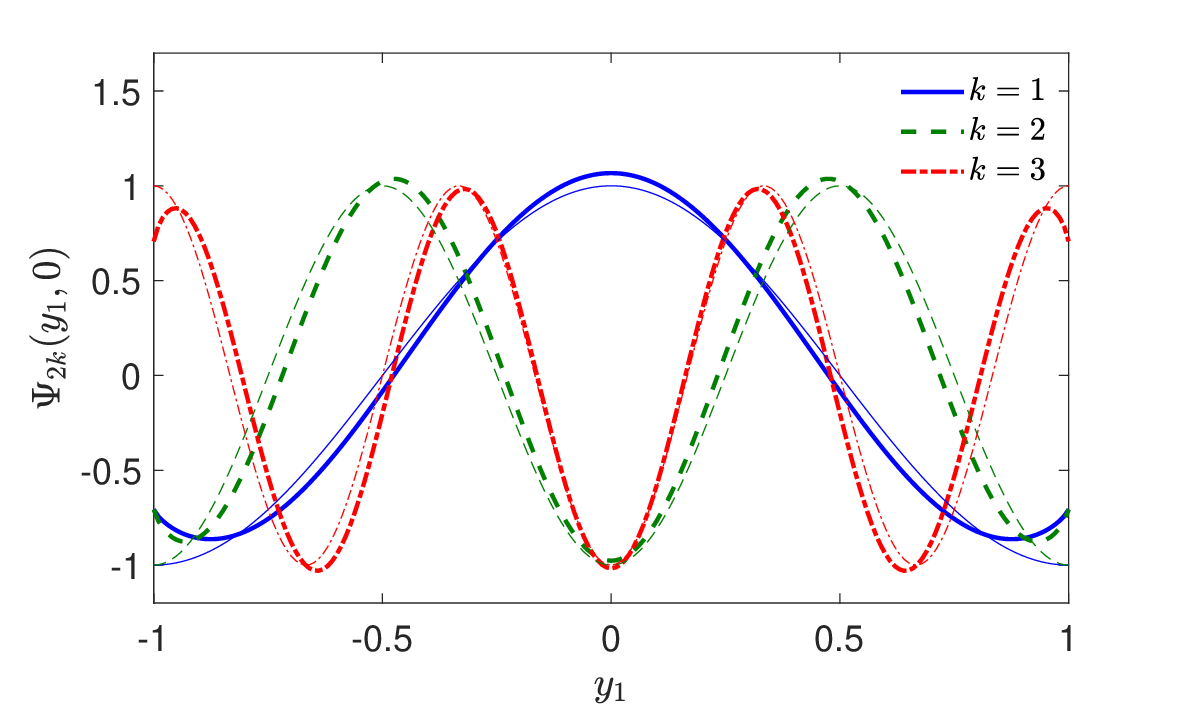} 
\end{center}
\caption{
The Steklov eigenfunctions $\Psi_{2k}(y_1,0)$, restricted onto the
interval $(-1,1)$, are shown by thick lines.  These
eigenfunctions are obtained by truncating the series in
Eq. (\ref{eq:Psik_ck}) to $n \leq 100$, with the coefficients
$c_{k,n}$ found by diagonalizing the truncated matrix $\MM$ from
Eq. (\ref{eq:M_ellipse}).  For comparison, functions $\cos(\pi k
y_1)$ are plotted by thin lines.}
\label{fig:Vkint}
\end{figure}

In practice, one can truncate the infinite-dimensional matrix $\MM$ to
a finite size $M\times M$ and then diagonalize it numerically.  Its
eigenvalues and eigenvectors approximate $1/\mu_k$ and $c_{k,n}$,
respectively.  We checked numerically that these approximations
converge very rapidly as $M$ increases.  We used this technique to
obtain the numerical values reported in Table \ref{tab:mu_int}.
Figure \ref{fig:Vkint} shows several eigenfunctions $\Psi_{2k}(y_1,0)$
and their approximations by $\cos(\pi k y_1)$.

\section{Asymptotic behavior of $\Cmu(\mu)$}
\label{sec:Cmu}

In this Appendix, we derive the exact values of the coefficients $C_1$
and $C_2$ of the Taylor expansion (\ref{eq:Cmu_Taylor}) of $\Cmu(\mu) -
\pi/(2\mu)$ as $\mu\to 0$.  

From the divergence theorem, there is no solution to
Eq. (\ref{eq:Steklov_inner}) for $\mu = 0$.  As such, for $\mu \ll 1$,
the solution should bifurcate from infinity, so that $\Cmu(\mu)$ is
expected to be large in this limit.  For this reason, we expand the
solution for $0 < \mu \ll 1$ as
\begin{equation}
g_\mu(\y) = - \frac{C_0}{\mu} + v_1(\y) - \mu v_2(\y) + \mu^2 v_3(\y) + \ldots
\end{equation}
Inserting this expansion into Eq. (\ref{eq:Steklov_inner}), we obtain
three BVPs:
\begin{subequations}  \label{eq:v1}
\begin{align}
\Delta v_1 & = 0  \quad \textrm{in}~ \H_2\,,   \\
\partial_n v_1 & = C_0 \quad \textrm{on}~ y_2 =0\,, ~ |y_1| < 1 \,; \qquad
\partial_n v_1  = 0  \quad \textrm{on}~ y_2 = 0\,, ~ |y_1| \geq 1\,,\\
v_1 & \sim \ln |\y| + {\mathcal O}(1) \quad \textrm{as}~ |\y| \to \infty\,,
\end{align}
\end{subequations}
for $v_1(\y)$,
\begin{subequations}  \label{eq:v2}
\begin{align}
\Delta v_2 & = 0  \quad \textrm{in}~ \H_2\,,   \\
\partial_n v_2 & = v_1 \quad \textrm{on}~ y_2 =0\,, ~ |y_1| < 1 \,; \qquad
\partial_n v_2  = 0  \quad \textrm{on}~ y_2 = 0\,, ~ |y_1| \geq 1\,,\\
v_2 & \sim {\mathcal O}(1) \quad \textrm{as}~ |\y| \to \infty\,,
\end{align}
\end{subequations}
for $v_2(\y)$, and
\begin{subequations}  \label{eq:v3}
\begin{align}
\Delta v_3 & = 0  \quad \textrm{in}~ \H_2\,,   \\
\partial_n v_3 & = v_2 \quad \textrm{on}~ y_2 =0\,, ~ |y_1| < 1 \,; \qquad
\partial_n v_3  = 0  \quad \textrm{on}~ y_2 = 0\,, ~ |y_1| \geq 1\,,\\
v_3 & \sim {\mathcal O}(1) \quad \textrm{as}~ |\y| \to \infty\,,
\end{align}
\end{subequations}
for $v_3(\y)$.  Note that the solutions $v_1$, $v_2$ and $v_3$ are
known only up to additive constants.  These constants are needed for
ensuring that Eqs. (\ref{eq:v1}, \ref{eq:v2}, \ref{eq:v3}) have
solutions.
In particular, by the divergence theorem, there exists a solution to
(\ref{eq:v1}) provided that 
\begin{equation}
  0 = \int\limits_{-1}^1 \partial_n v_1 \vert_{y_2=0} \, dy_1 +
  \lim\limits_{R\to \infty} \int\limits_{B_R \cap \H_2} \partial_n v_1 \, ds=
  2C_0 + \pi,
\end{equation}
where $B_R$ is a large disk of radius $R$.  This yields
\begin{equation} \label{eq:C0}
C_0 = - \frac{\pi}{2} \,.
\end{equation}

To get the next-order terms, we will use the following lemma, which
follows by using the method of images.

{\bf Lemma}:  {\it For a given function $f(\xi)$ such that
 $\left| \int\limits_{-\infty}^\infty f(\xi) d\xi \right| < \infty$,
the solution $w(\y)$ to the boundary value problem
\begin{equation}
  \Delta w = 0 \quad \textrm{in}~ \H_2\,, \qquad \partial_n w = f(y_1)
  \quad \textrm{on}~ y_2 =0 \,,
\end{equation}
is
\begin{equation}
  w(\y) = - \frac{1}{2\pi} \int\limits_{-\infty}^\infty
  \ln [(\xi - y_1)^2 + y_2^2] f(\xi)\, d\xi\,.
\end{equation}
In particular, one has
\begin{equation}
  w(\y) \sim \left(- \frac{1}{\pi} \int\limits_{-\infty}^\infty f(\xi) \,
    d\xi \right) \ln |\y| + o(1)  \qquad \mbox{as} \,\,\, |\y|\to \infty\,.
\end{equation}
}

To apply this lemma, we decompose the solution $v_1$ as $v_1(\y) =
u_1(\y) + C_1$, where $u_1(\y)$ is a solution for which $u_1(\y)
\sim\ln|\y| + o(1)$ as $|\y|\to\infty$.
In this way, using Eq. (\ref{eq:C0}) and setting 
\begin{equation}
f(\xi) = - (\pi/2) \Theta(1-|\xi|) \,,
\end{equation}
where $\Theta(z)$ is the Heaviside step function ($\Theta(z) = 1$ for
$z >0$ and $0$ otherwise), we get
\begin{equation}\label{app:U1}
u_1(\y) = \frac14 \int\limits_{-1}^1 \ln[(\xi - y_1)^2 + y_2^2] \, d\xi \,,
\end{equation}
which satisfies the required condition $u_1(\y) \sim \ln |\y| + o(1)$
as $|\y|\to \infty$.  

Now we turn to the next-order term $v_2$ satisfying Eq. (\ref{eq:v2}),
in which $v_1|_{y_2=0} = C_1 + u_1|_{y_2=0}$.  The divergence theorem
yields that the necessary and sufficient condition for $v_2$ to be bounded
(i.e., $v_2 \sim {\mathcal O}(1)$ as $|\y|\to \infty$)is
\begin{equation}
\int\limits_{-1}^1 v_1(y_1,0) \, dy_1 = 0 \,,
\end{equation}
from which
\begin{equation} \label{eq:C1_0}
C_1 = - \frac12 \int\limits_{-1}^1 u_1(y_1,0) \, dy_1 \,,
\end{equation}
where from an integration of Eq. (\ref{app:U1}) with $y_2=0$ we get
\begin{equation} \label{eq:u1_int2} 
u_1(y_1,0) = \frac12 \biggl[-2 +
  y_1 \ln\biggl(\frac{1+y_1}{1-y_1}\biggr) + \ln(1 - y_1^2)\biggr]\,.
\end{equation}
Substituting this expression into Eq. (\ref{eq:C1_0}), we get
\begin{equation} \label{eq:C1}
C_1 = \frac32 - \ln 2\,.
\end{equation}

Similarly, we write $v_2 = u_2 + C_2$ such that $u_2(\y) \sim o(1)$ as
$|\y|\to \infty$.  Using again the above lemma with $f(\xi) =
u_1(\xi,0) + C_1$, we conclude that
\begin{equation}
  u_2(\y) = - \frac{1}{2\pi} \int\limits_{-1}^1 \ln[(\xi-y_1)^2+y_2^2]
  (u_1(\xi,0) + C_1) \, d\xi \,.
\end{equation}
In particular, we find
\begin{equation*} 
  u_2(y_1,0)  = -\frac{1}{\pi} \int\limits_{-1}^1 \ln |\xi-y_1|
  \left(u_1(\xi,0) + C_1\right) \, d\xi = - \frac{2C_1}{\pi} u_1(y_1,0)
  - \frac{1}{\pi}
  \int\limits_{-1}^1 \ln|\xi-y_1| u_1(\xi,0)\,  d\xi \,,
\end{equation*}
where we used Eq. (\ref{app:U1}) with $y_2=0$ for the first term.
Rewriting Eq. (\ref{eq:u1_int2}) as
\begin{equation}
u_1(y_1,0) = - 1 + \frac12 w(y_1)
\end{equation}
with
\begin{equation}  \label{eq:w_u1}
w(y_1) = y_1 \ln\biggl(\frac{1+y_1}{1-y_1}\biggr) + \ln(1-y_1^2),
\end{equation}
we get 
\begin{align} 
  u_2(y_1,0) & = \frac{2(1-C_1)}{\pi} u_1(y_1,0) - \frac{1}{2\pi}
               \int\limits_{-1}^1 \ln|\xi-y_1| w(\xi) \, d\xi \,.
\end{align}

The divergence theorem applied to Eq. (\ref{eq:v3}) yields
\begin{equation}
\int\limits_{-1}^1 v_2(y_1,0) \, dy_1 = 0 , 
\end{equation}
from which
\begin{equation*}
C_2  = -\frac12 \int\limits_{-1}^1 u_2(y_1,0) \, dy_1 
 = \frac{-(1-C_1)}{\pi} \underbrace{\int\limits_{-1}^1 u_1(y_1,0) \, dy_1}_{=-2C_1}
 + \frac{1}{4\pi} \int\limits_{-1}^1 \biggl\{ \int\limits_{-1}^1 \ln |\xi-y_1|
 w(\xi) \, d\xi \biggr\} \, dy_1\,.
\end{equation*}
Exchanging the order of integrals in the second term and using
(\ref{app:U1}) with $y_2=0$, we get
\begin{align*}
C_2 & = \frac{2(1-C_1)C_1}{\pi} + \frac{1}{2\pi} \int\limits_{-1}^1 w(\xi) u_1(\xi,0) \, d\xi =
 \frac{2(1-C_1)C_1}{\pi} + \frac{1}{2\pi} \int\limits_{-1}^1 w(\xi) \biggl[-1 + \frac12 w(\xi)\biggr] \, d\xi \\
& = \frac{2(1-C_1)C_1}{\pi} - \frac{2(1-C_1)}{\pi}  + \frac{1}{4\pi}\int\limits_{-1}^1 \left[w(\xi)\right]^2 \, d\xi =
 -\frac{2(1-C_1)^2}{\pi} + \frac{1}{4\pi}\int\limits_{-1}^1 \left[w(\xi)\right]^2 \, d\xi .
\end{align*}
Substituting $w(\xi)$ from Eq. (\ref{eq:w_u1}) and $C_1$ from
Eq. (\ref{eq:C1}), we get after some simplifications that
\begin{equation}  \label{eq:C2}
C_2 = \frac{21 - 2\pi^2}{18 \pi} \,.
\end{equation}

\section{Neumann Green's functions}
\label{sec:NGreen}

In this Appendix, we summarize the available results on various
Neumann Green's for both interior and exterior settings when the
singularity is either on the boundary and in the bulk.  Table
\ref{tab:NGreen} collects their definitions and formulas for three
shapes: the unit disk, ellipses, and rectangles (some derivations are
provided below).  Even though the exterior (bulk) Neumann Green's
function $G_{\rm eb} (\x,\xxi)$ is not discussed in the main text, we
provide its definition for completeness.  For a fixed point $\xxi \in
\Omega_0 = \R^2 \backslash C$ in the exterior of a compact set $C$,
$G_{\rm eb} (\x,\xxi)$ satisfies
\begin{subequations} \label{eq:Neumann-Green_extB}
\begin{align}
\Delta_{\x} G_{\rm eb} & = -\delta(\x-\xxi)  \quad \textrm{in}~ \Omega_0\,,
   \qquad\qquad \partial_n G_{\rm eb} = 0 \quad \textrm{on}~ \pa_0  \,, \\
  \label{eq:Gext_inf}
  G_{\rm eb} (\x,\xxi) & \sim -\frac{1}{2\pi}\ln |\x| + o(1) \quad
                        \textrm{as}~ |\x|\to \infty\,, 
\end{align}
\end{subequations}
whereas the regular part $R_{\rm eb}(\xxi)$ characterizes its singular
behavior near $\xxi$: 
\begin{equation}  \label{eq:Reb_def}
  R_{\rm eb}(\xxi) = \lim\limits_{\x \to \xxi} \bigl(G_{\rm eb}(\x,\xxi) +
  \frac{1}{2\pi} \ln|\x-\xxi|\bigr) \,.
\end{equation}
Note that the $o(1)$ condition in the asymptotic behavior
(\ref{eq:Gext_inf}) determines this function uniquely.  For instance,
for the exterior of the unit disk, we can readily derive by summing an
eigenfunction expansion that
\begin{equation}  \label{eq:disk_NGreen_extB}
  G_{\rm eb} (\x,\xxi) = - \frac{1}{2\pi} \biggl(\ln |\x - \xxi| +
  \ln \bigl|\x - \xxi/|\xxi|^2\bigr| - \ln |\x|\biggr)\,,
\end{equation}
from which $R_{\rm eb}(\xxi) = - \ln (1 - 1/|\xxi|^2)/(2\pi)$.

\begin{table}
\begin{center}
\begin{tabular}{|c|c|c|c|c|}  \cline{1-5}
          & \multicolumn{2}{c|}{Interior problem} & \multicolumn{2}{c|}{Exterior problem}  \\  \cline{1-5}  
          & surface                & bulk  & surface & bulk \\  \cline{1-5}
Definition& (\ref{eq:surfNeumann}) & (\ref{eq:Neumann-Green}) & (\ref{eq:Neumann-Green_ext}) & (\ref{eq:Neumann-Green_extB})     \\  \cline{1-5} 
Disk      & (\ref{eq:Green_disk}) & (\ref{eq:NGreen_disk})   & (\ref{eq:Green_disk_ext1}), see Sec. \ref{sec:Gext_disk}   &  
(\ref{eq:disk_NGreen_extB}) \\   
Ellipse   & (\ref{eq:ellip_NGreen_surf}, \ref{eq:ellip_NGreenR_surf})
& (\ref{eq:ellipse_NGreen_intB}) & (\ref{eq:ellip_NGreen_extS}, \ref{eq:ellip_NGreenR_extS}) &  
(\ref{eq:Gext_ellipse}, \ref{eq:Rext_ellipse})    \\   
Rectangle &  accessible       & Eq. (4.13) from \cite{kolok_split}, and \cite{McCann01}      &  unavailable  & unavailable     \\  \cline{1-5} 
\end{tabular}
\end{center}
\vspace{0.2cm}
\caption{
Summary of available formulas for various Neumann Green's functions.
Note that the surface Neumann Green's function for rectangles can be
derived from the results in \cite{kolok_split,McCann01}. In turn, its
extension to the exterior problem is not available.}
\label{tab:NGreen}
\end{table}

\subsection{Surface Neumann Green's function for the exterior of the unit disk}
\label{sec:Gext_disk}

In this Appendix, we provide a rigorous derivation of
Eq.~(\ref{eq:Green_disk_ext1}) for the surface Neumann Green's
function $G_{\rm e}(\x,\xxi)$ for the exterior of the unit disk.

To establish the result in Eq.~(\ref{eq:Green_disk_ext1}), we let
$r=|\x|$ and we decompose $G_{\rm e}$ as
$G_{\rm e}(\x,\xxi)=-(2\pi)^{-1}\ln|\x|+H_{\rm e}(\x,\xxi)$, to obtain that
$H_{\rm e}(\x,\xxi)$ satisfies
\begin{subequations}
\begin{align}
  \Delta_{\x} H_{\rm e}= 0  \quad \mbox{in}\quad |\x|>1 \,; \qquad
  \partial_r H_{\rm e} =\frac{1}{2\pi}  \quad \mbox{on} \,\,\, |\x|=1\,,
  \,\, \x\neq \xxi\,,\\
  H_{\rm e}\sim -\frac{1}{\pi}\ln|\x-\xxi| \quad \mbox{as}\quad
  \x\to\xxi\,; \qquad H_{\rm e}\to 0 \quad \mbox{as}\quad |\x|\to \infty\,,
\end{align}
\end{subequations}
where $\Delta_{\x}$ is the Laplacian in the $\x$-variable.
Consider now the interior surface Neumann Green's function inside the
disk $|\y|\leq 1$ satisfying Eq. (\ref{eq:surfNeumann}), which we
decompose as $G(\y;\xxi)={|\y|^2/(4\pi)}+ H_{\rm i}(\y,\xxi)$, where
$H_{\rm i}(\y,\xxi)$ satisfies
\begin{subequations}
\begin{align}
  \Delta_{\y} H_{\rm i}= 0  \quad \mbox{in}\quad |\y|<1 \,; \qquad
  \partial_\rho H_{\rm i} =\frac{-1}{2\pi} \quad \mbox{on} \quad |\y|=1\,,
  \,\,\, \y\neq \xxi\,,\\
  H_{\rm i}\sim -\frac{1}{\pi}\ln|\y-\xxi|  \quad \mbox{as}\quad
  \y\to\xxi\,; \qquad H_{\rm i} \quad \mbox{bounded as}\quad \y \to
  {\boldsymbol 0}\,.
\end{align}
\end{subequations}
Here $\rho=|\y|$ and $\Delta_{\y}$ denotes the Laplacian in the $\y$
variable.  From Eq. (\ref{eq:Green_disk}) it follows that
\begin{equation}
  H_{\rm i}(\y,\xxi) = -\frac{1}{\pi} \ln|\y-\xxi| + C\,,
\end{equation}
where $C$ is a constant.  By using conformal invariance under Kelvin's
transformation $\y={\x/|\x|^2}$ in the unit disk, and noting that
$\partial_r=-\partial_\rho$ on $r=\rho=1$, it follows that $H_{\rm e}$
is given by
\begin{equation}\label{eq:he_ex}
  H_{\rm e}(\x,\xxi)= H_{\rm i}\left(\frac{\x}{|\x|^2},\xxi\right)= -
    \frac{1}{\pi} \ln {\biggl \vert} \frac{\x}{|\x|^2}-\xxi {\biggr\vert } \,,
\end{equation}
where we observe that $H_{\rm e}\to 0$ as $|\x|\to \infty$.  Finally, if
$|\xxi|=1$, we can readily calculate that
\begin{equation}\label{eq:he_ident}
  {\biggl\vert} \frac{\x}{|\x|^2}-\xxi{\biggr \vert} =
  \frac{|\x-\xxi|}{|\x|} \,.
\end{equation}
Upon subsituting Eq. (\ref{eq:he_ident}) into Eq. (\ref{eq:he_ex}) and
using $G_{\rm e}=-(2\pi)^{-1}\ln|\x|+H_{\rm e}$, we obtain
the result in Eq. (\ref{eq:Green_disk_ext1}).


\subsection{Neumann Green's functions for an ellipse}
\label{sec:NGreen_ellipse_int}

For an ellipse with semiaxes $a$ and $b$ ($a > b$),
$\Omega_0 = \{ \x = (x_1,x_2) \in \R^2 ~:~ (x_1/a)^2 + (x_2/b)^2 <
1\}$, a rapidly converging representation for the ``bulk'' Neumann
Green's function was derived in \cite{Iyaniwura21}.  In the elliptic
coordinates introduced in Eq. (\ref{eq:ellip_coord}), Eq. (5.21a) from
\cite{Iyaniwura21} reads
\begin{equation}  \label{eq:ellipse_NGreen_intB}
  G_b(\x,\xxi) = \frac{|\x|^2 + |\xxi|^2}{4\pi ab} -
  \frac{3(a^2+b^2)}{16\pi ab} + \frac{\alpha_b - \alpha_>}{2\pi} + S(\x,\xxi)\,,
\end{equation}
where $\x = (\alpha,\theta)$, $\xxi = (\alpha_0, \theta_0)$, $\alpha_>
= \max\{ \alpha,\alpha_0\}$, $\beta = (a-b)/(a+b)$, $\alpha_b =
\atanh(b/a) = \tfrac12 \ln(1/\beta)$ describes the boundary $\pa_0$,
\begin{equation}
S(\x,\xxi) = - \frac{1}{2\pi} \sum\limits_{n=0}^\infty
\sum\limits_{j=1}^8 \ln|1 - \beta^{2n} z_j| \,,
\end{equation}
and
\begin{align*}
z_1 & = e^{-|\alpha-\alpha_0| + i(\theta-\theta_0)} , \qquad
z_2  = e^{-4\alpha_b + |\alpha-\alpha_0| + i(\theta-\theta_0)} , \qquad
z_3  = e^{-2\alpha_b-\alpha-\alpha_0 + i(\theta-\theta_0)}, \qquad
z_4  = e^{-2\alpha_b+\alpha+\alpha_0 + i(\theta-\theta_0)}, \\
z_5 & = e^{-4\alpha_b+\alpha+\alpha_0 + i(\theta+\theta_0)}, \qquad
z_6  = e^{-\alpha-\alpha_0 + i(\theta+\theta_0)} ,\qquad
z_7  = e^{-2\alpha_b + |\alpha-\alpha_0| + i(\theta+\theta_0)}, \qquad
z_8  = e^{-2\alpha_b - |\alpha-\alpha_0| + i(\theta+\theta_0)} \,.
\end{align*}

Setting $\xxi$ to the boundary $\pa_0$, we get some simplifications.
In particular, we have $\alpha \leq \alpha_0 = \alpha_b$ that implies
\begin{align*}
z_1 & = z_4 = e^{-\alpha_b+\alpha + i(\theta-\theta_0)}  \,,\qquad
z_3 = z_2 = e^{-3\alpha_b - \alpha + i(\theta-\theta_0)} \,,\\
z_5 & = z_8 = e^{-3\alpha_b+\alpha + i(\theta+\theta_0)} \,,\qquad
z_7 = z_6 = e^{-\alpha_b-\alpha + i(\theta+\theta_0)} \,.
\end{align*}
As a consequence, we obtain a rapidly converging representation for
the surface Neumann Green's function:
\begin{align} \label{eq:ellip_NGreen_surf}
G(\x,\xxi) & = \frac{|\x|^2 + |\xxi|^2}{4\pi ab} - \frac{3(a^2+b^2)}{16\pi ab}  
             - \frac{1}{\pi} \sum\limits_{n=0}^\infty
             \sum\limits_{j=1}^4 \ln|1 - \beta^{2n} z_{2j-1}| \,.
\end{align}
The regular part of this function can be deduced as $\x\to\xxi$.
Setting $\theta = \theta_0$ and $\alpha = \alpha_0 -
\epsilon$, we find as $\epsilon\to 0$:
\begin{align} \nonumber
  & G(\x,\xxi) \approx \frac{|\x|^2 + |\xxi|^2}{4\pi ab} - \frac{3(a^2+b^2)}
    {16\pi ab} 
 -\frac{1}{\pi} \ln(1 - e^{-\epsilon}) 
 - \frac{1}{\pi} \sum\limits_{n=1}^\infty \ln(1 - \beta^{2n}) 
    - \frac{1}{\pi} \sum\limits_{n=0}^\infty \sum\limits_{j=2}^4 \ln|1 -
    \beta^{2n} z_{2j-1}|  \,,
\end{align}
where 
\begin{align*}
z_3 & = e^{-4\alpha_b} = \beta^2\,,\qquad
z_5 = z_7 = e^{-2\alpha_b + 2i\theta_0} = \beta e^{2i\theta_0} \,.
\end{align*}
On the other hand, we have for $\theta = \theta_0$ that
\begin{align*}  
  |\x-\xxi|^2 & = a_E^2\bigl[(\cosh\alpha \cos\theta -
                \cosh\alpha_0 \cos\theta_0)^2  
+ (\sinh\alpha \sin\theta - \sinh\alpha_0 \sin\theta_0)^2\bigr] \\    
              &  \approx a_E^2 \epsilon^2 (\sinh^2 \alpha_0 + \sin^2\theta_0)
                \qquad \mbox{as}
                \quad \epsilon\to 0 \,,
\end{align*}
so that
\begin{equation}  \label{eq:dist}
  |\x - \xxi| \approx \epsilon a_E \sqrt{\sinh^2 \alpha_0 + \sin^2\theta_0}
  \qquad \mbox{as} \quad \epsilon\to 0 \,.
\end{equation}
We can thus express $\epsilon$ in terms of $|\x-\xxi|$ to obtain the
following infinite series representation for the regular part in terms
of the aspect ratio $\beta={(a-b)/(a+b)}$:
\begin{align}  \label{eq:ellip_NGreenR_surf}
R(\xxi) & = \frac{|\xxi|^2}{2\pi ab} - \frac{3(a^2+b^2)}{16 \pi ab} 
+ \frac{1}{\pi} \ln \biggl(a_E \sqrt{\sinh^2\alpha_0 + \sin^2\theta_0}\biggr) 
 - \frac{2}{\pi} \sum\limits_{n=1}^\infty \biggl(\ln(1 - \beta^{2n})
+ \ln\bigl|1 - \beta^{2n-1} e^{2i\theta_0}\bigr|\biggr) \,.
\end{align}

\subsection{Neumann Green's functions for the exterior of an ellipse}
\label{sec:ellip_NGreen_ext}

Finally, we consider the exterior of an ellipse with semiaxes $a > b$:
$\Omega_0 = \{ (x_1,x_2) \in \R^2 ~:~ (x_1/a)^2 + (x_2/b)^2 > 1\}$.
We first derive the ``bulk'' Neumann Green's function
$G_{\rm eb}(\x,\xxi)$ for this domain and then let the singularity
point $\xxi$ tend to the boundary to get $G_{\rm e}(\x,\xxi)$ and its
regular part $R_{\rm e}(\xxi)$.

\subsubsection*{Bulk Neumann Green's function}

In elliptic coordinates introduced in Eq. (\ref{eq:ellip_coord}), we
search the bulk Neumann Green's function as the unique solution of
Eqs.  (\ref{eq:Neumann-Green_extB}) in the form
\begin{equation}
  G_{\rm eb}(\x,\xxi) = \sum\limits_{k=-\infty}^\infty
  A_k(\alpha) e^{ik(\theta-\theta_0)}\,,
\end{equation}
with unknown functions $A_k(\alpha)$.  Substitution of this form into
the equation for the Green's function yields
\begin{align*}
  -\Delta_{\x} G(\x,\xxi) &= - \frac{1}{h_\alpha^2}
                            (\partial_\alpha^2 + \partial_\theta^2) 
\sum\limits_{k=-\infty}^\infty A_k(\alpha) e^{ik(\theta-\theta_0)} 
 = \frac{1}{h_\alpha^2} \delta(\alpha-\alpha_0) \delta(\theta-\theta_0)\,,
\end{align*}
where $h_\alpha = a_E \sqrt{\cosh^2\alpha - \cos^2\theta}$ is the
scale factor, and $a_E = \sqrt{a^2-b^2}$.  Multiplying by
$e^{ik^{\prime}\theta}$ and integrating over $\theta$ from $-\pi$ to $\pi$, we
obtain a set of equations for $A_k(\alpha)$:
\begin{equation}
(-\partial_\alpha^2 + k^2) A_k(\alpha) = \frac{1}{2\pi} \delta(\alpha-\alpha_0).
\end{equation}
Note that $A_k(\alpha)$ should satisfy $A^{\prime}_k(\alpha_b) = 0$, where
$\alpha_b$ is the location of the elliptic boundary (i.e., $\tanh
\alpha_b = b/a$).  As a consequence, we can search for solutions
separately on $\alpha_b < \alpha < \alpha_0$ and $\alpha > \alpha_0$:
\begin{equation}
  A_k(\alpha) = \begin{cases} a_k \cosh(k(\alpha-\alpha_b))\,, \quad
    \alpha_b < \alpha < \alpha_0\,, \cr
b_k e^{-k\alpha} \,,\hskip 30mm \alpha > \alpha_0 \,,\end{cases}
\end{equation} 
where we assumed that $k > 0$.  The unknown coefficients $a_k$ and
$b_k$ are obtained by requiring the continuity of $A_k(\alpha)$ together
with the jump condition for the derivative at $\alpha = \alpha_0$. This
yields that
\begin{align}
b_k e^{-k\alpha_0} &= a_k \cosh(k(\alpha_0 - \alpha_b)) ,\qquad
-k b_k e^{-k\alpha_0} - k a_k \sinh(k(\alpha_0-\alpha_b)) = -1/(2\pi) \,,
\end{align}
from which we determine
%
\begin{align}
a_k & = \frac{1}{2\pi k} e^{-k(\alpha_0-\alpha_b)} \,,  \qquad
b_k = \frac{1}{2\pi k} e^{k\alpha_b} \cosh(k(\alpha_0-\alpha_b)) \,.
\end{align}
We conclude that
\begin{equation}
A_k = \frac{1}{2\pi k} e^{-k (\alpha_> - \alpha_b)} \cosh(k(\alpha_< - \alpha_b)) \,,
\end{equation}
where $\alpha_< = \min\{\alpha,\alpha_0\}$ and
$\alpha_> = \max\{\alpha,\alpha_0\}$. We can further simplify this expression
to
\begin{equation}
  A_k = \frac{1}{4\pi k} \left[e^{-k |\alpha-\alpha_0|} + e^{-k(\alpha + \alpha_0  -
      2\alpha_b)}\right] \,.
\end{equation}
We recall that this solution holds for $k > 0$.  By symmetry, it also
holds for $k < 0$, if $k$ is replaced by $|k|$.

For $k = 0$, a general solution of the Laplace equation on the
interval $\alpha_b < \alpha < \alpha_0$ is $A_0 = a_0 + c_0 \alpha$,
where we must set $c_0 = 0$ to ensure that the Neumann condition at
$\alpha_b$ is satisfied.  In turn, we have $A_0 = d_0 + b_0 \alpha$
for $\alpha > \alpha_0$.  We set $a_0 = d_0 + b_0\alpha_0$ to ensure
the continuity of $A_0$ at $\alpha_0$.  The coefficient $b_0$ is
determined by the jump of the derivative, which yields
$b_0 = - 1/(2\pi)$.  We conclude that, in terms of a constant
$d_0$ to be fixed, $A_0$ has the form
\begin{equation*}
A_0 = d_0 - \frac{1}{2\pi} \alpha_> \,.
\end{equation*}

Combining these results, we get
\begin{align*}
G(\x,\xxi) & = d_0 - \frac{\alpha_>}{2\pi} + \frac{1}{4\pi} \sum\limits_{k\ne 0} 
 \frac{e^{ik(\theta-\theta_0)} }{|k|}
   \bigl(e^{-|k| |\alpha-\alpha_0|} + e^{-|k|(\alpha + \alpha_0  - 2\alpha_b)}\bigr)\,, \\
           & = d_0 - \frac{\alpha_>}{2\pi} - \frac{1}{4\pi}
             \biggl[\ln(1 - e^{-|\alpha-\alpha_0|+i(\theta-\theta_0)}) 
+ \ln(1 - e^{-|\alpha-\alpha_0|-i(\theta-\theta_0)})  \\
& \qquad + \ln(1 - e^{-(\alpha+\alpha_0-2\alpha_b)+i(\theta-\theta_0)}) 
+ \ln(1 - e^{-(\alpha+\alpha_0-2\alpha_b)-i(\theta-\theta_0)}) \biggr].
\end{align*}
Since
$|\x|^2 = a_E^2 (\cosh^2\alpha - \sin^2\theta) \approx a_E^2
\cosh^2\alpha$ when $|\x| \gg a_E$, we have
$\alpha_> = \alpha \approx \ln(2|\x|/a_E)$ as $|\x|\to \infty$, where
we used that $\cosh^{-1}(z) = \ln(z + \sqrt{z^2-1})$.  Since the
Neumann Green's function behaves at infinity according to
Eq. (\ref{eq:Gext_inf}), the constant term $d_0$ must compensate the
constant contribution from $\alpha_>$. This condition yields
$d_0 = \ln\left({2/a_E}\right)/(2\pi)$, and so we conclude that that
\begin{align}  \label{eq:Gext_ellipse}
G_{\rm eb}(\x,\xxi) & = \frac{\ln(2/a_E) - \alpha_>}{2\pi} 
   - \frac{1}{4\pi} \biggl[\ln\bigl(1 - 2\cos(\theta-\theta_0)
      e^{-|\alpha-\alpha_0|} + e^{-2|\alpha-\alpha_0|}\bigr) \\  \nonumber
   & \qquad + \ln\bigl(1 - 2\cos(\theta-\theta_0) e^{-(\alpha+\alpha_0-2\alpha_b)}
               + e^{-2(\alpha+\alpha_0-2\alpha_b)}\bigr) \biggr]\,.
\end{align}

To evaluate the regular part, we set $\theta = \theta_0$ and $\alpha =
\alpha_0 + \epsilon$, so that
\begin{align*}
  G_{\rm eb}(\x,\xxi) & \approx \frac{1}{2\pi} \left[\ln(2/a_E) -
  \alpha_0 - \ln(1 - e^{-\epsilon}) - \ln\bigl(1 - e^{-2(\alpha_0-\alpha_b)}\bigr)
                        \right] \,.
\end{align*}
Using again Eq. (\ref{eq:dist}) to express $\epsilon$ in terms of
$|\x-\xxi|$, we conclude that
\begin{equation}  \label{eq:Rext_ellipse}
  R_{\rm eb}(\xxi) = \frac{1}{2\pi} \left[\frac12 \ln(\sinh^2 \alpha_0 +
    \sin^2\theta_0) + \ln(2) - \alpha_0 -
    \ln\bigl(1 - e^{-2(\alpha_0-\alpha_b)}\bigr)\right]\,.
\end{equation}  

\subsubsection*{Surface Neumann Green's function}

Setting $\xxi\in \pa$ in Eq. (\ref{eq:Gext_ellipse}), we get a
simplification for the surface Neumann Green's function:
\begin{equation}  \label{eq:ellip_NGreen_extS}
G_{\rm e}(\x,\xxi) = \frac{\ln(2/a_E)-\alpha}{2\pi} 
- \frac{1}{2\pi} \ln\bigl(1 - 2\cos(\theta-\theta_0)e^{-(\alpha-\alpha_b)}
+ e^{-2(\alpha-\alpha_b)}\bigr) \,.
\end{equation}
In turn, setting $\theta = \theta_0$ and
$\alpha = \alpha_b + \epsilon$, and using Eq. (\ref{eq:dist}), we find
that
\begin{equation}  \label{eq:ellip_NGreenR_extS}
  R_{\rm e}(\xxi) = \frac{1}{2\pi} \left[
    \ln(2a_E) - \alpha_0 + \ln(\sinh^2\alpha_0 + \sin^2\theta_0)\right]\,.
\end{equation}  

\bibliographystyle{plain}
\bibliography{references}

\end{document}